\documentclass[a4paper,12pt,amsart,frenchb]{article}

\usepackage{amsmath,amsbsy,amsfonts,amssymb}
\usepackage[french]{babel}
\newcommand{\ass}[2]{\vskip0.3cm\noindent
{\bf {#1}}. { \sl {#2}}\vskip0.3cm\noindent
}

\begin{document}

  \title{ Repr\'esentations  de r\'eduction unipotente pour $SO(2n+1)$, II: endoscopie}
\author{J.-L. Waldspurger}
\date{3 d\'ecembre 2016}
\maketitle

{\bf Introduction}

Cet article est la suite de \cite{W5}, dont on  utilise les d\'efinitions et notations. Ces derni\`eres sont rappel\'ees dans l'index \`a la fin de l'article. Le corps de base $F$ est local, non-archim\'edien et de caract\'eristique nulle. On note $p$ sa caract\'eristique r\'esiduelle. Un  entier $n\geq1$ est fix\'e pour tout l'article. On suppose $p\geq 6n+4$.  On s'int\'eresse aux groupes $G_{iso}$ et $G_{an}$ d\'efinis en \cite{W5} 1.1.   Le groupe $G_{iso}$ est le groupe sp\'ecial orthogonal d'un espace $V_{iso}$ de dimension $2n+1$ sur $F$ muni d'une forme quadratique $Q_{iso}$ et $G_{an}$  est le groupe sp\'ecial orthogonal d'un espace $V_{an}$ de dimension $2n+1$ sur $F$ muni d'une forme quadratique $Q_{an}$. Le groupe $G_{iso}$ est 
d\'eploy\'e et $G_{an}$ en est la forme int\'erieure non d\'eploy\'ee. Pour un indice $\sharp=iso$ ou $an$, on note $Irr_{tunip,\sharp}$ l'ensemble des classes d'isomorphismes de repr\'esentations admissibles irr\'eductibles de $G_{\sharp}(F)$ qui sont temp\'er\'ees et de r\'eduction unipotente, cf. \cite{W5} 1.3 pour la d\'efinition de cette propri\'et\'e. On note $Irr_{tunip}$ la r\'eunion disjointe de $Irr_{tunip,iso}$ et $Irr_{tunip,an}$. Pour une partition symplectique $\lambda$ de $2n$, fixons un homomorphisme alg\'ebrique $\rho_{\lambda}:SL(2;{\mathbb C})\to Sp(2n;{\mathbb C})$ param\'etr\'e par $\lambda$, cf. \cite{W5} 1.3. On note $Z(\lambda)$ le commutant dans $Sp(2n;{\mathbb C})$ de l'image de $\rho_{\lambda}$. Soit $s\in Z(\lambda)$ un \'el\'ement semi-simple dont toutes les valeurs propres sont de module $1$. On note $Z(s,\lambda)$ le commutant de $s$ dans $Z(\lambda)$, ${\bf Z}(\lambda,s)$ son groupe de composantes connexes et ${\bf Z}(\lambda,s)^{\vee}$ le groupe des caract\`eres de ${\bf Z}(\lambda,s)$. On a vu en \cite{W5} 1.3 que l'on pouvait pr\'esenter la param\'etrisation   de Langlands sous la forme suivante: $Irr_{tunip}$ est param\'etr\'e par l'ensemble des classes de conjugaison (en un sens facile \`a pr\'eciser) de triplets $(\lambda,s,\epsilon)$, o\`u $\lambda$ et $s$ sont comme ci-dessus et $\epsilon\in {\bf Z}(s,\lambda)^{\vee}$. Le param\'etrage a \'et\'e obtenu par diff\'erents auteurs: Lusztig, cf.  \cite{L}; Moeglin, cf. \cite{M} th\'eor\`eme 5.2; Arthur, cf. \cite{Ar1} th\'eor\`eme 2.2.1.   Nous utilisons les constructions de Lusztig. Le but  de l'article est de prouver que les repr\'esentations qu'il a construites  v\'erifient 
  les propri\'et\'es  relatives \`a l'endoscopie qui caract\'erisent la param\'etrisation.   Pour un triplet $(\lambda,s,\epsilon)$ comme ci-dessus, on note $\pi(\lambda,s,\epsilon)$ la repr\'esentation qui lui est associ\'ee par Lusztig. 

Soit $(\lambda,s)$ comme ci-dessus et soit $h\in Z(\lambda,s)$ un \'el\'ement  tel que $h^2=1$. Pour $\epsilon\in {\bf Z}(\lambda,s)^{\vee}$ on note $\epsilon(h)$ la valeur de $\epsilon$ en l'image de $h$ dans $Z(\lambda,s)$. On d\'efinit la repr\'esentation virtuelle
$$\Pi(\lambda,s,h)=\sum_{\epsilon\in {\bf Z}(\lambda,s)^{\vee}}\pi(\lambda,s,\epsilon)\epsilon(h).$$
Pour $\sharp=iso$ ou $an$, on note $\Pi_{\sharp}(\lambda,s,h)$ la sous-somme o\`u l'on se restreint aux $\epsilon$ tels que $\pi(\lambda,s,\epsilon)$ soit une repr\'esentation de $G_{\sharp}(F)$.

Soient $n_{1},n_{2}\in {\mathbb N}$ deux entiers tels que $n_{1}+n_{2}=n$. Un tel couple d\'etermine une donn\'ee endoscopique de $G_{iso}$ ou $G_{an}$ dont le groupe endoscopique est $G_{n_{1},iso}\times G_{n_{2},iso}$ (ce sont les groupes similaires \`a $G_{iso}$ quand on remplace $n$ par $n_{1}$ ou $n_{2}$). Soient $(\lambda_{1},s_{1},h_{1})$ et $(\lambda_{2},s_{2},h_{2})$ des triplets similaires au triplet $(\lambda,s,h)$ ci-dessus, relatifs aux entiers $n_{1}$ et $n_{2}$. Supposons $h_{1}=1$ et $h_{2}=1$. Comme on l'a prouv\'e dans \cite{MW}, les caract\`eres de $\Pi_{iso}(\lambda_{1},s_{1},1)$ et $\Pi_{iso}(\lambda_{2},s_{2},1)$ sont des distributions stables et, en identifiant les repr\'esentations \`a leurs caract\`eres, on d\'efinit le transfert \`a $G_{iso}(F)$ ou $G_{an}(F)$ de $\Pi_{iso}(\lambda_{1},s_{1},1)\otimes \Pi_{iso}(\lambda_{2},s_{2},1)$. A l'aide des triplets $(\lambda_{1},s_{1},h_{1})$ et $(\lambda_{2},s_{2},h_{2})$ et du couple $(n_{1},n_{2})$, on construit naturellement un triplet $(\lambda,s,h)$ pour l'entier $n$, cf. \cite{W5} 2.1. En particulier, $h$ a  pour valeurs propres $1$ avec la multiplicit\'e $2n_{1}$  et $-1$ avec la multiplicit\'e $2n_{2}$. 

\ass{Th\'eor\`eme}{Le transfert de $\Pi_{iso}(\lambda_{1},s_{1},1)\otimes \Pi_{iso}(\lambda_{2},s_{2},1)$ \`a $G_{iso}(F)$ est $\Pi_{iso}(\lambda,s,h)$. Le transfert de $\Pi_{iso}(\lambda_{1},s_{1},1)\otimes \Pi_{iso}(\lambda_{2},s_{2},1)$ \`a $G_{an}(F)$ est $-\Pi_{an}(\lambda,s,h)$. }

On a \'enonc\'e ce r\'esultat en \cite{W5} 2.1 et on va maintenant le prouver. 

La preuve reprend \'etroitement celle de \cite{MW}. On se ram\`ene par induction au cas o\`u les repr\'esentations consid\'er\'ees sont elliptiques au sens d'Arthur. D'apr\`es le r\'esultat de \cite{Ar2}, il suffit alors de comparer les valeurs de leurs caract\`eres sur des \'el\'ements elliptiques fortement r\'eguliers de nos diff\'erents groupes. Ce sont les int\'egrales orbitales de leurs pseudo-coefficients. On a d\'ej\`a calcul\'e ces derniers dans \cite{MW}. A l'aide des travaux de Lusztig et de l'involution que l'on a d\'efinie dans \cite{MW} (on a repris sa d\'efinition dans \cite{W5}), on obtient une expression des caract\`eres en termes d'int\'egrales orbitales de fonctions de Green. Ces fonctions ont \'et\'e d\'efinies sous leur forme la plus g\'en\'erale par Lusztig. Elle vivent sur des groupes finis mais on en d\'eduit par un proc\'ed\'e usuel d'inflation des fonctions sur nos groupes $G_{iso}(F)$ et $G_{an}(F)$ (ou leurs groupes endoscopiques). Il reste \`a d\'emontrer un r\'esultat de transfert entre fonctions de Green. Ces fonctions \'etant \`a support topologiquement unipotent, on descend facilement aux alg\`ebres de Lie. On se rappelle que, d'apr\`es Harish-Chandra, la transform\'ee de Fourier d'une int\'egrale orbitale est une distribution localement int\'egrable et peut donc \^etre consid\'er\'ee comme une fonction. Selon une id\'ee qui remonte \`a Springer, on a exprim\'e dans \cite{W3} puis dans \cite{MW} les int\'egrales orbitales des fonctions de Green comme combinaisons lin\'eaires explicites de transform\'ees de Fourier de certaines int\'egrales orbitales relatives \`a des \'el\'ements semi-simples r\'eguliers. Puisque le transfert endoscopique commute \`a des constantes pr\`es aux transformations de Fourier, on est ramen\'e \`a transf\'erer de telles int\'egrales orbitales, ce qui est \`a peu pr\`es tautologique.

Comme on l'a dit, ce sch\'ema de d\'emonstration se trouve d\'ej\`a dans \cite{MW} o\`u on s'\'etait limit\'e \`a prouver des r\'esultats de stabilit\'e. Pour obtenir les \'egalit\'es endoscopiques, il suffit de poursuivre les calculs un peu plus loin. Ce n'est sans doute pas tr\`es passionnant mais les calculs sont suffisamment compliqu\'es pour m\'eriter, \`a ce qu'ils nous semblent, d'\^etre pr\'esent\'es en d\'etail. Ces calculs sont en particulier compliqu\'es par la pr\'esence de constantes terrifiques. A l'\'evidence, les normalisations que nous avons choisies de nos divers objets ne sont pas les bonnes. Il nous a toutefois sembl\'e pr\'ef\'erable de les conserver, cela permet de rendre compr\'ehensibles les r\'ef\'erences \`a nos articles ant\'erieurs.

 \bigskip

\section{Les  r\'esultats}

\bigskip

\subsection{Calcul des caract\`eres sur les \'el\'ements compacts}
Pour tout groupe r\'eductif d\'efini sur $F$ ou ${\mathbb F}_{q}$ (qui est le corps r\'esiduel de $F$), not\'e $H$, on note par la lettre gothique $\mathfrak{h}$ son alg\`ebre de Lie. Fixons un caract\`ere continu $\psi_{F}:F\to {\mathbb C}^{\times}$ de conducteur l'id\'eal maximal $\varpi \mathfrak{o}$. Pour $\sharp=iso$ ou $an$, l'alg\`ebre de Lie $\mathfrak{g}_{\sharp}(F)$ s'identifie naturellement \`a un sous-espace de l'alg\`ebre des endomorphismes de $V_{\sharp}$. Ainsi, $trace (XY)$ est bien d\'efinie pour $X,Y\in \mathfrak{g}_{\sharp}(F)$. On introduit une transformation de Fourier $f\mapsto \hat{f}$ dans $C_{c}^{\infty}(\mathfrak{g}_{\sharp}(F))$ par la formule usuelle
$$\hat{f}(X)=\int_{\mathfrak{g}_{\sharp}(F)}f(Y)\psi_{F}(trace(XY))\,dY.$$
Le mesure $dY$ est la mesure de Haar sur $\mathfrak{g}_{\sharp}(F)$ qui est autoduale, c'est-\`a-dire que $\hat{\hat{f}}(X)=f(-X)$ pour tous $f$, $X$. L'exponentielle $exp$ est bien d\'efinie dans un voisinage de $0$ dans $\mathfrak{g}_{\sharp}(F)$, \`a valeurs dans un voisinage de $1$ dans $G_{\sharp}(F)$. On munit $G_{\sharp}(F)$ de la mesure de Haar telle que jacobien de l'exponentielle vaille $1$ au point $0$.

Soit $x\in G_{\sharp}(F)$ un \'el\'ement fortement r\'egulier et compact (cf. \cite{W5}1.2). Soit $f\in C_{c}^{\infty}(G_{\sharp}(F))$. On d\'efinit l'int\'egrale orbitale
$$J(x,f)=D^{G_{\sharp}}(x)^{1/2}\int_{G_{\sharp}(F)}f(g^{-1}xg)\,dg,$$
o\`u $D^{G_{\sharp}}$ est le module de Weyl usuel. Si maintenant $f=(f_{iso},f_{an})\in C_{c}^{\infty}(G_{iso}(F))\oplus C_{c}^{\infty}(G_{an}(F))$ et si $x\in G_{\sharp}(F)$ est comme ci-dessus, on pose 
$J(x,f)=J(x,f_{\sharp})$. 

On a d\'efini un espace ${\cal R}^{par}$ en \cite{W5} 1.5. On d\'efinit une application lin\'eaire $\Psi:{\cal R}^{par}\to C_{c}^{\infty}(G_{iso}(F))\oplus C_{c}^{\infty}(G_{an}(F))$ de la fa\c{c}on suivante. Par lin\'earit\'e, il suffit de fixer $\sharp=iso$ ou $an$ et $(n',n'')\in D_{\sharp}(n)$   et de la d\'efinir sur la composante $C'_{n'}\otimes C''_{n'',\sharp}$ de ${\cal R}^{par}$. Soit donc $\varphi\in C'_{n'}\otimes C''_{n'',\sharp}$. C'est une fonction sur ${\bf SO}(2n'+1;{\mathbb F}_{q})\times {\bf O}(2n'')_{\sharp}({\mathbb F}_{q})$.  On introduit le sous-groupe $K_{n',n''}^{\pm}$ de $G_{\sharp}(F)$, cf. \cite{W5} 1.2. Le groupe pr\'ec\'edent s'identifie \`a $K_{n',n''}^{\pm}/K_{n',n''}^{u}$. Alors $\varphi$ appara\^{\i}t comme une fonction sur $K_{n',n''}^{\pm}$, invariante par $K_{n',n''}^{u}$. On la prolonge par $0$ hors de $K_{n',n''}^{\pm}$ et on la multiplie par $mes(K_{n',n''}^{\pm})^{-1}$. On obtient un \'el\'ement de $ C_{c}^{\infty}(G_{\sharp}(F))$ qui est l'une des composantes de $\Psi(\varphi)$. L'autre composante est nulle. 

Soit $\sharp=iso$ ou $an$ et soit $\pi\in Irr_{unip,\sharp}$, cf. \cite{W5} 1.3. Pour $f\in C_{c}^{\infty}(G_{\sharp}(F))$, on d\'efinit l'op\'erateur $\pi(f)$. Il est de rang fini et poss\`ede une trace. D'apr\`es Harish-Chandra, il existe une fonction localement int\'egrable $\Theta_{\pi}$ sur $G_{\sharp}(F)$, localement constante sur les \'el\'ements fortement r\'eguliers, de sorte que
$$trace(\pi(f))=\int_{G_{\sharp}(F)}\Theta_{\pi}(x)f(x)\,dx$$
pour tout $f$. Posons $f_{\pi}=\Psi\circ proj_{cusp}\circ Res(\pi)$, avec  les notations de \cite{W5} 1.5.

\ass{Proposition}{Soit $x\in G_{\sharp}(F)$ un \'el\'ement fortement r\'egulier et compact. On a l'\'egalit\'e
$$D^{G_{\sharp}}(x)^{1/2}\overline{\Theta_{\pi}(x)}=J(x,f_{\pi}).$$}

Cela r\'esulte de \cite{MW} th\'eor\`eme 1.9 et lemme 4.2.

\subsection{Un r\'esultat de transfert de fonctions}

Pour $\sharp=iso$ ou $an$ et pour $f\in C_{c}^{\infty}(G_{\sharp}(F))$, on dit que $f$ est cuspidale si et seulement si ses int\'egrales orbitales sont nulles en tout point fortement r\'egulier et non elliptique de $G_{\sharp}(F)$. On note $C_{cusp}^{\infty}(G_{\sharp}(F))$ l'espace des fonctions cuspidales. L'application $\Psi:{\cal R}^{par}\to C_{c}^{\infty}(G_{iso}(F))\oplus C_{c}^{\infty}(G_{an}(F))$ introduite en 1.1 envoie l'espace ${\cal R}^{par}_{cusp}$ de \cite{W5} 1.5 dans $C_{cusp}^{\infty}(G_{iso}(F))\oplus C_{cusp}^{\infty}(G_{an}(F))$. Pour un \'el\'ement $f$ de ce dernier espace, on note $f_{iso}$ et $f_{an}$ ses deux composantes. 

Il nous faut adapter les d\'efinitions usuelles de l'endoscopie \`a notre cas o\`u nous travaillons simultan\'ement avec nos deux groupes $G_{iso}$ et $G_{an}$. Pour $x,y\in G_{iso}(F)\cup G_{an}(F)$, on dit que $x$ et $y$ sont conjugu\'es si et seulement s'il existe $\sharp=iso$ ou $an$ tels que $x,y\in G_{\sharp}(F)$ et $x$ et $y$ sont conjugu\'es par un \'el\'ement de $G_{\sharp}(F)$. Il y a une correspondance bijective entre les classes de conjugaison stable dans $G_{iso}(F)$ form\'ees d'\'el\'ements elliptiques et fortement r\'eguliers et  les classes de conjugaison stable dans $G_{an}(F)$ form\'ees d'\'el\'ements elliptiques et fortement r\'eguliers. Appelons classe totale de conjugaison stable (form\'ee d'\'el\'ements elliptiques et fortement r\'eguliers) la r\'eunion d'une telle classe dans $G_{iso}(F)$ et de celle qui lui correspond dans $G_{an}(F)$. On sait d'ailleurs que chacune de ces classes contient le m\^eme nombre de classes de conjugaison (ordinaires).

Soit $f\in C_{cusp}^{\infty}(G_{iso}(F))\oplus C_{cusp}^{\infty}(G_{an}(F))$ et soit $x\in G_{iso}(F)$ un \'el\'ement fortement r\'egulier et elliptique. On d\'efinit l'int\'egrale orbitale stable $S(x,f)$ par
$$S(x,f)=\sum_{y}J(y,f) ,$$
o\`u $y$ parcourt les \'el\'ements de la classe totale de conjugaison stable de $x$, \`a conjugaison pr\`es.   On note 
$z(x)$  le nombre de termes $y $   intervenant dans cette somme. Pour $\sharp=iso$ ou $an$ et $f\in C_{c}^{\infty}(G_{\sharp}(F))$, on d\'efinit $S(x,f)$ en consid\'erant que $f$ est un \'el\'ement de $C_{cusp}^{\infty}(G_{iso}(F))\oplus C_{cusp}^{\infty}(G_{an}(F))$ dont la composante sur $G_{\sharp}(F)$ est la fonction donn\'ee et dont l'autre composante est nulle.

Soit $(n_{1},n_{2})\in D(n)$. Ce couple d\'etermine (\`a conjugaison pr\`es) un \'el\'ement $h\in Sp(2n;{\mathbb C})$ tel que $h^2=1$ et une donn\'ee endoscopique de $G_{iso}$ ou $G_{an}$ dont le groupe endoscopique est $G_{n_{1},iso}\times G_{n_{2},iso}$, cf. \cite{W5} 2.1. On note $\Delta_{h}$ le facteur de transfert relatif \`a cette donn\'ee (il est uniquement d\'etermin\'e). Soit $f\in C_{cusp}^{\infty}(G_{iso}(F))\oplus C_{cusp}^{\infty}(G_{an}(F))$ et $(x_{1},x_{2})\in G_{n_{1},iso}(F)\times G_{n_{2},iso}(F)$ un couple $n$-r\'egulier form\'e  d'\'el\'ements   elliptiques. On entend par $n$-r\'egulier le fait que la classe de conjugaison stable de  $(x_{1},x_{2})$ correspond par endoscopie \`a une classe de conjugaison stable fortement r\'eguli\`ere dans $G_{iso}(F)$. On pose
$$J^{ endo}(x_{1},x_{2},f)=\sum_{x}\Delta_{h}((x_{1},x_{2}),x)J(x,f) ,$$
o\`u $x$ parcourt les \'el\'ements de $G_{iso}(F)\cup G_{an}(F)$ fortement r\'eguliers et elliptiques, \`a conjugaison pr\`es. Bien s\^ur, on peut se limiter aux \'el\'ements dans la classe  totale de conjugaison stable qui correspond \`a $(x_{1},x_{2})$ (en dehors, les facteurs de transfert sont nuls).  La somme est donc finie. Dans le cas o\`u $n_{1}=n$ et $n_{2}=0$, auquel cas $(x_{1},x_{2})$ se r\'eduit \`a un seul \'el\'ement $x=x_{1}$, on retrouve la d\'efinition pr\'ec\'edente: $J^{endo}(x_{1},x_{2},f)=S(x,f)$. Comme ci-dessus, pour $\sharp=iso$ ou $an$ et $f\in C_{c}^{\infty}(G_{\sharp}(F))$, on d\'efinit $J^{ endo}(x_{1},x_{2},f)$ en consid\'erant $f$ comme un \'el\'ement de $C_{cusp}^{\infty}(G_{iso}(F))\oplus C_{cusp}^{\infty}(G_{an}(F))$ dont une composante est nulle. 

Soient $(n_{1},n_{2})\in D(n)$, $f\in C_{cusp}^{\infty}(G_{iso}(F))\oplus C_{cusp}^{\infty}(G_{an}(F))$, $f_{1}\in C_{cusp}^{\infty}(G_{n_{1},iso}(F))\oplus C_{cusp}^{\infty}(G_{n_{1},an}(F))$ et $f_{2}\in C_{cusp}^{\infty}(G_{n_{2},iso}(F))\oplus C_{cusp}^{\infty}(G_{n_{2}, an}(F))$. On dit que $f_{1}\otimes f_{2}$ est un transfert de $f$ relatif \`a $(n_{1},n_{2})$ (ou \`a $h$) si et seulement si on a l'\'egalit\'e
$$S(x_{1},f_{1})S(x_{2},f_{2})=J^{endo}(x_{1},x_{2},f)$$
pour tout couple $n$-r\'egulier $(x_{1},x_{2})\in G_{n_{1},iso}(F)\times G_{n_{2},iso}(F)$   form\'e de deux \'el\'ements  elliptiques. 

Soit $(r',r'',N',N'')\in \Gamma$, cf. \cite{W5}1.8, soient $\varphi'\in {\mathbb C}[\hat{W}_{N'}]_{cusp}$ et $\varphi''\in {\mathbb C}[\hat{W}_{N''}]_{cusp}$.  Pour un nombre r\'eel $x$, on note $[x]$ sa partie enti\`ere. Posons
$$r'_{1}=sup([\frac{r'+r''}{2}],-[\frac{r'+r''}{2}]-1),\,\, r''_{1}=\vert [\frac{r'+r''+1}{2}]\vert ,$$
$$r'_{2}=sup([\frac{r'-r''}{2}],-[\frac{r'-r''}{2}]-1),\,\, r''_{2}=\vert [\frac{r'-r''+1}{2}]\vert .$$
$$n_{1}=r_{1}^{'2}+r'_{1}+r^{''2}_{1}+N',\,\,n_{2}=r_{2}^{'2}+r'_{2}+r^{''2}_{2}+N'',$$
$$(N'_{1},N''_{1})=(N',0),\,,\ (N'_{2},N''_{2})=(N'',0),$$
$$\gamma_{1}=(r'_{1},r''_{1},N'_{1},N''_{1}),\,\, \gamma_{2}=(r'_{2},r''_{2},N'_{2},N''_{2}).$$
 D'apr\`es la relation (2) du paragraphe 4 ci-dessous, la d\'efinition de $n_{1}$ et $n_{2}$ peut  se r\'ecrire
$$n_{1}=\frac{(r'+r'')^2+(r'+r''+1)^2-1}{4}+N',\,\, n_{2}=\frac{(r'-r'')^2+(r'-r''+1)^2-1}{4}+N''.$$

On v\'erifie que $(n_{1},n_{2})\in D(n)$, que $\gamma_{1}\in \Gamma_{n_{1}}$ et $\gamma_{2}\in \Gamma_{n_{2}}$. L'\'el\'ement $\varphi'$ peut \^etre consid\'er\'e comme un \'el\'ement de ${\cal R}_{cusp}(\gamma_{1})$ et l'\'el\'ement $\varphi''$ peut \^etre consid\'er\'e comme un \'el\'ement de ${\cal R}_{cusp}(\gamma_{2})$. Posons $\varphi=\varphi'\otimes \varphi''$. C'est un \'el\'ement de ${\cal R}_{cusp}(\gamma)$. Posons
$$f=\Psi\circ k\circ\rho\iota(\varphi),\,\, f_{1}=\Psi\circ k\circ\rho\iota(\varphi'),\,\, f_{2}=\Psi\circ k\circ\rho\iota(\varphi''),$$
avec les notations de \cite{W5} 1.9, 1.10. 
Bien s\^ur, dans les deux derni\`eres d\'efinitions, ce sont les applications $\Psi$, $k$ et $\rho\iota$ relatives \`a $n_{1}$ et $n_{2}$ qui interviennent. 

\ass{Proposition}{(i) Sous les hypoth\`eses ci-dessus, $f_{1}\otimes f_{2}$ est un transfert de $f$ relatif \`a $(n_{1},n_{2})$.

(ii) Pour $(\bar{n}_{1},\bar{n}_{2})\in D(n)$ diff\'erent de $(n_{1},n_{2})$, le transfert de $f$ relatif \`a $(\bar{n}_{1},\bar{n}_{2})$ est nul.}

Nous d\'emontrerons cette proposition dans la section 3. Nous l'admettons pour la suite de la pr\'esente section.

\bigskip

\subsection{Transfert de repr\'esentations elliptiques}
Soit $(n_{1},n_{2})\in D(n)$, auquel on associe un \'el\'ement $h\in Sp(2n;{\mathbb C})$ tel que $h^2=1$. Soient $(\lambda_{1},s_{1};1)\in \mathfrak{St}_{n_{1},unip,disc}$ et $(\lambda_{2},s_{2},1)\in \mathfrak{St}_{n_{2},unip,disc}$, cf. \cite{W5} 2.4.  Comme en \cite{W5} 2.1, on associe \`a ces donn\'ees un triplet $(\lambda,s,h)\in \mathfrak{Endo}_{unip-disc}$. 

\ass{Th\'eor\`eme}{On a les \'egalit\'es 

$transfert_{h,iso}(\Pi_{iso}(\lambda_{1},s_{1},1)\otimes \Pi_{iso}(\lambda_{2},s_{2},1))=\Pi_{iso}(\lambda,s,h)$;

$transfert_{h,an}(\Pi_{iso}(\lambda_{1},s_{1},1)\otimes \Pi_{iso}(\lambda_{2},s_{2},1))=-\Pi_{an}(\lambda,s,h)$.}

Preuve. 
Notons $\Theta_{1,iso}$, $\Theta_{1,an}$, $\Theta_{2,iso}$, $\Theta_{2,an}$, $\Theta_{iso}$, $\Theta_{an}$ les caract\`eres-fonctions des repr\'esentations
 $\Pi_{iso}(\lambda_{1},s_{1},1)$, $\Pi_{an}(\lambda_{1},s_{1},1)$, $\Pi_{iso}(\lambda_{2},s_{2},1)$, $\Pi_{an}(\lambda_{2},s_{2},1)$, $\Pi_{iso}(\lambda,s,h)$, $\Pi_{an}(\lambda,s,h)$. Les \'egalit\'es \`a d\'emontrer se traduisent par des \'egalit\'es entre ces caract\`eres. Toutes les repr\'esentations intervenant sont elliptiques. D'apr\`es \cite{Ar2} th\'eor\`eme 6.2, il suffit donc de d\'emontrer ces \'egalit\'es restreintes aux \'el\'ements elliptiques des diff\'erents groupes. Ecrivons ces \'egalit\'es. Posons $sgn_{iso}=1$ et $sgn_{an}=-1$. Pour $\sharp=iso$ ou $an$ et pour $x\in G_{\sharp}(F)$ fortement r\'egulier et elliptique, on doit avoir
$$(1) \qquad D^{G_{\sharp}}(x)^{1/2}\Theta_{\sharp}(x)=sgn_{\sharp}\sum_{x_{1},x_{2}}D^{G_{n_{1},iso}}(x_{1})^{1/2}D^{G_{n_{2},iso}}(x_{2})^{1/2}\Delta_{h}((x_{1},x_{2}),x)$$
$$\Theta_{1,iso}(x_{1})\Theta_{2,iso}(x_{2}),$$
o\`u  $x_{1}$, resp. $x_{2}$, parcourt les \'el\'ements fortement r\'eguliers et elliptiques de  $G_{n_{1},iso}(F)$, resp. $G_{n_{2},iso}(F)$, \`a conjugaison stable pr\`es. Ici encore, on peut se limiter aux couples $(x_{1},x_{2})$ correspondant par endoscopie \`a la classe de conjugaison stable de $x$. La somme est finie.

Posons $f_{1}=\Psi\circ proj_{cusp}\circ Res\circ D\circ \Pi(\lambda_{1},s_{1},1)$, $f_{2}=\Psi\circ proj_{cusp}\circ Res\circ D\circ \Pi(\lambda_{2},s_{2},1)$, $f=\Psi\circ proj_{cusp}\circ Res\circ D\circ \Pi(\lambda,s,h)$.  
Pour $\sharp=iso$ ou $an$, on a l'\'egalit\'e $proj_{cusp}\circ Res\circ D\circ \Pi_{\sharp}(\lambda,s,h)=(-1)^nsgn_{\sharp}proj_{cusp}\circ Res\circ \Pi_{\sharp}(\lambda,s,h)$, cf. \cite{MW} corollaire 5.7. D'apr\`es la proposition 1.1, on a donc l'\'egalit\'e 
$$D^{G_{\sharp}}(x)^{1/2}\overline{\Theta_{\sharp}(x)}=(-1)^n sgn_{\sharp}J(x,f ).$$
On a des \'egalit\'es analogues pour $\Theta_{1}(x_{1})$ et $\Theta_{2}(x_{2})$.  Ces derni\`eres \'egalit\'es  entra\^{\i}nent 

(2) soient $j=1,2$, $\sharp_{j}=iso$ ou $an$, $x_{j}\in G_{n_{j},iso}(F)$ et $y_{j}\in G_{n_{j},\sharp_{j}}(F)$; supposons que $x_{j}$ et $y_{j}$  sont fortement r\'eguliers et elliptiques et que leurs classes de conjugaison stable sont \'egales si $\sharp_{j}=iso$ ou se correspondent si $\sharp_{j}=an$; alors $J(x_{j},f_{j})=J(y_{j},f_{j})$. 

Cela traduit le fait que $\Theta_{j,iso}$ est stable et que $-\Theta_{j,an}$ en est son transfert, cf. \cite{W5} 2.1.

 Parce que les facteurs de transfert valent $\pm 1$ et sont donc \`a valeurs r\'eelles, l'\'egalit\'e (1)  se r\'ecrit
$$(3) \qquad J(x,f)=\sum_{x_{1},x_{2}}\Delta_{h}((x_{1},x_{2}),x) J(x_{1},f_{1})J(x_{2},f_{2})$$
pour tout $x\in G_{iso}(F)\cup G_{an}(F)$ fortement r\'egulier et elliptique.

La propri\'et\'e de base de l'endoscopie est que les int\'egrales orbitales $J(x,f)$ sont d\'etermin\'ees   par les int\'egrales endoscopiques $J^{endo}(\bar{x}_{1},\bar{x}_{2},f)$ quand 
$(\bar{n}_{1},\bar{n}_{2})$ d\'ecrit $ D(n)$ et $(\bar{x}_{1},\bar{x}_{2})\in G_{\bar{n}_{1},iso}(F)\times G_{\bar{n}_{2},iso}(F)$ d\'ecrit les couples $n$-r\'eguliers form\'es d'\'el\'ements elliptiques. Donc (3) \'equivaut \`a ce que, pour toutes  donn\'ees $(\bar{n}_{1},\bar{n}_{2})$ et $(\bar{x}_{1},\bar{x}_{2})$ comme ci-dessus, on ait l'\'egalit\'e
$$J^{endo}(\bar{x}_{1},\bar{x}_{2},f)=\sum_{x}\Delta_{\bar{h}}((\bar{x}_{1},\bar{x}_{2}),x)J(x,f)$$
$$=\sum_{x}\Delta_{\bar{h}}((\bar{x}_{1},\bar{x}_{2}),x)\sum_{x_{1},x_{2}}\Delta_{h}((x_{1},x_{2}),x) J(x_{1},f_{1})J(x_{2},f_{2}),$$
o\`u on somme sur les \'el\'ements $x\in G_{iso}(F)\cup G_{an}(F)$ fortement r\'eguliers et elliptiques modulo conjugaison. On a not\'e $\bar{h}$ l'\'el\'ement de $Sp(2n;{\mathbb C})$ tel que $\bar{h}^2=1$ associ\'e \`a $(\bar{n}_{1},\bar{n}_{2})$. On r\'ecrit l'\'egalit\'e ci-dessus:
$$J^{endo}(\bar{x}_{1},\bar{x}_{2},f)=\sum_{x_{1},x_{2}}  J(x_{1},f_{1})J(x_{2},f_{2})\sum_{x}\Delta_{\bar{h}}((\bar{x}_{1},\bar{x}_{2}),x) \Delta_{h}((x_{1},x_{2}),x).$$
La somme int\'erieure en $x$ est nulle sauf si $(\bar{n}_{1},\bar{n}_{2})=(n_{1},n_{2})$ et, \`a conjugaison stable pr\`es, $(\bar{x}_{1},\bar{x}_{2})=(x_{1},x_{2})$. Si ces conditions sont v\'erifi\'ees, la somme vaut le nombre de  $x$, pris \`a conjugaison pr\`es, qui interviennent dans la somme. C'est le nombre $z(x)$ d\'efini en 1.2 pour l'un quelconque de ces $x$. Or on v\'erifie facilement l'\'egalit\'e $z(x)=z(\bar{x}_{1})z(\bar{x}_{2})$. 
 Dans le cas o\`u $(\bar{n}_{1},\bar{n}_{2})=(n_{1},n_{2})$, on  obtient donc
 $$J^{endo}(\bar{x}_{1},\bar{x}_{2},f)=z(\bar{x}_{1})z(\bar{x}_{2})J(\bar{x}_{1},f_{1})J(\bar{x}_{2},f_{2}).$$
 Mais la propri\'et\'e (2) entra\^{\i}ne que, pour $j=1,2$, $z(\bar{x}_{j})J(\bar{x}_{j},f_{j})=S(\bar{x}_{j},f_{j})$. L'\'egalit\'e ci-dessus se transforme en
 $$J^{endo}(\bar{x}_{1},\bar{x}_{2},f)=S(\bar{x}_{1},f_{1})S(\bar{x}_{2},f_{2}).$$
 Autrement dit, $f_{1}\otimes f_{2}$ est un transfert de $f$. R\'esumons: l'\'egalit\'e (3) \'equivaut aux assertions
 
 (4)(a) $f_{1}\otimes f_{2}$ est un transfert de $f$ relatif \`a $(n_{1},n_{2})$;
 
 (4)(b) pour $(\bar{n}_{1},\bar{n}_{2})\in D(n)$ diff\'erent de $(n_{1},n_{2})$, le transfert de $f$ relatif \`a $(\bar{n}_{1},\bar{n}_{2})$ est nul.

Calculons $f$. D'apr\`es le lemme 2.3 de \cite{W5} et parce que ${\cal F}^{par}$ est une involution, on a $proj_{cusp}={\cal F}^{par}\circ proj_{cusp}\circ \mathfrak{F}^{par}$. Par d\'efinition de $\mathfrak{F}^{par}$, on a aussi $\mathfrak{F}^{par}\circ Res\circ D=Res\circ D\circ {\cal F}$. D'o\`u
$$f=\Psi\circ {\cal F}^{par}\circ proj_{cusp}\circ Res\circ D\circ {\cal F}\circ \Pi(\lambda,s,h))=\Psi\circ {\cal F}^{par}\circ proj_{cusp}\circ Res\circ D\circ \Pi(\lambda,h,s).$$
La d\'ecomposition de $\lambda$ associ\'ee \`a $h$ est $\lambda=\lambda_{1}\cup \lambda_{2}$. Il r\'esulte des d\'efinitions que 
$$\Pi(\lambda,h,s)=\sum_{\epsilon_{1},\epsilon_{2}}\pi(\lambda_{1},\epsilon_{1},\lambda_{2},\epsilon_{2})\epsilon_{1}(s_{1})\epsilon_{2}(s_{2}),$$
o\`u $(\epsilon_{1},\epsilon_{2})$ parcourt $\{\pm 1\}^{Jord_{bp}(\lambda_{1})}\times \{\pm 1\}^{Jord_{bp}(\lambda_{2})}$ et les termes $\epsilon_{1}(s_{1})$ et $\epsilon_{2}(s_{2})$ sont d\'efinis en interpr\'etant $\epsilon_{1}$ et $\epsilon_{2}$ comme des \'el\'ements de ${\bf Z}(\lambda_{1},1)^{\vee}$ et ${\bf Z}(\lambda_{2},1)^{\vee}$, cf. \cite{W5} 1.3. La proposition \cite{W5} 1.11 calcule $Res\circ D(\pi(\lambda_{1},\epsilon_{1},\lambda_{2},\epsilon_{2}) )$.

{\bf Changement de notations.} Dans la formule de cette proposition figure un isomorphisme $j$ entre deux espaces, l'un d'eux \'etant l'espace ${\cal R}$. Distinguer ces deux espaces nous a \'et\'e utile dans la deuxi\`eme section de \cite{W5}. Maintenant, cela ne nous sert plus. Pour simplifier, on identifie par $j$ les deux espaces en question et on fait dispara\^{\i}tre $j$ de la notation. 

\bigskip

On obtient
$$Res\circ D\circ \Pi(\lambda,h,s)=\sum_{\epsilon_{1},\epsilon_{2}}\epsilon_{1}(s_{1})\epsilon_{2}(s_{2})Rep\circ \rho\iota(\boldsymbol{\rho}_{\lambda_{1},\epsilon_{1}}\otimes \boldsymbol{\rho}_{\lambda_{2},\epsilon_{2}}).$$
Les applications $Rep$, $k$, ${\cal F}^{par}$ et $\rho\iota$ commutent aux projections cuspidales. De plus ${\cal F}^{par}\circ Rep=k$. Il en r\'esulte que
$$f=\sum_{\epsilon_{1},\epsilon_{2}}\epsilon_{1}(s_{1})\epsilon_{2}(s_{2})\Psi\circ k\circ \rho\iota\circ proj_{cusp}(\boldsymbol{\rho}_{\lambda_{1},\epsilon_{1}}\otimes \boldsymbol{\rho}_{\lambda_{2},\epsilon_{2}}).$$
Un calcul analogue s'applique \`a $f_{1}$ et $f_{2}$. On a cette fois ${\cal F}(\lambda_{j},s_{j},1)=(\lambda_{j},1,s_{j})$ pour $j=1,2$ et la d\'ecomposition de $\lambda_{j}$ associ\'ee \`a $1$ est $\lambda_{j}=\lambda_{j}\cup \emptyset$. Autrement dit, les deuxi\`emes composantes de la formule ci-dessus disparaissent. On obtient
$$f_{1}=\sum_{\epsilon_{1}}\epsilon_{1}(s_{1})\Psi\circ k\circ \rho\iota\circ proj_{cusp}(\boldsymbol{\rho}_{\lambda_{1},\epsilon_{1}}),$$
$$f_{2}=\sum_{\epsilon_{2}}\epsilon_{2}(s_{2})\Psi\circ k\circ \rho\iota\circ proj_{cusp}(\boldsymbol{\rho}_{\lambda_{2},\epsilon_{2}}),$$
 o\`u  les repr\'esentations $\boldsymbol{\rho}_{\lambda_{1},\epsilon_{2}}$ et $\boldsymbol{\rho}_{\lambda_{2},\epsilon_{2}}$  sont assimil\'ees \`a des produits tensoriels dont le deuxi\`eme terme est trivial. Pour $(\epsilon_{1},\epsilon_{2})\in \{\pm 1\}^{Jord_{bp}(\lambda_{1})}\times \{\pm 1\}^{Jord_{bp}(\lambda_{2})}$, posons
 $$f_{\epsilon_{1},\epsilon_{2}}=\Psi\circ k\circ proj_{cusp}(\boldsymbol{\rho}_{\lambda_{1},\epsilon_{1}}\otimes \boldsymbol{\rho}_{\lambda_{2},\epsilon_{2}}),$$
 $$f_{\epsilon_{1}}=\Psi\circ k\circ proj_{cusp}(\boldsymbol{\rho}_{\lambda_{1},\epsilon_{1}}),$$
 $$f_{\epsilon_{2}}=\Psi\circ k\circ proj_{cusp}(\boldsymbol{\rho}_{\lambda_{2},\epsilon_{2}}).$$
 Les propri\'et\'es (4) r\'esultent des propri\'et\'es suivantes, pour tout $(\epsilon_{1},\epsilon_{2})$:
 
 (5)(a) $f_{\epsilon_{1}}\otimes f_{\epsilon_{2}}$ est un transfert de $f_{\epsilon_{1},\epsilon_{2}}$ relatif \`a $(n_{1},n_{2})$;
 
  (5)(b) pour $(\bar{n}_{1},\bar{n}_{2})\in D(n)$ diff\'erent de $(n_{1},n_{2})$, le transfert de $f_{\epsilon_{1},\epsilon_{2}}$ relatif \`a $(\bar{n}_{1},\bar{n}_{2})$ est nul.

Fixons $\epsilon_{1}$ et $\epsilon_{2}$. A $(\lambda_{1},\epsilon_{1})$ et $(\lambda_{2},\epsilon_{2})$  sont associ\'es des couples $(k_{1},N_{1})$ et $(k_{2},N_{2})$, puis un triplet $\gamma=(r',r'',N_{1},N_{2})\in \Gamma$, cf. \cite{W5} 1.11. L'\'el\'ement $proj_{cusp}(\boldsymbol{\rho}_{\lambda_{1},\epsilon_{1}}\otimes \boldsymbol{\rho}_{\lambda_{2},\epsilon_{2}})$ appartient \`a ${\cal R}_{cusp}(\gamma)$. De m\^eme, en rempla\c{c}ant $n$ par $n_{j}$ pour $j=1,2$, \`a $(\lambda_{j},\epsilon_{j})$ sont associ\'es des termes $\gamma_{1}$ et $\gamma_{2}$ analogues. Ainsi qu'on l'a remarqu\'e ci-dessus, il faut consid\'erer que $(\lambda_{j},\epsilon_{j})$ sont le premier couple  d'un quadruplet $(\lambda_{j}^+,\epsilon_{j}^+,\lambda_{j}^-,\epsilon_{j}^-)$ dont le second couple est trivial.  On voit que $\gamma_{j}$ est de la forme $(r'_{j},r''_{j},N_{j},0)$. L'\'el\'ement $proj_{cusp}(\boldsymbol{\rho}_{\lambda_{j},\epsilon_{j}})$ appartient \`a ${\cal R}_{cusp}(\gamma_{j})$. En consid\'erant la recette de \cite{W5} 1.11 qui calcule $r',r'',r'_{1},r''_{1},r'_{2},r''_{2}$ en fonction de $k_{1}$ et $k_{2}$, on constate que $r'_{1},r''_{1},r'_{2},r''_{2}$ se d\'eduisent de $(r',r'')$ par les formules de 1.2. On est alors dans la situation de ce paragraphe, les fonctions $\varphi'$ et $\varphi''$ \'etant respectivement $proj_{cusp}(\boldsymbol{\rho}_{\lambda_{1},\epsilon_{1}})$ et $proj_{cusp}(\boldsymbol{\rho}_{\lambda_{2},\epsilon_{2}})$. La proposition 1.2  affirme que les propri\'et\'es (5)(a) et (5)(b) sont v\'erifi\'ees. Cela ach\`eve la d\'emonstration. $\square$

\bigskip

\subsection{D\'emonstration du th\'eor\`eme 2.1 de \cite{W5}}
Soit $h\in Sp(2n;{\mathbb C})$ tel que $h^2=1$, auquel est associ\'e un couple $(n_{1},n_{2})\in D(n)$. Pour $j=1,2$, soit $(\lambda_{j},s_{j},1)\in \mathfrak{St}_{n_{j},tunip}$. Pour $j=1,2$, notons $\cup_{b\in B_{j}}\{s_{j,b},s_{j,b}^{-1}\}$ l'ensemble des valeurs propres de $s_{j}$ autres que $\pm 1$ et notons 
$$\lambda_{j}=\lambda_{j}^+\cup \lambda^-_{j}\cup\bigcup_{b\in B_{j}}(\lambda_{j,b}\cup\lambda_{j,b})$$
la d\'ecomposition de $\lambda_{j}$ associ\'ee \`a $s_{j}$, cf. \cite{W5} 2.1. Posons 
$\lambda_{j,0}=\lambda_{j}^+\cup \lambda_{j}^-$, $n_{j,0}=S(\lambda_{j,0})/2$ et $m_{j,b}=S(\lambda_{j,b})$ pour $b\in B_{j}$. Introduisons un sous-groupe parabolique $P_{j}$ de $G_{n_{j},iso}$ de composante de Levi
$$M_{j}=(\prod_{b\in B_{j}}GL(m_{j,b}))\times G_{n_{j,0},iso}.$$
Pour tout $b\in B_{j}$, soit $\chi_{j,b}$ le caract\`ere non ramifi\'e de $F^{\times}$ tel que $\chi_{j,b}(\varpi)=s_{j,b}$. L'\'el\'ement $s_{j}$ se restreint en un \'el\'ement  $s_{j,0}\in Sp(2n_{j,0};{\mathbb C})$ de sorte que $\lambda_{j,0}=\lambda_{j}^+\cup \lambda_{j}^-$ soit la d\'ecomposition de $\lambda_{j,0}$ associ\'ee \`a cet \'el\'ement. 
Introduisons la repr\'esentation de $M_{j}(F)$
$$\sigma_{j,b}=(\otimes_{b\in B_{j}}(\chi_{j,b}\circ det)st_{\lambda_{j,b}}))\otimes \Pi_{iso}(\lambda_{j,0},s_{j,0},1)$$
(on rappelle que $st_{m}$ est la repr\'esentation de Steinberg de $GL(m;F)$). 
Il est connu que
$$\Pi_{iso}(\lambda_{j},s_{j},1)=Ind_{P_{j}}^{G_{n_{j},iso}}(\sigma_{j}).$$

 Introduisons le triplet $(\lambda,s,h)$ associ\'e \`a $h$, $(\lambda_{1},s_{1},1)$ et $(\lambda_{2},s_{2},1)$. Posons $n_{0}=n_{0,1}+n_{0,2}$. Soit $\sharp=iso$ ou $an$. Dans le cas $\sharp=an$, supposons provisoirement $n_{0}\geq1$. On introduit un sous-groupe parabolique $P$ de $G_{\sharp}$ de composante de Levi
$$M=(\prod_{j=1,2,b\in B_{j}}GL(m_{j,b}))\times G_{n_{0},\sharp}.$$
Parce que l'\'el\'ement $h$ commute \`a $s$, il se restreint en un \'el\'ement $h_{0}\in Sp(2n_{0};{\mathbb C})$ et la repr\'esentation $\Pi_{\sharp}(\lambda_{0},s_{0},h_{0})$ de $G_{n_{0},\sharp}(F)$ est bien d\'efinie. On introduit la repr\'esentation de $M(F)$
$$\sigma=(\otimes_{j=1,2,b\in B_{j}}(\chi_{j,b}\circ det)st_{\lambda_{j,b}}))\otimes \Pi_{\sharp}(\lambda_{0},s_{0},h_{0}).$$
 
On v\'erifie que
$$\Pi_{\sharp}(\lambda,s,h)=Ind_{P}^{G_{\sharp}}(\sigma).$$
Pour $j=1,2$, le triplet $(\lambda_{j,0},s_{j,0},1)$ appartient \`a $\mathfrak{St}_{n_{j,0},unip-quad}$. Le th\'eor\`eme 1.3 entra\^{\i}ne que 
$$\Pi_{\sharp}(\lambda_{0},s_{0},h_{0})=sgn_{\sharp}transfert_{h_{0},\sharp}(\Pi_{iso}(\lambda_{0,1},s_{0,1},1)\otimes \Pi_{iso}(\lambda_{0,2},s_{0,2},1)).$$
Le transfert $transfert_{h_{0},\sharp}$ se prolonge en un transfert entre les groupes $M_{1}\times M_{2}$ et $M$ (le transfert \'etant l'identit\'e sur les composantes $GL(m_{j,b})$).  Pour celui-ci, $\sigma_{1}\otimes \sigma_{2}$ se transf\`ere en $sgn_{\sharp}\sigma$.  Mais on sait que le transfert commute \`a l'induction. Il r\'esulte alors des descriptions ci-dessus  que
$$\Pi_{\sharp}(\lambda,s,h)=sgn_{\sharp}transfert_{h,\sharp}(\Pi_{iso}(\lambda_{1},s_{1},1)\otimes \Pi_{iso}(\lambda_{2},s_{2},1)).$$
  Dans le cas particulier o\`u $\sharp=an$ et $n_{0}=0$, il n'y a plus de sous-groupe parabolique $P$. Le membre de droite ci-dessus est nul car c'est le transfert d'une induite \`a partir  d'un sous-groupe parabolique qui n'est pas relevant pour $G_{\sharp}$. Le membre de gauche est nul lui-aussi car, d'apr\`es la relation \cite{W5} 1.3(1), toutes les composantes irr\'eductibles de $\Pi(\lambda,s,h)$ vivent sur l'unique groupe $G_{iso}(F)$. L'\'egalit\'e ci-dessus est donc aussi  v\'erifi\'ee dans ce cas particulier.
  Elle   d\'emontre le th\'eor\`eme 2.1 de \cite{W5}. $\square$

\bigskip

\section{Trois lemmes de transfert pour des alg\`ebres de Lie}

\subsection{ Le cas sp\'ecial orthogonal impair}
 Cette section et la suivante   s'appuient sur les calculs d\'ej\`a faits dans \cite{MW}. On n'a gu\`ere envie de reproduire ces calculs, ni les d\'efinitions compliqu\'ees des objets qui y apparaissent. On se contentera d'indiquer les r\'ef\'erences n\'ecessaires.  Dans cette section, on va \'enoncer  des lemmes se transfert pour des alg\`ebres de Lie. Il y a trois cas: le cas sp\'ecial orthogonal impair, le cas sp\'ecial orthogonal pair et le cas unitaire. Les preuves sont similaires dans les trois cas (et beaucoup plus faciles dans le dernier). On n'\'ecrira que celle concernant le cas sp\'ecial orthogonal pair, qui est le plus subtil. Dans ces trois paragraphes, on oublie les objets $n$, $G_{iso}$ et $G_{an}$ pr\'ec\'edemment fix\'es afin de lib\'erer ces notations. On va introduire de nouveaux entiers $n$ et on conserve l'hypoth\`ese $p>6n+4$.

Dans ce paragraphe, on consid\`ere un entier $n\geq1$ et un \'el\'ement $\eta\in F^{\times}/F^{\times2}$. On construit comme en \cite{W5} 1.1 les deux espaces quadratiques $(V_{iso},Q_{iso})$ et $(V_{an},Q_{an})$ d\'efinis sur $F$ tels que $dim_{F}(V_{iso})=dim_{F}(V_{an})=2n+1$ et $\eta(Q_{iso})=\eta(Q_{an})=\eta$. On note $G_{iso}$ et $G_{an}$ leurs groupes sp\'eciaux orthogonaux. On consid\`ere quatre entiers $r',r'',N',N''\in {\mathbb N}$ tels que
$$r^{'2}+r^{''2}+2N'+2N''=2n+1,$$
$$r'\equiv 1+val_{F}(\eta)\,\,mod\,\,2{\mathbb Z},\,\, r''\equiv val_{F}(\eta)\,\,mod\,\,2{\mathbb Z}.$$
On pose $N=N'+N''$. On a d\'efini en \cite{MW} 3.11 et 3.12 une application
$${\cal Q}(r',r'')^{Lie}\circ \rho_{N}^*\circ \iota_{N',N''}:{\mathbb C}[\hat{W}_{N'}]_{cusp}\otimes {\mathbb C}[\hat{W}_{N''}]_{cusp}\to C_{cusp}^{\infty}(\mathfrak{g}_{iso}(F))\oplus C_{cusp}^{\infty}(\mathfrak{g}_{an}(F))$$
(en \cite{MW}, l'indice $cusp$ \'etait remplac\'e par $ell$ mais le sens \'etait le m\^eme).

Rappelons que, pour $m\in {\mathbb N}$, les classes de conjugaison dans $W_{m}$ sont param\'etr\'es par les paires de partitions $(\alpha,\beta)$ telles que $S(\alpha)+S(\beta)=m$.  Pour $w\in W_{m}$, notons $\varphi_{w}$ la fonction caract\'eristique de la classe de conjugaison de $w$. Alors ${\mathbb C}[\hat{W}_{m}]_{cusp}$ a pour base les $\varphi_{w}$ quand $w$ d\'ecrit, \`a conjugaison pr\`es, les \'el\'ements dont la classe est param\'etr\'ee par un couple de partitions de la forme $(\emptyset,\beta)$. On fixe des \'el\'ements $w'\in W_{N'}$ et $w''\in W_{N''}$ dont les classes sont param\'etr\'ees par des couples de cette forme. On pose
$$f={\cal Q}(r',r'')^{Lie}\circ\rho_{N}^*\circ\iota_{N',N''}(\varphi_{w'}\otimes \varphi_{w''}).$$

Soit $(n_{1},n_{2})\in D(n)$. A ce couple est associ\'ee une donn\'ee endoscopique de $G_{iso}$ et $G_{an}$. On utilise les m\^emes  d\'efinitions qu'en  1.2. Celles-ci se descendent aux alg\`ebres de Lie, c'est-\`a-dire  qu'il y a un transfert de $ C_{cusp}^{\infty}(\mathfrak{g}_{iso}(F))\oplus C_{cusp}^{\infty}(\mathfrak{g}_{an}(F))$ dans
$$\left( C_{cusp}^{\infty}(\mathfrak{g}_{n_{1},iso}(F))\oplus C_{cusp}^{\infty}(\mathfrak{g}_{n_{1},an}(F))\right)\otimes \left( C_{cusp}^{\infty}(\mathfrak{g}_{n_{2},iso}(F))\oplus C_{cusp}^{\infty}(\mathfrak{g}_{n_{2},an}(F))\right).$$
On le note $transfert_{n_{1},n_{2}}$. 
Le facteur de transfert est uniquement d\'efini. Les formules en sont donn\'ees en \cite{W3} proposition X.8. Il convient de supprimer les discriminants de Weyl de ces formules, que l'on a incorpor\'es aux int\'egrales orbitales. Le $\eta$ de la d\'efinition X.7 de \cite{W3} n'est pas notre pr\'esent $\eta$, c'est $(-1)^n\eta$. 

Fixons $\eta_{1},\eta_{2}\in F^{\times}/F^{\times2}$ tels que $\eta_{1}\eta_{2}=\eta$.  Notons $t'_{1}$ l'\'el\'ement de l'ensemble $\{\frac{r'+r''+1}{2},\frac{r'+r''-1}{2}\}$ qui est de la m\^eme parit\'e que $1+val_{F}(\eta_{1})$. Notons $t''_{1}$ l'autre \'el\'ement. Notons $t'_{2}$ l'\'el\'ement de l'ensemble $\{\frac{\vert r'-r''\vert +1}{2},\frac{\vert r'-r''\vert -1}{2}\}$ qui est de la m\^eme parit\'e que $1+val_{F}(\eta_{2})$. Notons $t''_{2}$ l'autre \'el\'ement. 
D\'efinissons deux entiers $n_{1}$ et $n_{2}$ par les formules
$$2n_{1}+1=t^{'2}_{1}+t^{''2}_{1}+2N',\,\,2n_{2}+1=t^{'2}_{2}+t^{''2}_{2}+2N'',$$
qui sont \'equivalentes \`a
$$n_{1}=\frac{(r'+r'')^2-1}{4}+N',\,\,n_{2}=\frac{(r'-r'')^2-1}{4}+N''.$$
On v\'erifie que $n_{1}+n_{2}=n$. Pour $j=1,2$, on peut consid\'erer $G_{n_{j},iso}$ et $G_{n_{j},an}$ comme les groupes sp\'eciaux orthogonaux d'espaces quadratiques de discriminants $\eta_{j}$.  On peut donc appliquer la m\^eme construction que ci-dessus, o\`u le couple $(N',N'')$ est remplac\'e par $(N',0)$ si $j=1$, par $(N'',0)$ si $j=2$. Cela nous permet de d\'efinir les \'el\'ements
$$f_{1}={\cal Q}(t'_{1},t''_{1})^{Lie}\circ\rho^*_{N'}\circ\iota_{N',0}(\varphi_{w'})\in  C_{cusp}^{\infty}(\mathfrak{g}_{n_{1},iso}(F))\oplus C_{cusp}^{\infty}(\mathfrak{g}_{n_{1},an}(F)),$$
$$f_{2}={\cal Q}(t'_{2},t''_{2})^{Lie}\circ\rho^*_{N''}\circ\iota_{N'',0}(\varphi_{w''})\in  C_{cusp}^{\infty}(\mathfrak{g}_{n_{2},iso}(F))\oplus C_{cusp}^{\infty}(\mathfrak{g}_{n_{2},an}(F)).$$

Notons $sgn$ l'unique caract\`ere non trivial de $\mathfrak{o}^{\times}/\mathfrak{o}^{\times2}$. D\'efinissons une constante $C$ par les formules suivantes:

si $r''\leq r'$, $C=sgn(-1)^{\frac{r'+r''-1}{2}}sgn(-\eta_{2}\varpi^{-val_{F}(\eta_{2})})^{val_{F}(\eta)}$;

 si $r'<r''$, $C=sgn(-1)^{\frac{r'+r''-1}{2}}sgn_{CD}(w'')sgn(-\eta_{2}\varpi^{-val_{F}(\eta_{2})}) ^{1+val_{F}(\eta)}$.

\ass{Lemme}{(i) On a l'\'egalit\'e $transfert_{n_{1},n_{2}}(f)=Cf_{1}\otimes f_{2}$.

(ii) Soit $(\bar{n}_{1},\bar{n}_{2})\in D(n)$ un couple diff\'erent de $(n_{1},n_{2})$. Alors  $transfert_{\bar{n}_{1},\bar{n}_{2}}(f)=0$.}

\bigskip

\subsection{Le cas sp\'ecial orthogonal pair}

Dans ce paragraphe, on consid\`ere un entier $n\geq1$ et un \'el\'ement $\eta\in F^{\times}/F^{\times2}$. On exclut le cas $n=1$, $\eta=1$. On construit comme en \cite{W5} 1.1 les deux espaces quadratiques $(V_{iso},Q_{iso})$ et $(V_{an},Q_{an})$ d\'efinis sur $F$ tels que $dim_{F}(V_{iso})=dim_{F}(V_{an})=2n$ et $\eta(Q_{iso})=\eta(Q_{an})=\eta$. On note $G_{\eta,iso}$ et $G_{\eta,an}$ leurs groupes sp\'eciaux orthogonaux. On consid\`ere quatre entiers $r',r'',N',N''\in {\mathbb N}$ tels que
$$r^{'2}+r^{''2}+2N'+2N''=2n,$$
$$r'\equiv r''\equiv\,\,val_{F}(\eta)\,\,mod\,\,2{\mathbb Z} .$$
On pose $N=N'+N''$. On a d\'efini en \cite{MW} 3.11 et 3.12 une application
$${\cal Q}(r',r'')^{Lie}\circ \rho_{N}^*\circ \iota_{N',N''}:{\mathbb C}[\hat{W}_{N'}]_{cusp}\otimes {\mathbb C}[\hat{W}_{N''}]_{cusp}\to C_{cusp}^{\infty}(\mathfrak{g}_{\eta,iso}(F))\oplus C_{cusp}^{\infty}(\mathfrak{g}_{\eta,an}(F)).$$
 
  On fixe des \'el\'ements $w'\in W_{N'}$ et $w''\in W_{N''}$ dont les classes de conjugaison sont param\'etr\'ees par des couples de partitions  de la forme $(\emptyset,\beta')$ et $(\emptyset,\beta'')$. On pose
$$f={\cal Q}(r',r'')^{Lie}\circ\rho_{N}^*\circ\iota_{N',N''}(\varphi_{w'}\otimes \varphi_{w''}).$$

{\bf Remarque.} Dans le cas o\`u $r'=r''=0$, cette fonction est nulle si le couple $(w',w'')$ ne v\'erifie pas la relation $sgn(\eta\varpi^{-val_{F}(\eta)})sgn_{CD}(w')sgn_{CD}(w'')=1$.
\bigskip

Soit $(n_{1},n_{2})\in D(n)$ et soient $\eta_{1},\eta_{2}\in F^{\times}/F^{\times2}$ tels que $\eta=\eta_{1}\eta_{2}$.   A ces objets est associ\'ee une donn\'ee endoscopique de $G_{\eta,iso}$ et $G_{\eta,an}$. Le groupe endoscopique de cette donn\'ee est $G_{n_{1},\eta_{1},iso}\otimes G_{n_{2},\eta_{2},iso}$.  De nouveau, il y a un transfert de $ C_{cusp}^{\infty}(\mathfrak{g}_{\eta,iso}(F))\oplus C_{cusp}^{\infty}(\mathfrak{g}_{\eta,an}(F))$ dans
$$\left( C_{cusp}^{\infty}(\mathfrak{g}_{n_{1},\eta_{1},iso}(F))\oplus C_{cusp}^{\infty}(\mathfrak{g}_{n_{1},\eta_{1},an}(F))\right)\otimes \left( C_{cusp}^{\infty}(\mathfrak{g}_{n_{2},\eta_{2},iso}(F))\oplus C_{cusp}^{\infty}(\mathfrak{g}_{n_{2},\eta_{2},an}(F))\right).$$
On le note $transfert_{n_{1},\eta_{1},n_{2},\eta_{2}}$. 
Il y a un choix naturel de facteur de transfert. Les formules en sont donn\'ees en \cite{W3} proposition X.8. Il convient encore d'en supprimer les discriminants de Weyl. Le $\eta$ de la d\'efinition X.7 de \cite{W3} n'est pas notre pr\'esent $\eta$, c'est $(-1)^n\eta$. On note $\Delta_{n_{1},\eta_{1},n_{2},\eta_{2}}$ ce facteur de transfert. 

Posons
$$t'_{1}=t''_{1}=\frac{r'+r''}{2},\,\, t'_{2}=t''_{2}=\frac{\vert r'-r''\vert }{2}.$$
D\'efinissons deux entiers $n_{1}$ et $n_{2}$ par les formules
$$2n_{1}=t^{'2}_{1}+t^{''2}_{1}+2N',\,\,2n_{2}=t^{'2}_{2}+t^{''2}_{2}+2N'',$$
qui sont \'equivalentes \`a
$$n_{1}=\frac{(r'+r'')^2}{4}+N',\,\, n_{2}=\frac{(r'-r'')^2}{4}+N''.$$
On v\'erifie que $n_{1}+n_{2}=n$. Fixons $\eta_{1},\eta_{2}\in F^{\times}/F^{\times2}$ tels que $\eta_{1}\eta_{2}=\eta$ et
$$(1) \qquad val_{F}(\eta_{1})\equiv t'_{1}=t''_{1}\,,\,mod\,\,2{\mathbb Z},\,\,val_{F}(\eta_{2})\equiv t'_{2}=t''_{2}\,,\,mod\,\,2{\mathbb Z}.$$
  Pour $j=1,2$, on peut  appliquer la m\^eme construction que ci-dessus \`a $n_{j}$ et $\eta_{j}$, o\`u le couple $(N',N'')$ est remplac\'e par $(N',0)$ si $j=1$, par $(N'',0)$ si $j=2$. Cela nous permet de d\'efinir les \'el\'ements
$$f_{1,\eta_{1}}={\cal Q}(t'_{1},t''_{1})^{Lie}\circ\rho^*_{N'}\circ\iota_{N',0}(\varphi_{w'})\in  C_{cusp}^{\infty}(\mathfrak{g}_{n_{1},\eta_{1},iso}(F))\oplus C_{cusp}^{\infty}(\mathfrak{g}_{n_{1},\eta_{1},an}(F)),$$
$$f_{2,\eta_{2}}={\cal Q}(t'_{2},t''_{2})^{Lie}\circ\rho^*_{N''}\circ\iota_{N'',0}(\varphi_{w''})\in  C_{cusp}^{\infty}(\mathfrak{g}_{n_{2},\eta_{2},iso}(F))\oplus C_{cusp}^{\infty}(\mathfrak{g}_{n_{2},\eta_{2},an}(F)).$$

 D\'efinissons une constante $C_{\eta_{1},\eta_{2}}$ par les formules suivantes:

si $r''\leq r'$, $C_{\eta_{1},\eta_{2}}= sgn(\eta_{2}\varpi^{-val_{F}(\eta_{2})})^{val_{F}(\eta)}$;

si $r'<r''$, $C_{\eta_{1},\eta_{2}}=sgn(-1)^{ val_{F}(\eta_{2})}sgn_{CD}(w'') sgn(\eta_{2}\varpi^{-val_{F}(\eta_{2})})^{1+val_{F}(\eta)} $.

\ass{Lemme}{Soient $\eta_{1},\eta_{2}\in F^{\times}/F^{\times2}$ tels que $\eta_{1}\eta_{2}=1$. 

(i) Si (1) est v\'erifi\'e,   on a l'\'egalit\'e $transfert_{n_{1},\eta_{1},n_{2},\eta_{2}}(f)=C_{\eta_{1},\eta_{2}}f_{1,\eta_{1}}\otimes f_{2,\eta_{2}}$.

(ii) Soit $(\bar{n}_{1},\bar{n}_{2})\in D(n)$. Supposons  que $(\bar{n}_{1},\bar{n}_{2})\not=(n_{1},n_{2})$ ou que (1) ne soit pas v\'erifi\'e.  Alors  $transfert_{\bar{n}_{1},\eta_{1},\bar{n}_{2},\eta_{2}}(f)=0$.}

La d\'emonstration de ce lemme occupe les paragraphes 2.4 \`a 2.8. 

\bigskip

\subsection{Le  cas unitaire}

Dans ce paragraphe, on fixe une tour d'extensions  finies non ramifi\'ees $E/E^{\natural}/F$, avec $[E:E^{\natural}]=2$, et un entier $d\geq1$. On consid\`ere les espaces vectoriels $V$ sur $E$, de dimension $d$, munis d'une forme hermitienne non d\'eg\'en\'er\'ee $Q$, relativement \`a l'extension $E/E^{\natural}$. De nouveau, il y a deux classes d'isomorphie de couples $(V,Q)$, qui se distinguent par la valuation du d\'eterminant de $Q$ (ce d\'eterminant est un \'el\'ement de $E^{\natural,\times}/norme_{E/E^{\natural}}(E^{\times })$, sa valuation est l'image de ce terme dans ${\mathbb Z}/2{\mathbb Z}$ par l'application $val_{E^{\natural}}$). On note $(V_{iso},Q_{iso})$ celui pour lequel cette valuation est paire et $(V_{an},Q_{an})$ celui pour lequel cette valuation est impaire. On note $G_{iso}$ et $G_{an}$ les groupes unitaires de ces espaces hermitiens. 

Pour $m\in {\mathbb N}$, on sait que les classes de conjugaison dans $\mathfrak{S}_{m}$ sont param\'etr\'ees par les partitions de $m$.  On note ${\mathbb C}[\hat{\mathfrak{S}}_{m}]_{U-cusp}$ l'espace des fonctions sur $\mathfrak{S}_{m}$, invariantes par conjugaison, \`a support dans des classes de conjugaison param\'etr\'ees par des partitions dont tous les termes non nuls sont impairs.

Soient $(d',d'')\in D(d)$. Dans \cite{MW} 3.1, 3.2 et 3.3, on a d\'efini une application
$${\cal Q}(d',d'')^{Lie}\circ\rho^*_{d}\circ \iota_{d',d''}:{\mathbb C}[\hat{\mathfrak{S}}_{d'}]_{U-cusp}\otimes {\mathbb C}[\hat{\mathfrak{S}}_{d''}]_{U-cusp}\to C_{cusp}^{\infty}(\mathfrak{g}_{iso}(E^{\natural}))\oplus C_{cusp}^{\infty}(\mathfrak{g}_{an}(E^{\natural}))
.$$
Fixons $w'\in \mathfrak{S}_{d'}$ et $w''\in \mathfrak{S}_{d''}$ dont les classes de conjugaison sont param\'etr\'ees par des partitions dont tous les termes non nuls sont impairs. On note $\varphi_{w'}$ et $\varphi_{w''}$ les fonctions caract\'eristiques de ces classes de conjugaison. On  pose
$$f={\cal Q}(d',d'')^{Lie}\circ\rho^*_{d}\circ \iota_{d',d''}(\varphi_{w'}\otimes \varphi_{w''}).$$

Soit $(d_{1},d_{2})\in D(d)$. Ce couple d\'etermine une donn\'ee endoscopique de $G_{iso}$ et $G_{an}$. Le groupe endoscopique de cette donn\'ee est $G_{d_{1},iso}\times G_{d_{2},iso}$. Il y a un transfert de $C_{cusp}^{\infty}(\mathfrak{g}_{iso}(E^{\natural}))\oplus C_{cusp}^{\infty}(\mathfrak{g}_{an}(E^{\natural}))$ dans
$$\left(C_{cusp}^{\infty}(\mathfrak{g}_{d_{1},iso}(E^{\natural}))\oplus C_{cusp}^{\infty}(\mathfrak{g}_{d_{1},an}(E^{\natural}))\right)\otimes \left(C_{cusp}^{\infty}(\mathfrak{g}_{d_{2},iso}(E^{\natural}))\oplus C_{cusp}^{\infty}(\mathfrak{g}_{d_{2},an}(E^{\natural}))\right).$$
On le note $transfert_{d_{1},d_{2}}$. 
Il y a un choix naturel de facteur de transfert, cf. \cite{W3} proposition X.8. Comme pr\'ec\'edemment, on supprime les discriminants de Weyl. Le $\eta$ de la d\'efinition X.7 de \cite{W3} est une unit\'e de $E^{\natural}$ si $d$ est impair, une unit\'e de $E$ de trace nulle dans $E^{\natural}$ si $d$ est pair. 

Dans le cas o\`u $(d_{1},d_{2})=(d',d'')$, on peut appliquer la construction ci-dessus en rempla\c{c}ant $d$ par $d'$, resp. $d''$, et $(d',d'')$ par $(d',0)$, resp. $(d'',0)$. On d\'efinit ainsi
$$f_{1}={\cal Q}(d',0)^{Lie}\circ\rho^*_{d'}\circ \iota_{d',0}(\varphi_{w'}),$$
$$f_{2}={\cal Q}(d'',0)^{Lie}\circ\rho^*_{d''}\circ \iota_{d'',0}(\varphi_{w''}).$$

\ass{Lemme}{(i) On a l'\'egalit\'e $transfert_{d',d''}(f)=f_{1}\otimes f_{2}$.

(ii) Pour un couple $(d_{1},d_{2})\in D(d)$ diff\'erent de $(d',d'')$, on a l'\'egalit\'e $transfert_{d_{1},d_{2}}(f)=0$.}

\bigskip

\subsection{Une expression \`a l'aide de transform\'ees de Fourier d'int\'egrales orbitales}
On commence la d\'emonstration du lemme 2.2. On  utilise les notations de ce paragraphe et les donn\'ees que l'on y a fix\'ees. Rappelons que, pour $\sharp=iso$ ou $an$, on d\'efinit dans l'espace $\mathfrak{g}_{\sharp}(F)$ la notion d'\'el\'ement topologiquement nilpotent. En identifiant $\mathfrak{g}_{\sharp}(F)$ \`a un sous-ensemble de l'alg\`ebre des endomorphismes de $V_{\sharp}$, un \'el\'ement $Y\in \mathfrak{g}_{\sharp}(F)$ est topologiquement nilpotent si et seulement si la suite des puissances $Y^{m}$ tend vers $0$ quand $m$ tend vers l'infini. Par construction de l'application ${\cal Q}(r',r'')^{Lie}$, la fonction $f$ est \`a support topologiquement nilpotent.

On note $(\emptyset,\beta')$ et $(\emptyset,\beta'')$ les couples de partitions param\'etrant les classes de conjugaison de $w'$ et $w''$.  Posons $\beta=\beta'\cup\beta''$, $\beta=(\beta_{1}\geq...\geq \beta_{t}>0)$.   Pour tout $k\in \{1,...,t\}$, consid\'erons la tour d'extensions non ramifi\'ees $E_{k}/E_{k}^{\natural}/F$ telle que $[E_{k}:E_{k}^{\natural}]=2$ et $[E_{k}^{\natural}:F]=\beta_{k}$. On fixe un \'el\'ement $X_{k}\in E_{k}^{\times}$ tel que $val_{E_{k}}(X_{k})=0$,  que  $trace_{E_{k}/E_{k}^{\natural}}(X_{k})=0$ et que, en notant $\bar{X}_{k}$ la r\'eduction de $X_{k}$ dans le corps r\'esiduel ${\mathbb F}_{q^{2\beta_{k}}}$, tous les conjugu\'es de $\bar{X}_{k}$ par le groupe de Galois de l'extension ${\mathbb F}_{q^{2\beta_{k}}}/{\mathbb F}_{q}$ soient distincts. On suppose de plus que, pour $k,k'\in \{1,...,t\}$ avec $k\not=k'$, $\bar{X}_{k}$ n'est pas conjugu\'e \`a $\bar{X}_{k'}$. Cette hypoth\`ese est loisible puisque  $p>6n+4$. 

Posons $R=sup(r',r'')$, $r=inf(r',r'')$, $J=\{1,...,R\}$, $\hat{J}=\{j\in J; j\leq R-r, j\,\,pair\}$. Fixons un ensemble de repr\'esentants $\Gamma_{0}$ de $\mathfrak{o}^{\times}/(1+\varpi\mathfrak{o})$.Pour tout $j\in J$ tel que $j>R-r$, fixons un ensemble de repr\'esentants $\Gamma_{j}$ de $\mathfrak{o}^{\times}/\mathfrak{o}^{\times2}$. Pour tout $j\in J$ tel que $j\leq R-r$,  posons $\Gamma_{j}=\Gamma_{0}$. 
Notons $\Gamma$ le sous-ensemble des $\gamma=(\gamma_{j})_{j\in J}\in \prod_{j\in J}\Gamma_{j}$ tels que

pour tout $j\in \hat{J}$, $\gamma_{j-1}\not\in \gamma_{j}+\varpi\mathfrak{o}$;

(1) $sgn(\eta\varpi^{-val_{F}(\eta)}\prod_{j\in J}\gamma_{j})=sgn_{CD}(w')sgn_{CD}(w'')$.

{\bf Remarques.} L'ensemble $\Gamma$ n'a \'evidemment rien \`a voir avec celui de \cite{W5} 1.8. D'autre part, dans le cas o\`u $r'=r''=0$, on a $J=\emptyset$ et $\prod_{j\in J}\Gamma_{j}$ a un unique \'el\'ement. Cet \'el\'ement v\'erifie (1) si et seulement si $sgn(\eta\varpi^{-val_{F}(\eta)} =sgn_{CD}(w')sgn_{CD}(w'')$. Donc $\Gamma$ est vide si cette \'egalit\'e n'est pas v\'erifi\'ee. Dans l'\'egalit\'e du lemme ci-dessous, le membre de droite est nul. C'est coh\'erent avec la remarque de 2.2 qui nous dit que $f_{\sharp}=0$. 

\bigskip

On note $\sigma:\Gamma\to {\mathbb C}$ la fonction d\'efinie par
$$\sigma(\gamma)=\left(\prod_{j\in \hat{J}}(q-2+sgn(\gamma_{j-1}\gamma_{j}))sgn(\gamma_{j-1}\gamma_{j}(\gamma_{j-1}-\gamma_{j}))\right)$$
$$\prod_{j=R-r+1,...R; \,j\,\, impair}sgn(-\gamma_{j}).$$

Notons $\bar{F}$ une cl\^oture alg\'ebrique de $F$. Pour $\gamma\in \Gamma$ et pour $j\in J$, on fixe $a_{j}(\gamma)\in \bar{F}^{\times}$ tel que

si $j\in \hat{J}$, $a_{j}(\gamma)^{2j-2}=\varpi\gamma_{j}$;

si $j\leq R-r$ et $j$ est impair, $a_{j}(\gamma)^{2j}=\varpi\gamma_{j}$;

si $j> R-r$, $a_{j}(\gamma)^{2(2j-R+r-1)}=\varpi\gamma_{j}$.
On note $F_{j}(\gamma)$ l'extension $F[a_{j}(\gamma)]$, $F_{j}^{\natural}(\gamma)=F[a_{j}(\gamma)^2]$ la sous-extension d'indice $2$, $\mathfrak{p}_{j}(\gamma)$ l'id\'eal maximal de l'anneau des entiers de $F_{j}(\gamma)$ et
$$A_{j}(\gamma)=\{a_{j}\in a_{j}(\gamma)+\mathfrak{p}_{j}(\gamma)^2; trace_{F_{j}(\gamma)/F_{j}^{\natural}(\gamma)}(a_{j})=0\}.$$
On pose $A(\gamma)=\prod_{j\in J}A_{j}(\gamma)$. Remarquons que $A(\gamma)$ est naturellement un espace principal homog\`ene sous un certain groupe. On munit $A(\gamma)$ d'une mesure invariante par l'action de ce groupe et de masse totale $1$.

Consid\'erons des \'el\'ements $\gamma\in \Gamma$ et $a=(a_{j})_{j\in J}\in A(\gamma)$. Consid\'erons des familles $c=(c_{j})_{j\in J}$ et $u=(u_{k})_{k=1,...,t}$ telles que $c_{j}\in F_{j}^{\natural}(\gamma)^{\times}$ pour $j\in J$ et $u_{k}\in E^{\natural,\times}_{k}$ pour $k=1,...,t$. Pour $j\in J$, on construit un espace quadratique $(V_{j},Q_{j})$: $V_{j}=F_{j}(\gamma)$ et, pour $v,v'\in V_{j}$, $Q_{j}(v,v')=[F_{j}(\gamma)/F]^{-1}trace_{F_{j}(\gamma)/F}(\tau_{j}(v)v'c_{j})$, o\`u  $\tau_{j}$ est l'unique \'el\'ement non trivial du groupe de Galois de $F_{j}(\gamma)/F_{j}^{\natural}(\gamma)$. Pour $k=1,...,t$, on construit  un espace quadratique $(V_{k},Q_{k})$: $V_{k}=E_{k}$ et, pour $v,v'\in V_{k}$, $Q_{k}(v,v')=[ E_{k}/F]^{-1}trace_{ E_{k}/F}(\tau_{k}(v)v'u_{k})$, o\`u  $\tau_{k}$ est l'unique \'el\'ement non trivial du groupe de Galois de $ E_{k}/E_{k}^{\natural}$. 
En utilisant la relation $r'\equiv r''\equiv val_{F}(\eta)$ impos\'ee en 2.2 et la relation (1) ci-dessus, on montre que la somme directe des $(V_{j},Q_{j})$ et des $(V_{k},Q_{k})$ est isomorphe \`a l'un des deux espaces   de 2.2, disons \`a $(V_{\sharp},Q_{\sharp})$. On fixe un isomorphisme. On d\'efinit alors un \'el\'ement $X\in \mathfrak{g}_{\sharp}(F)$: il respecte chacun des sous-espaces $V_{j}$ et $V_{k}$; pour $j\in J$, il agit sur $V_{j}$ par multiplication par $a_{j}$ et, pour $k=1,...,t$, il agit sur $V_{k}$ par multiplication par $X_{k}$. La classe de conjugaison de $X$ par le groupe orthogonal $O(Q_{\sharp};F)$ est bien d\'etermin\'ee par nos donn\'ees de d\'epart. Il y a une petite difficult\'e caus\'ee par  la parit\'e de la dimension de $V_{\sharp}$: cette classe de conjugaison par $O(Q_{\sharp};F)$ se d\'ecompose en deux classes de conjugaison par $G_{\sharp}(F)$. Au lieu de $X$, nous fixons donc des \'el\'ements $X^+$ et $X^-$ dans chacune de ces deux classes. On voit facilement que les constructions ne d\'ependent que des images des $c_{j}$ dans les groupes $F_{j}^{\natural}(\gamma)^{\times}/norme_{F_{j}(\gamma)/F_{j}^{\natural}(\gamma)}(F_{j}(\gamma)^{\times})$ et des images des $u_{k}$ dans les groupes $E_{k}^{\natural,\times}/norme_{E_{k}/E_{k}^{\natural}}(E_{k}^{\times})$. Parce que les $F_{j}(\gamma)/F$ sont totalement ramifi\'es et les $E_{k}/F$ sont non ramifi\'es, ces groupes s'identifient respectivement \`a $\mathfrak{o}^{\times}/\mathfrak{o}^{\times2}$ et ${\mathbb Z}/2{\mathbb Z}$. En posant ${\cal E}=(\mathfrak{o}^{\times}/\mathfrak{o}^{\times2})^J$ et ${\cal U}=({\mathbb Z}/2{\mathbb Z})^t$, on peut donc remplacer les donn\'ees $c=(c_{j})_{j\in J}$ et $u=(u_{k})_{k=1,...,t}$ ci-dessus par des familles $c=(c_{j})_{j\in J}\in {\cal E}$ et $u=(u_{k})_{k=1,...,t}\in {\cal U}$. Notons plus pr\'ecis\'ement  $X^+(a,c,u)$ et $X^-(a,c,u)$ les \'el\'ements associ\'es ci-dessus \`a de telles familles.

  Soient  $\gamma\in \Gamma$ et $e\in {\cal E}$. On d\'efinit une nouvelle famille $c [\gamma,e]=(c[\gamma,e]_{j})_{j\in J}\in {\cal E}$ par les \'egalit\'es suivantes, pour $j\in J$:
 
  si $r''\leq r'$, $c[\gamma,e]_{j}=(-1)^j\gamma_{j}e_{j}$;
  
  si $r'<r''$, $c[\gamma,e]_{j}=(-1)^{j+1}e_{j}$
  
\noindent (il s'agit plut\^ot  des images modulo $\mathfrak{o}^{\times2}$ des expressions indiqu\'ees).
   Pour $a\in A(\gamma)$ et $\zeta=\pm$, on note d\'esormais $X^{\zeta}(a,e,u)$ l'\'el\'ement $X^{\zeta}(a,c[\gamma,e],u)$ introduit ci-dessus. Il appartient \`a $\mathfrak{g}_{\sharp}(F)$ pour un certain indice $\sharp$, que l'on note plus pr\'ecisement $\sharp(a,e,u)$.

Soit $\sharp=iso$ ou $an$. De m\^eme qu'en 1.1, on introduit une transformation de Fourier dans $C_{c}^{\infty}(\mathfrak{g}_{\sharp}(F))$. Rappelons un r\'esultat d'Harish-Chandra. Soit $X\in \mathfrak{g}_{\sharp}(F)$ un \'el\'ement r\'egulier et elliptique. Il existe une fonction $Y\mapsto \hat{i}_{\sharp}(X,Y)$,  d\'efinie et localement constante sur l'ensemble des \'el\'ements elliptiques r\'eguliers de $\mathfrak{g}_{\sharp}(F)$ de sorte que, pour tout tel \'el\'ement $Y$ et pour toute $\phi\in C_{c}^{\infty}(\mathfrak{g}_{\sharp}(F))$, on ait l'\'egalit\'e
$$J(X ,\hat{\phi})=\int_{\mathfrak{g}_{\sharp}(F)}\hat{i}_{\sharp} (X,Y)\phi(Y)D^{G_{\sharp}}(Y)^{-1/2}dY.$$
Les \'el\'ements $X^{\zeta}(a,e,u)$ ci-dessus, tels que $\sharp(a,e,u)=\sharp$,  sont elliptiques r\'eguliers dans $\mathfrak{g}_{\sharp}(F)$. On note $\hat{i}^{\zeta}_{\sharp}[a,e,u]$ la fonction $Y\mapsto \hat{i}_{\sharp}(X^{\zeta}(a,e,u),Y)$. On pose
$$\hat{i}_{\sharp}[a,e,u]=\frac{1}{2}(\hat{i}_{\sharp}^+[a,e,,u]+\hat{i}_{\sharp}^-[a,e,u]).$$
Pour un triplet $(a,e,u)$ tel que $\sharp(a,e,u)\not=\sharp$, on pose $\hat{i}_{\sharp}[a,e,u]=0$.

 On fixe une  d\'ecomposition $\{1,...,t\}=K'\sqcup K''$ de sorte que $\beta'$, resp. $\beta''$, soit form\'ee des $\beta_{k}$ pour $k\in K'$, resp. $k\in K''$.  On d\'efinit un caract\`ere $\kappa^{{\cal U}}$ de ${\cal U}$  par
$$\kappa^{{\cal U}}(u)=(-1)^{\sum_{k\in K''}u_{k}}.$$
 On va se limiter  \`a un sous-groupe  ${\cal E}^0$ de ${\cal E}$.  C'est le sous-groupe des familles $e=(e_{j})_{j\in J}\in {\cal E}$ telles que

 pour tout $j\in \hat{J}$, tel que $j<R$,  $e_{j-1}=e_{j}$.
 
 {\bf Remarque.} La condition  $j<R$ est automatique si $r>0$.   
 
 \bigskip 

On d\'efinit le caract\`ere $\kappa^0$ du groupe ${\cal E}^0$ par 
$$\kappa^0(e)=\prod_{j\in \hat{J}}sgn(e_{j-1}).$$

On doit enfin d\'efinir quelques constantes. On note ${\cal O}(w')$ et ${\cal O}(w'')$ les classes de conjugaison de $w'$ dans $W_{N'}$ et de $w''$ dans $W_{N''}$. On pose
$$C(w')=\vert {\cal O}(w')\vert \vert W_{N'}\vert ^{-1}q^{N'/2}\prod_{k\in K'}(q^{\beta_{k}}+1)^{-1},$$
et on d\'efinit $C(w'') $ de fa\c{c}on similaire. On pose

$\alpha(r',r'',w',w'')=sgn((-1)^{(r'+r'')/2}\eta\varpi^{-val_{F}(\eta)})$   si $val_{F}(\eta)$ est pair; 

$\alpha(r',r'',w',w'')=  sgn((-1)^{(r'+r'')/2}\eta\varpi^{-val_{F}(\eta)})sgn_{CD}(w')sgn_{CD}(w'')$ si $val_{F}(\eta)$ est impair;

$$C(r',r'')=2^{1-r'-r''}\left((q-1)^2(q-3)\right)^{(r-R)/2}.$$

\ass{Lemme}{ Pour $\sharp=iso$ ou $an$ et pour tout \'el\'ement  $Y\in \mathfrak{g}_{\sharp}(F)$ qui est elliptique, r\'egulier et topologiquement nilpotent, on a l'\'egalit\'e
$$J(Y,f_{\sharp})=C(w')C(w'')C(r',r'')\alpha(r',r'',w',w'')\sum_{\gamma\in \Gamma}\sigma(\gamma)$$
$$\int_{A(\gamma)}\sum_{e\in {\cal E}^0}\sum_{u\in {\cal U}}\kappa^0(e)\kappa^{{\cal U}}(u)\hat{i}_{\sharp}[a,e,u](Y) \,da.$$}

Preuve. Posons $N=N'+N''$. En \cite{MW} 3.15 (1), on a \'ecrit une \'egalit\'e
$$(2) \qquad f_{\sharp}=\sum_{(n',n'')\in D(N)}\sum_{(v',v'')\in {\cal V}_{n',n'',\sharp}} b(v',v''){\cal Q}(r',r'';v',v'')_{\sharp}^{Lie}.$$
Les ${\cal Q}(r',r'';v',v'')_{\sharp}^{Lie}$ sont des fonctions sur $\mathfrak{g}_{\sharp}(F)$ et les $b(v',v'')$ sont des coefficients. L'ensemble ${\cal V}_{n',n'',\sharp}$ est un ensemble de repr\'esentants des classes de conjugaison dans $W_{n'}\times W_{n''}$  param\'etr\'ees par des couples de partitions de la forme $(\emptyset,\delta')$, $(\emptyset,\delta'')$ et tels que $\delta'\cup \delta''=\beta$. Il y a une restriction (l'hypoth\`ese 3.13(1) de \cite{MW}), que nous levons en posant  
${\cal Q}(r',r'';v',v'')_{\sharp}^{Lie}=0$ si elle n'est pas v\'erifi\'ee (dans \cite{MW}, cette fonction n'est d\'efinie que sous cette hypoth\`ese 3.13(1)).  En supprimant ainsi cette restriction, l'ensemble ${\cal V}_{n',n'',\sharp}$ devient ind\'ependant de l'indice $\sharp$ et nous supprimons cet indice. Il y a une application naturelle 
$$\delta:{\cal U}\to \cup_{(n',n'')\in D(N)}{\cal V}_{n',n''}.$$
Pour $u=(u_{k})_{k=1,...,t}\in {\cal U}$,  on note $\delta'(u)$ la partition form\'ee des $\beta_{k}$ pour $k$ tel que $u_{k}=0$  et $\delta''(u)$ celle form\'ee des $\beta_{k}$ pour $k$ tel que $u(k)=1$. L'application ci-dessus envoie $u$ sur le couple de classes de conjugaison param\'etr\'ees par $(\emptyset,\delta'(u))$ et $(\emptyset,\delta''(u))$. L'\'egalit\'e 3.15(6) de \cite{MW}  et les calculs qui la suivent montrent que, pour $(v',v'')\in   \cup_{(n',n'')\in D(N)}{\cal V}_{n',n''}$, on a l'\'egalit\'e
$$(3) \qquad b(v',v'')=C_{1}\sum_{u\in \delta^{-1}(v',v'')}\kappa^{{\cal U}}(u),$$
o\`u 
$$(4) \qquad C_{1}=\vert {\cal O}(w')\vert \vert W_{N'}\vert ^{-1}\vert {\cal O}(w'')\vert \vert W_{N''}\vert ^{-1}.$$
Fixons $(v',v'')\in  \cup_{(n',n'')\in D(N)}{\cal V}_{n',n''}$. L'int\'egrale orbitale $J(Y,{\cal Q}(r',r'';v',v'')_{\sharp}^{Lie})$ est calcul\'ee par une application successive des propositions 3.13 et 3.14 de \cite{MW}, o\`u l'on remplace $(w',w'')$ par $(v',v'')$. On obtient
$$(5) \qquad J(Y,{\cal Q}(r',r'';v',v'')_{\sharp}^{Lie})=C_{2}(r',r'',v',v'')\sum_{\gamma\in \Gamma}\sigma(\gamma)$$
$$\int_{A(\gamma)}\sum_{e\in {\cal E}^0_{MW}}\kappa^0_{MW}(e)\hat{i}_{\sharp,MW}[a,e,v',v''](Y)\,da,$$
o\`u 

$C_{2}(r',r'',v',v'')$ est une certaine constante sur laquelle nous allons revenir;

$\Gamma$ et, pour $\gamma\in \Gamma$, $A(\gamma)$ et $\sigma(\gamma)$ sont les termes que l'on a d\'efinis ci-dessus;

les autres termes ${\cal E}^0_{MW}$ etc...  sont quelque peu diff\'erents des n\^otres, on y a ajout\'e un indice $MW$ pour les en distinguer.  

 On constate que le couple $(v',v'')$ intervient dans  (5) d'une part par la constante, d'autre part par la fonction  $\hat{i}_{\sharp,MW}[a,e,v',v'']$.  Pour ce qui est de la constante, le couple $(v',v'')$ intervient via des termes $\vert T'\vert \vert T''\vert $ et des produits $sgn_{CD}(v')sgn_{CD}(v'')$, cf. les formules de \cite{MW}. Il r\'esulte de la d\'efinition ci-dessus que ce dernier produit vaut $sgn_{CD}(w')sgn_{CD}(w'')$. En se reportant \`a la d\'efinition des tores $T'$ et $T''$, on calcule explicitement
 $$(6) \qquad \vert T'\vert \vert T''\vert =\prod_{k=1,...,t}(q^{\beta_{k}}+1).$$
 En conclusion, $C_{2}(r',r'',v',v'')$  ne d\'epend pas de $(v',v'')$ mais seulement  de $(w',w'')$. On la note plut\^ot $C_{2}(r',r'',w',w'')$. Quant \`a la fonction $\hat{i}_{\sharp,MW}[a,e,v',v'']$,  elle ne d\'epend de $(v',v'')$ que par un \'el\'ement $X$ d\'efini en \cite{MW} 3.13. On a un certain choix pour cet \'el\'ement. On constate que, pour tout $u\in \delta^{-1}(v',v'')$, on peut prendre pour $X$ un \'el\'ement construit comme ci-dessus \`a l'aide des donn\'ees $(X_{k})_{k=1,...,t}$ et $u$. Notons 
$\hat{i}_{\sharp,MW}[a,e,u]$ la fonction associ\'ee \`a ce choix. L'\'egalit\'e suivante  r\'esulte alors  alors de (2), (3) et (5):
$$(7) \qquad J(Y,f_{\sharp})=C_{1}C_{2}(r',r'',w',w'')\sum_{\gamma\in \Gamma}\sigma(\gamma)\int_{A(\gamma)}\sum_{e\in {\cal E}^0_{MW}}\sum_{u\in {\cal U}}\kappa^0_{MW}(e)\kappa^{{\cal U}}(u)\hat{i}_{\sharp,MW}[a,e,u](Y)\,da.$$
 Supposons $(r',r'')\not=(0,0)$. Il y a une application naturelle de notre ensemble ${\cal E}^0$ dans l'ensemble ${\cal E}^0_{MW}$. A $e\in {\cal E}^0$, elle associe l'unique \'el\'ement $e_{MW}\in {\cal E}^0_{MW}$ tel que $e_{MW,j}=e_{j}$ pour $j=1,...,R-1$. L'\'el\'ement $e_{MW,R}$ est d\'efini par les \'egalit\'es
 
 $\prod_{j=R-r+1,...,R}e_{MW,j}=1$ si $r>0$;
 
 $e_{MW,R}=e_{MW,R-1}$ si $r=0$. 
 
 On constate que $\kappa^0(e)=\kappa^0_{MW}(e_{MW})$. 
 Fixons $\gamma\in \Gamma$ et $a\in A(\gamma)$. Dans l'\'enonc\'e du lemme, on peut \'evidemment limiter la somme en $e$ et $u$ \`a la somme sur les couples $(e,u)$ tels que   $\sharp(a,e,u)=\sharp$: si cette condition n'est pas v\'erifi\'ee, $\hat{i}_{\sharp}[a,e,u]=0$. Notons $({\cal E}^0\times{\cal U})_{a,\sharp}$ ce sous-ensemble. On constate alors que l'application
 $$\begin{array}{ccc}({\cal E}^0\times{\cal U})_{a,\sharp}&\to&{\cal E}_{MW}^0\times {\cal U}\\ (e,u)&\mapsto& (e_{MW},u)\\ \end{array}$$
 est bijective et que, pour $(e,u)\in ({\cal E}^0\times {\cal U})_{a,\sharp}$, on a l'\'egalit\'e
 $\hat{i}_{\sharp,MW}[a,e_{MW},u]=\hat{i}_{\sharp}[a,e,u]$.
 En utilisant ces propri\'et\'es, l'\'egalit\'e (7) devient alors l'\'egalit\'e de l'\'enonc\'e, aux constantes pr\`es. Il ne reste donc plus qu'\`a d\'emontrer l'\'egalit\'e
 $$(8) \qquad C_{1}C_{2}(r',r'',w',w'')=C(w')C(w'')C(r',r'')\alpha(r',r'',w',w'').$$
 On a suppos\'e $(r',r'')\not=(0,0)$. Si $(r',r'')=(0,0)$, les termes $\Gamma$, ${\cal E}^0$ etc... disparaissent et la formule (7) se compare directement \`a celle de l'\'enonc\'e. De nouveau, il ne reste plus qu'\`a  comparer les constantes.
 
 La constante $C_{2}(r',r'',w',w'')$ est le produit de la constante $C\alpha(w',w'')_{\sharp}$ de la proposition 3.14 de \cite{MW} et de l'inverse de la constante de la proposition 3.13 de cette r\'ef\'erence. Malheureusement, il y a des erreurs dans ces constantes. On a d\'ej\`a effectu\'e quelques corrections dans \cite{W2} page 335:
 
  la constante de la proposition 3.13 est
  $$(9)\qquad q^{-N/2}2^{-\beta(r',r'')}\gamma(r',r'')_{\sharp}\vert T'\vert \vert T''\vert $$
  (en fait, $T'$ et $T''$ sont associ\'es \`a un couple $(v',v'')\in {\cal V}_{n',n''}$ mais,  d'apr\`es (5), le produit de leurs nombres d'\'el\'ements en est ind\'ependant); 
  
  on a l'\'egalit\'e
  $$C=2^{1-r'-r''-\beta(r',r'')}\left((q-1)^2(q-3)\right)^{(r-R)/2}.$$
  
  A l'aide de (4) et (6), on voit que le produit de $C_{1}$, $C$ et de l'inverse de (9) vaut
  $C(w')C(w'')C(r',r'')\gamma(r',r'')_{\sharp}^{-1}$.  Rappelons qu'en \cite{W5} 1.1, on a associ\'e \`a $Q_{\sharp}$ des discriminants $\eta'(Q_{\sharp})$ et $\eta''(Q_{\sharp})$, notons-les simplement $\eta'_{\sharp}$ et $\eta''_{\sharp}$.  D'apr\`es \cite{MW} 3.13, on a 
  
  $\gamma(r',r'')_{\sharp}=sgn((-1)^{(r'+r'')/2} \eta'_{\sharp}\eta''_{\sharp})$, si $r'$ et $r''$ sont pairs, c'est-\`a-dire si $val_{F}(\eta)$ est pair;
  
  $\gamma(r',r'')_{\sharp}=1$ si $val_{F}(\eta)$ est impair.
  
  Il y a une autre correction \`a apporter \`a \cite{MW}. On a
  
   $\alpha(w',w'')_{\sharp}=1$  si $val_{F}(\eta)$ est pair;
   
   $\alpha(w',w'')_{\sharp}=sgn((-1)^{(r'-r'')/2}\eta'_{\sharp}\eta''_{\sharp})sgn_{CD}(w')sgn_{CD}(w'')$ si $val_{F}(\eta)$ est impair. 
   
   Cette constante provient en fait de \cite{W3} paragraphes VII.25 et VII.26. Dans cette r\'ef\'erence, on avait suppos\'e $r'$ et $r''$ pairs et obtenu la constante $1$. Dans \cite{MW}, on a abusivement adopt\'e cette valeur m\^eme quand $r'$ et $r''$ sont impairs. En reprenant la preuve de \cite{W3} VII.26 et en y supposant $r'$ et $r''$ impairs, on obtient la constante ci-dessus. Signalons que, dans le cas sp\'ecial orthogonal impair, il n'y a (\`a notre connaissance) pas de correction \`a apporter \`a la formule de \cite{MW} pour la constante $\alpha(w',w'')_{\sharp}$.

\bigskip

D'apr\`es \cite{W5} 1.1, on a l'\'egalit\'e $\eta'_{\sharp}\eta''_{\sharp}(-1)^{r''}=\eta\varpi^{-val_{F}(\eta)}$. En la reportant dans la formule ci-dessus, on obtient l'\'egalit\'e
$$\gamma(r',r'')_{\sharp}^{-1}\alpha(w',w'')_{\sharp}=\alpha(r',r'',w',w'').$$
En mettant ces calculs bout-\`a-bout, on obtient l'\'egalit\'e (8), ce qui ach\`eve la d\'emonstration. $\square$

\bigskip

\subsection{Transformation de l'expression pr\'ec\'edente}
On note ${\cal L}$ l'ensemble des couples $(L_{1},L_{2})$ de sous-ensembles de $\{1,...,R-r\}$ tels que

$\{1,...,R-r\} $ est r\'eunion disjointe de $L_{1}$ et $L_{2}$;

pour tout $j\in \hat{J}$, les intersections $L_{1}\cap \{j-1,j\}$ et $L_{2}\cap\{j-1,j\}$  ont un unique \'el\'ement; on note cet \'el\'ement $l_{1}(j/2)$, resp. $l_{2}(j/2)$;

si $r=0$ et $R>0$, $l_{2}(R/2)=R-1$.

Pour un tel couple, on d\'efinit un caract\`ere $\kappa^{L_{2}}$ de ${\cal E}$ par
$$\kappa^{L_{2}}(e)=\prod_{j=1,...,(R-r)/2}sgn(e_{l_{2}(j)}).$$

\ass{Lemme}{ Pour $\sharp=iso$ ou $an$ et pour tout \'el\'ement  $Y\in \mathfrak{g}_{\sharp}(F)$ qui est elliptique, r\'egulier et topologiquement nilpotent, on a l'\'egalit\'e
$$J(Y,f_{\sharp})=C(w')C(w'')C(r',r'')\alpha(r',r'',w',w'')\vert {\cal L}\vert ^{-1}\sum_{(L_{1},L_{2})\in {\cal L}}\sum_{\gamma\in \Gamma}\sigma(\gamma)$$
$$\int_{A(\gamma)}\sum_{e\in {\cal E}}\sum_{u\in {\cal U}}\kappa^{L_{2}}(e)\kappa^{{\cal U}}(u)\hat{i}_{\sharp}[a,e,u](Y) \,da.$$}

Preuve. En vertu de l'\'enonc\'e pr\'ec\'edent, il suffit de prouver que, pour $e\in {\cal E}$, on a les \'egalit\'es

(1) $\sum_{(L_{1},L_{2})\in {\cal L}}\kappa^{L_{2}}(e)=0$ si $e\not\in {\cal E}^0$;

(2) $\sum_{(L_{1},L_{2})\in {\cal L}}\kappa^{L_{2}}(e)=\vert {\cal L}\vert \kappa^0(e)$, si $e\in {\cal E}^0$.

Supposons $e\not\in {\cal E}^0$. Alors il existe $j\in  \{1,...,(R-r)/2\}$, avec $j<R/2$ si $r=0$ et $R>0$, de sorte que $e_{2j-1}\not=e_{2j}$. On fixe un tel $j$. On regroupe les \'el\'ements de ${\cal L}$ en paires. Si $(L_{1},L_{2})$ est un \'el\'ement d'une paire, l'autre \'el\'ement $(L'_{1},L'_{2})$ est obtenu en \'echangeant $l_{1}(j)$
 et $l_{2}(j)$ et en ne changeant pas les autres $l_{1}(h)$, $l_{2}(h)$ pour $h\not=j$. On a alors
$$\kappa^{L_{2}}(e)=sgn(e_{l_{2}(j)})\prod_{h\not=j}sgn(e_{l_{2}(h)}),$$
$$\kappa^{L'_{2}}(e)=sgn(e_{l_{1}(j)})\prod_{h\not=j}sgn(e_{l_{2}(h)}).$$
Puisque $e_{l_{2}(j)}\not=e_{l_{1}(j)}$, la somme de ces termes est nulle, ce qui prouve (1). 

Supposons $e\in {\cal E}^0$. Pour tout $j=1,...,(R-r)/2$ et tout $(L_{1},L_{2})\in {\cal L}$, on a l'\'egalit\'e $sgn(e_{l_{2}(j)})=sgn(e_{2j-1})$. Cela r\'esulte de l'\'egalit\'e $e_{2j-1}=e_{2j}$ impos\'ee aux \'el\'ements de ${\cal E}^0$, sauf dans le cas o\`u $r=0$, $R>0$ et $j=R/2$. Mais dans ce cas, $l_{2}(j)=2j-1$ par d\'efinition de ${\cal L}$. En cons\'equence, $\kappa^{L_{2}}(e)=\prod_{j=1,...,(R-r)/2}sgn(e_{2j-1})=\kappa^0(e)$. D'o\`u (2). $\square$

\bigskip

\subsection{Calcul de facteurs de transfert}
Pour $\gamma\in \Gamma$, $a\in A(\gamma)$, $e\in {\cal E}$, $u\in {\cal U}$ et $\zeta=\pm$, on a introduit un \'el\'ement $X^{\zeta}(a,e,u)$. Fixons $\gamma$ et $a$. Quand $e$, $u$ et $\zeta$ varient, ces \'el\'ements $X^{\zeta}(a,e,u)$ d\'ecrivent un ensemble de repr\'esentants des classes de conjugaison contenues dans deux classes totales de conjugaison stable  dans $\mathfrak{g}_{iso}(F)\cup \mathfrak{g}_{an}(F)$. Ces deux classes se d\'eduisent l'une de l'autre par conjugaison par  des \'el\'ements de d\'eterminant $-1$ des groupes orthogonaux. 
On n'a pas donn\'e de crit\`ere pour distinguer $X^+(a,e,u)$ de $X^-(a,e,u)$. On voit que l'on peut effectuer ces choix de signes de sorte que, pour $\zeta$ fix\'e et quand $e$ et $u$ varient, les \'el\'ements $X^{\zeta}(a,e,u)$ d\'ecrivent un ensemble de repr\'esentants des classes de conjugaison contenues dans l'une de nos deux classes totales de conjugaison stable. 
 On fixe  un \'el\'ement   dans cette classe  qui appartienne \`a  $\mathfrak{g}_{iso}(F)$, on le note $X^{\zeta}(a,w',w'')$ .

On a d\'efini $t'_{1}$, $t''_{1}$, $t'_{2}$, $t''_{2}$ en 2.2. On a les \'egalit\'es $t'_{1}=t''_{1}=(R+r)/2$ et $t'_{2}=t''_{2}=(R-r)/2$. Fixons $(L_{1},L_{2})\in {\cal L}$.  
Soit  $\gamma\in \Gamma$. On d\'efinit deux suites $\gamma^{L_1}=(\gamma^{L_1}_{j})_{j=1,...,t'_{1}}$ et $\gamma^{L_2}=(\gamma^{L_2}_{j})_{j=1,...,t'_{2}}$ par les formules  suivantes:

pour $j\in\{1,...,R-r\} $, $\gamma^{L_1}_{j}=\gamma_{l_{1}(j)}$, $\gamma^{L_2}_{j}=\gamma_{l_{2}(j)}$;

pour $j\in \{(R-r)/2+1,...,(R+r)/2\}$, $\gamma^{L_1}_{j}=\gamma_{j+(R-r)/2}$. 

Pour $a\in A(\gamma)$, on d\'efinit de m\^eme des suites $a^{L_1}=(a^{L_1}_{j})_{j=1,...,t'_{1}}$ et $a^{L_2}=(a^{L_2}_{j})_{j=1,...,t'_{2}}$. On d\'efinit un \'el\'ement $\eta[L_{2},\gamma]\in F^{\times}/F^{\times2}$ par les deux propri\'et\'es:

$val_{F}(\eta[L_{2},\gamma])\equiv t'_{2}\,\,mod\,\,2{\mathbb Z}$;

$sgn(\eta[L_{2},\gamma]\varpi^{-val_{F}(\eta[L_{2},\gamma])}\prod_{j=1,...,t'_{2}}\gamma^{L_2}_{j})=sgn_{CD}(w'')$.

On d\'efinit $\eta[L_{1},\gamma]$ par l'\'egalit\'e $\eta[L_{1},\gamma]\eta[L_{2},\gamma]=\eta$. On constate que les donn\'ees $\gamma^{L_1}$, $a^{L_1}$, $(N',0)$, $w'_{1}=w'$ et $w''_{1}=\emptyset$ v\'erifient les m\^emes conditions que $\gamma$, $a$, $(N',N'')$, $w'$ et $w''$,  le couple $(n,\eta)$ \'etant remplac\'e par $(n_{1},\eta[L_{1},\gamma])$. De m\^eme pour les donn\'ees $\gamma^{L_2}$, $a^{L_2}$, $(N'',0)$, $w'_{2}=w''$ et $w''_{2}=\emptyset$,  le couple $(n,\eta)$ \'etant remplac\'e par $(n_{2},\eta[L_{2},\gamma])$.
De m\^eme que l'on a construit ci-dessus des \'el\'ements $X^{\zeta}(a,w',w'')$, on construit dans $\mathfrak{g}_{n_{1},\eta[L_{1},\gamma],iso}(F)$ un \'el\'ement $X^{\zeta}(a^{L_1},w')$ et dans $\mathfrak{g}_{n_{2},\eta[L_{2},\gamma],iso}(F)$ un \'el\'ement $X^{\zeta}(a^{L_2},w'')$.  La correspondance endoscopique entre classes de conjugaison stable est un  peu perturb\'ee par les signes $\zeta$. Pr\'ecis\'ement, pour $\zeta_{1},\zeta_{2}=\pm$, il existe $\zeta=\pm$ tel que la classe de conjugaison stable de $(X^{\zeta_{1}}(a^{L_1},w'),X^{\zeta_{2}}(a^{L_2},w''))$ dans $\mathfrak{g}_{n_{1},\eta[L_{1},\gamma],iso}(F)\times\mathfrak{g}_{n_{2},\eta[L_{2},\gamma],iso}(F)$ corresponde \`a la classe totale de conjugaison stable de $X^{\zeta}(a,w',w'')$. Ainsi, pour $e\in {\cal E}$ et $u\in {\cal U}$, le facteur de transfert $\Delta_{n_{1},\eta[L_{1},\gamma],n_{2},\eta[L_{2},\gamma]} ((X^{\zeta_{1}}(a^{L_1},w'),X^{\zeta_{2}}(a^{L_2},w'')),X^{\zeta}(a,e,u))$ est d\'efini et non nul. 

{\bf Remarque.} Il peut arriver que $n_{1}=0$ ou $n_{2}=0$. Dans ce cas, les termes  $X^{\zeta_{1}}(a^{L_1},w')$ ou $X^{\zeta_{2}}(a^{L_2},w'')$ disparaissent. 

\bigskip

Posons
$$d(r',r'',\gamma,L_{2})=sgn(\eta\varpi^{-val_{F}(\eta)})^{(R-r)/2}sgn_{CD}(w')^{(R-r)/2}sgn_{CD}(w'')^{(R-r)/2+val_{F}(\eta)}$$
$$\left(\prod_{j\in \hat{J}}sgn(\gamma_{j-1}\gamma_{j})^{j/2-1}sgn(\gamma_{j-1}-\gamma_{j})\right)\left(\prod_{j=R-r+1,...,R}sgn(\gamma_{j})^{(R-r)/2}\right)$$
$$sgn(-1)^{val_{F}(\eta[L_{2},\gamma])B}sgn(\eta[L_{2},\gamma]\varpi^{-val_{F}(\eta[L_{2},\gamma])})^Bsgn_{CD}(w'')^B,$$
o\`u $B=0$ si $r'\geq r''$, $B=1$ si $r'<r''$. 

\ass{Lemme}{Sous les hypoth\`eses ci-dessus, on a l'\'egalit\'e
$$\Delta_{n_{1},\eta[L_{1},\gamma],n_{2},\eta[L_{2},\gamma]}  ((X^{\zeta_{1}}(a^{L_1},w'),X^{\zeta_{2}}(a^{L_2},w'')),X^{\zeta}(a,e,u))=d(r',r'',\gamma,L_{2})\kappa^{L_{2}}(e)\kappa^{{\cal U}}(u).$$
}

Preuve. Pour tout polyn\^ome $Q$, on note $Q'$ son polyn\^ome d\'eriv\'e. Notons $P$ le polyn\^ome caract\'eristique de $X^{\zeta}(a,e,u)$, vu comme endomorphisme de $V_{\sharp(a,e,u)}$.   Pour $j\in \{1,...,t'_{2}\}$, posons
$$C_{l_{2}(j)}=(-1)^n\eta[F_{l_{2}(j)}(\gamma):F]c[\gamma,e]_{l_{2}(j)}a_{l_{2}(j)}^{-1}P'(a_{l_{2}(j)}),$$
et notons $sgn_{l_{2}(j)}$ le caract\`ere quadratique de $F^{\natural}_{l_{2}(j)}(\gamma)^{\times}$ associ\'e \`a l'extension quadratique  $F_{l_{2}(j)}(\gamma)$. Pour $k\in K''$, posons
$$C_{k}=(-1)^n\eta[E_{k}:F] \varpi^{u_{k}}X_{k}^{-1}P'(X_{k}),$$
et notons $sgn_{k}$ le caract\`ere quadratique de $E^{\natural}_{k}$ associ\'e \`a l'extension quadratique $E_{k}$.

{\bf Remarque.} Les formules ci-dessus ne d\'efinissent $C_{l_{2}(j)}$ et $C_{k}$ que modulo les groupes 
$norme_{F_{l_{2}(j)}(\gamma)/F^{\natural}_{l_{2}(j)}(\gamma)}(F_{l_{2}(j)}(\gamma)^{\times})$, resp. 
$norme_{E_{k}/E_{k}^{\natural}}(
 E^{\times}_{k})$. C'est sans importance pour la suite. 

\bigskip

 D'apr\`es \cite{W3} proposition X.8, on a l'\'egalit\'e
$$(1) \qquad \Delta_{n_{1},\eta[L_{1},\gamma],n_{2},\eta[L_{2},\gamma]}  ((X^{\zeta_{1}}(a^{L_1},w'),X^{\zeta_{2}}(a^{L_2},w'')),X^{\zeta}(a,e,u))=$$
$$\left(\prod_{j=1,...,t'_{2}}sgn_{l_{2}(j)}(C_{l_{2}(j)})\right)\left(\prod_{k\in K''}sgn_{k}(C_{k})\right).$$

 Pour $k\in K''$, le caract\`ere $sgn_{k}$ est non ramifi\'e. Pour calculer $sgn_{k}(C_{k})$, il suffit de calculer la valuation de $C_{k}$ modulo $2{\mathbb Z}$. On voit facilement que la valuation de $(-1)^n[E_{k}:F]  X_{k}^{-1}P'(X_{k})$ est nulle. Celle de $\eta\varpi{u_{k}}$ est $val_{F}(\eta)+u_{k}$. D'o\`u $sgn_{k}(C_{k})=(-1)^{val_{F}(\eta)+u_{k}}$. Le produit des $(-1)^{u_{k}}$ pour $k\in K''$ est \'egal \`a $\kappa^{{\cal U}}(u)$. Le produit des $(-1)^{val_{F}(\eta)}$ est $(-1)^{\vert K''\vert val_{F}(\eta)}$. Mais $(-1)^{\vert K''\vert }=sgn_{CD}(w'')$. D'o\`u
 $$(2) \qquad \prod_{k\in K''}sgn_{k}(C_{k})=\kappa^{{\cal U}}(u)sgn_{CD}(w'')^{val_{F}(\eta)}.$$

Fixons $j=1,...,t'_{2}$. Posons pour simplifier $l_{1}=l_{1}(j)$, $l_{2}=l_{2}(j)$ et notons $l$ l'\'el\'ement de $\hat{J}$ tel que $\{l_{1},l_{2}\}=\{l-1,l\}$. Le polyn\^ome $P$ se d\'ecompose en produit des polyn\^omes caract\'eristiques $P_{k}$ des $X_{k}$, pour $k=1,...,t$, et des polyn\^omes caract\'eristiques $P_{h}$ des $a_{h}$ pour $h=1,...,R$. On a
$$P'(a_{l_{2}})=P'_{l_{2}}(a_{l_{2}})\left(\prod_{h=1,...,R; h\not=l_{2}}P_{h}(a_{l_{2}})\right)\left(\prod_{k=1,...,t}P_{k}(a_{l_{2}})\right).$$
On v\'erifie que tous les termes ci-dessus sauf le premier appartiennent \`a $F^{\natural}_{l_{2}}(\gamma)^{\times}$. Pour le premier,  son produit avec  $a_{l_{2}}$ appartient \`a ce groupe. On a calcul\'e les images par $sgn_{l_{2}}$ de presque tous ces termes en \cite{W3} page 400. Pour $k=1,...,t$, on a $sgn_{l_{2}}(P_{k}(a_{l_{2}}))=-sgn(-1)^{[E_{k}:F]/2}$. Parce que $[E_{k}:F]/2=\beta_{k}$, on a $\sum_{k=1,...,t}[E_{k}:F]/2=N$. On a aussi $(-1)^t=sgn_{CD}(w')sgn_{CD}(w'')$. On obtient 
$$sgn_{l_{2}}(\prod_{k=1,...,t}P_{k}(a_{l_{2}}))=sgn(-1)^Nsgn_{CD}(w')sgn_{CD}(w'').$$
Soit $h\in \{1,...,l-2\}$. On a $sgn_{l_{2}}(P_{h}(a_{l_{2}}))=sgn(-1)^{[F_{h}(\gamma):F]/2}$. Par construction, tous les $[F_{h}(\gamma):F]$ sont impairs et le terme ci-dessus vaut $sgn(-1)$. Puisque $l$ est pair, le produit de ces termes sur $h=1,...,l-2$ vaut  $1$. Pour $h\in \{l+1,...,R\}$, on a  $sgn_{l_{2}}(P_{h}(a_{l_{2}}))=sgn(-1)^{1+[F_{l_{2}}(\gamma):F]/2}sgn(\gamma_{l_{2}}\gamma_{h})$. Comme ci-dessus, $[F_{l_{2}}(\gamma):F]/2$ est impair et le terme ci-dessus vaut $sgn(\gamma_{l_{2}}\gamma_{h})$. Puisque $l$ est pair, le produit de ces termes sur $h=l+1,...,R$ vaut $sgn(\gamma_{l_{2}})^R\prod_{h=l+1,...,,R}sgn(\gamma_{h})$. On a $sgn_{l_{2}}(a_{l_{2}}P'_{l_{2}}(a_{l_{2}}))=sgn((-1)^{[F_{l_{2}}/F]/2}[F_{l_{2}}:F])$, ou encore $sgn(-[F_{l_{2}}:F])$. Il reste un  dernier terme $P_{l_{1}}(a_{l_{2}})$, dont l'image par $sgn_{l_{2}}$ n'est pas calcul\'ee dans \cite{W3}. Le polyn\^ome caract\'eristique de $a_{l_{1}}$ est $P_{l_{1}}(Y)=Y^{2l-2}-\varpi\gamma_{l_{1}}$. Donc $P_{l_{1}}(a_{l_{2}})=a_{l_{2}}^{2l-2}-\varpi\gamma_{l_{1}}=\varpi(\gamma_{l_{2}}-\gamma_{l_{1}})$. Mais $-\varpi^{-1}\gamma_{l_{2}}$ est une norme de l'extension $F_{l_{2}}(\gamma)/F_{l_{2}}^{\natural}(\gamma)$ (c'est la norme de $\varpi^{-1}a_{l_{2}}^{l-1}$). On a donc 
$$sgn_{l_{2}}(P_{l_{1}}(a_{l_{2}}))=sgn_{l_{2}}(-\varpi^{-1}\gamma_{l_{2}}\varpi(\gamma_{l_{2}}-\gamma_{l_{1}}))=sgn_{l_{2}}(-\gamma_{l_{2}}(\gamma_{l_{2}}-\gamma_{l_{1}}))=sgn(-\gamma_{l_{2}}(\gamma_{l_{2}}-\gamma_{l_{1}})),$$
puisque $sgn_{l_{2}}$ co\"{\i}ncide avec $sgn$ sur $\mathfrak{o}^{\times}$. 
En rassemblant ces calculs, on obtient
$$sgn_{a_{l_{2}}}(a_{l_{2}}P'(a_{l_{2}}))=sgn(-1)^Nsgn_{CD}(w')sgn_{CD}(w'')sgn([F_{l_{2}}(\gamma):F])sgn(\gamma_{l_{2}}-\gamma_{l_{1}})$$
$$sgn(\gamma_{l_{2}})^{R+1}\prod_{h=l+1,...,,R}sgn(\gamma_{h}).$$
Le terme $C_{l_{2}}$ est le produit de $a_{l_{2}}P'(a_{l_{2}})$ et de 
$(-1)^n\eta[F_{l_{2}}(\gamma):F]c[\gamma,e]_{l_{2}}a_{l_{2}}^{-2}$. Le terme $a_{l_{2}}^{-2}$ se remplace par $-1$ car $-a_{l_{2}}^{-2}$ est la norme de $a_{l_{2}}^{-1}$. Le facteur $[F_{l_{2}}(\gamma):F]$ compense celui intervenant dans la formule pr\'ec\'edente. On a 
$$sgn_{l_{2}}(\eta)=sgn_{l_{2}}(\eta\varpi^{-val_{F}(\eta)})sgn_{l_{2}}(\varpi^{val_{F}(\eta)}).$$
Le premier facteur se remplace par $sgn(\eta\varpi^{-val_{F}(\eta)})$.  Puisque 
$-\varpi^{-1}\gamma_{l_{2}}$ est une norme, le deuxi\`eme facteur se remplace par $sgn(-\gamma_{l_{2}})^{val_{F}(\eta)}$, ou encore, puisque $val_{F}(\eta)$ est de la m\^eme parit\'e que $R$, par $sgn(-\gamma_{l_{2}})^R$. En se reportant \`a la d\'efinition de $c[\gamma,e]_{l_{2}}$, on peut \'ecrire $c[\gamma,e]_{l_{2}}=(-1)^{l_{2}}\gamma_{l_{2}}e_{l_{2}}(-\gamma_{l_{2}})^{B}$, o\`u $B $ a \'et\'e d\'efini avant l'\'enonc\'e. 
 D'o\`u 
$$sgn_{l_{2}}(c[\gamma,e]_{l_{2}})=sgn((-1)^{l_{2}}\gamma_{l_{2}}e_{l_{2}})sgn(\gamma_{l_{2}})^B.$$
On obtient alors
$$sgn_{l_{2}}(C_{l_{2}})=sgn(-1)^{N+n+R+1+l_{2}}sgn(\eta\varpi^{-val_{F}(\eta)})sgn_{CD}(w')sgn_{CD}(w'') sgn(e_{l_{2}}) sgn(\gamma_{l_{2}}-\gamma_{l_{1}})$$
$$ sgn(-\gamma_{l_{2}})^B \prod_{h=l+1,...,,R}sgn(\gamma_{h}).$$
On a $n-N=(r^{'2}+r^{''2})/2$.  Puisque $r'$ et $r''$ sont de m\^eme parit\'e, on constate que $n-N$  est de la m\^eme parit\'e que $r'$ et $r''$, ou encore de $R$. Le premier terme se simplifie donc en $sgn(-1)^{1+l_{2}}$. On constate aussi que $(-1)^{1+l_{2}}(\gamma_{l_{2}}-\gamma_{l_{1}})=\gamma_{l-1}-\gamma_{l}$. D'o\`u 
$$(3) \qquad sgn_{l_{2}}(C_{l_{2}})= sgn(\eta\varpi^{-val_{F}(\eta)})sgn_{CD}(w')sgn_{CD}(w'') sgn(e_{l_{2}}) sgn(\gamma_{l-1}-\gamma_{l})$$ 
$$ sgn(-\gamma_{l_{2}})^B \prod_{h=l+1,...,,R}sgn(\gamma_{h}).$$

R\'etablissons l'indice $j$ et faisons le produit de ces expressions sur $j=1,...,t'_{2}=(R-r)/2$. Les premiers termes donnent $sgn(\eta\varpi^{-val_{F}(\eta)})^{(R-r)/2}sgn_{CD}(w')^{(R-r)/2}sgn_{CD}(w'')^{(R-r)/2}$. Le produit des $sgn(e_{l_{2}(j)})$ est \'egal \`a $\kappa^{L_{2}}(e)$. Le produit des termes suivants est $\prod_{j\in \hat{J}}sgn(\gamma_{j-1}-\gamma_{j})$. On a
$$\prod_{j=1,...,t'_{2}}sgn(-\gamma_{l_{2}(j)})=sgn(-1)^{t'_{2}}sgn(\prod_{j=1,...,t'_{2}}\gamma^{L_2}_{j}).$$
Puisque $t'_{2}$ est de la m\^eme parit\'e que $val_{F}(\eta[L_{2},\gamma])$, le premier terme vaut $sgn(-1)^{val_{F}(\eta[L_{2},\gamma])}$. Par d\'efinition de $\eta[L_{2},\gamma]$, le deuxi\`eme terme vaut $sgn(\eta[L_{2},\gamma]\varpi^{-val_{F}(\eta[L_{2},\gamma])})sgn_{CD}(w'')$. Consid\'erons le dernier produit de la formule (3). Pour $h\in \{R-r+1,...,R\}$, le terme $sgn(\gamma_{h})$ intervient pour tout $j$. Le produit en $j$ donne $sgn(\gamma_{h})^{(R-r)/2}$.
Pour $h\in \hat{J}$, les termes $sgn(\gamma_{h-1})$ et $sgn(\gamma_{h})$ interviennent pour les $j< h/2$. Le produit donne $sgn(\gamma_{h-1}\gamma_{h})^{h/2-1}$. D'o\`u
$$(4)\qquad \prod_{j=1,...,t'_{2}}sgn_{l_{2}(j)}(C_{l_{2}(j)})=\kappa^{L_{2}}(e)\left(sgn(\eta\varpi^{-val_{F}(\eta)})sgn_{CD}(w')sgn_{CD}(w'')\right)^{(R-r)/2}$$
$$\left(\prod_{j\in \hat{J}}sgn(\gamma_{j-1}\gamma_{j})^{j/2-1}sgn(\gamma_{j-1}-\gamma_{j}))\right)\left(\prod_{j=R-r+1,...,R}sgn(\gamma_{j})^{(R-r)/2}\right)$$
$$sgn(-1)^{val_{F}(\eta[L_{2},\gamma])B}sgn(\eta[L_{2},\gamma]\varpi^{-val_{F}(\eta[L_{2},\gamma])})^Bsgn_{CD}(w'')^B.$$
Le lemme r\'esulte de (1), (2) et (4). $\square$

\bigskip

\subsection{ D\'emonstration du (ii) du lemme 2.2}
Soient $(\bar{n}_{1},\bar{n}_{2})\in D(n)$ et $\eta_{1},\eta_{2}\in F^{\times}/F^{\times2}$ tels que $\eta_{1}\eta_{2}=\eta$. Soient $\bar{Y}_{1}$, resp. $\bar{Y}_{2}$ un \'el\'ement r\'egulier, elliptique et  topologiquement nilpotent de $\mathfrak{g}_{\bar{n}_{1},\eta_{1}}(F)$, resp. $\mathfrak{g}_{\bar{n}_{2},\eta_{2}}(F)$.  On va calculer  l'int\'egrale orbitale endoscopique $J^{endo}(\bar{Y}_{1},\bar{Y}_{2},f)$, cf. 1.2 dont les d\'efinitions s'adaptent aux alg\`ebres de Lie.  

Soit $\sharp=iso$ ou $an$.
 Pour $(L_{1},L_{2})\in {\cal L}$, $\gamma\in \Gamma$ et  $a\in A(\gamma)$, posons
$$(1) \qquad E_{\sharp}[a,L_{1},L_{2}](\bar{Y}_{1},\bar{Y}_{2})=d(r',r'',\gamma,L_{2})\sum_{Y}\Delta_{\bar{n}_{1},\eta_{1},\bar{n}_{2},\eta_{2}}((\bar{Y}_{1},\bar{Y}_{2}),Y)$$
$$\sum_{e\in {\cal E}}\sum_{u\in {\cal U}}\kappa^{L_{2}}(e)\kappa^{\cal U}(u)\hat{i}_{\sharp}[a,e,u](Y).$$
En utilisant le lemme 2.5, on calcule
 $$  J^{endo}(\bar{Y}_{1},\bar{Y}_{2},f_{\sharp})=\vert {\cal L}\vert ^{-1}\sum_{(L_{1},L_{2})\in {\cal L}}\sum_{\gamma\in \Gamma} C(\gamma,L_{2})\int_{A(\gamma)} E_{\sharp}[a,L_{1},L_{2}](\bar{Y}_{1},\bar{Y}_{2})\,da,$$
o\`u on a pos\'e
$$C(\gamma,L_{2}) =d(r',r'',\gamma,L_{2})C(w')C(w'')C(r',r'')\alpha(r',r'',w',w'')\sigma(\gamma).$$
 Par d\'efinition, on a l'\'egalit\'e
$$ J^{endo}(\bar{Y}_{1},\bar{Y}_{2},f)=J^{endo}(\bar{Y}_{1},\bar{Y}_{2},f_{iso})+J^{endo}(\bar{Y}_{1},\bar{Y}_{2},f_{an}).$$
D'o\`u l'\'egalit\'e
$$(2)\qquad J^{endo}(\bar{Y}_{1},\bar{Y}_{2},f)=\vert {\cal L}\vert ^{-1}\sum_{(L_{1},L_{2})\in {\cal L}}\sum_{\gamma\in \Gamma} C(\gamma,L_{2})\int_{A(\gamma)} E[a,L_{1},L_{2}](\bar{Y}_{1},\bar{Y}_{2})\,da,$$
o\`u on a pos\'e
$$E[a,L_{1},L_{2}](\bar{Y}_{1},\bar{Y}_{2})=E_{iso}[a,L_{1},L_{2}](\bar{Y}_{1},\bar{Y}_{2})+E_{an}[a,L_{1},L_{2}](\bar{Y}_{1},\bar{Y}_{2}).$$

Fixons $ \sharp=iso$ ou $an$, $(L_{1},L_{2})\in {\cal L}$, $\gamma\in \Gamma$ et  $a\in A(\gamma)$. 
  Rappelons que, pour $e\in {\cal E}$ et $u\in {\cal U}$, on a l'\'egalit\'e $\hat{i}_{\sharp}[a,e,u]=\frac{1}{2}(\hat{i}^+_{\sharp}[a,e,u]+\hat{i}^-_{\sharp}[a,e,u])$. Supposons   $n_{1}\not=0$ et $n_{2}\not=0$. On a d\'efini en 2.6 une application qui, \`a $\zeta_{1},\zeta_{2}=\pm$ associe $\zeta=\pm$ de sorte que les classes de conjugaison stable de $(X^{\zeta_{1}}(a^{L_1},w''),X^{\zeta_{2}}(a^{L_2},w''))$ et $X^{\zeta}(a,e,u)$ se correspondent. Notons-la   $Z$. Elle est surjective et ses fibres ont deux \'el\'ements. On a donc
$$\hat{i}_{\sharp}[a,e,u]=\frac{1}{4}\sum_{\zeta_{1},\zeta_{2}=\pm}\hat{i}^{Z(\zeta_{1},\zeta_{2})}_{\sharp}[a,e,u].$$
Utilisons le lemme 2.6. Alors
$$ \kappa^{L_{2}}(e)\kappa^{{\cal U}}(u)\hat{i}_{\sharp}[a,e,u]=\frac{1}{4}\sum_{\zeta_{1},\zeta_{2}=\pm}d(r',r'',\gamma,L_{2})$$
$$\Delta_{n_{1},\eta[L_{1},\gamma],n_{2},\eta[L_{2},\gamma]}((X^{\zeta_{1}}(a^{L_1},w'),X^{\zeta_{2}}(a^{L_2},w'')),X^{Z(\zeta_{1},\zeta_{2})}(a,e,u))\hat{i}^{Z(\zeta_{1},\zeta_{2})}_{\sharp}[a,e,u].$$
 L'\'egalit\'e (1) se r\'ecrit
$$(3) \qquad E_{\sharp}[a,L_{1},L_{2}](\bar{Y}_{1},\bar{Y}_{2})=\frac{1}{4}\sum_{\zeta_{1},\zeta_{2}=\pm} \sum_{Y}\Delta_{\bar{n}_{1},\eta_{1},\bar{n}_{2},\eta_{2}}((\bar{Y}_{1},\bar{Y}_{2}),Y) E_{\sharp}^{\zeta_{1},\zeta_{2}}[a,L_{1},L_{2}](Y),$$
o\`u
$$E_{\sharp}^{\zeta_{1},\zeta_{2}}[a,L_{1},L_{2}](Y)=$$
$$
\sum_{e\in {\cal E}}\sum_{u\in {\cal U}}\Delta_{n_{1},\eta[L_{1},\gamma],n_{2},\eta[L_{2},\gamma]}((X^{\zeta_{1}}(a^{L_1},w'),X^{\zeta_{2}}(a^{L_2},w'')),X^{Z(\zeta_{1},\zeta_{2})}(a,e,u))\hat{i}^{Z(\zeta_{1},\zeta_{2})}_{\sharp}[a,e,u](Y).$$
Fixons $\zeta_{1}$ et $\zeta_{2}$. Dans la formule  d\'efinissant $E_{\sharp}^{\zeta_{1},\zeta_{2}}[a,L_{1},L_{2}](Y)$, la somme est en fait limit\'ee aux $(e,u)$ tels que $\sharp(a,e,u)=\sharp$ (pour les autres, $\hat{i}^{Z(\zeta_{1},\zeta_{2})}_{\sharp}[a,e,u]$ est nulle). Les $X^{Z(\zeta_{1},\zeta_{2})}(a,e,u)$ parcourent un ensemble de repr\'esentants des classes de conjugaison par $G_{\sharp}(F)$ dans la classe de conjugaison stable correspondant \`a celle de $(X^{\zeta_{1}}(a^{L_1},w'),X^{\zeta_{2}}(a^{L_2},w''))$. On voit qu'aux diff\'erences de notations pr\`es et \`a une constante pr\`es, $E_{\sharp}^{\zeta_{1},\zeta_{2}}[a,L_{1},L_{2}](Y)$ n'est autre que le membre de droite de l'\'egalit\'e de la conjecture 2 de \cite{W4} VIII.7. Cette conjecture est d\'emontr\'ee depuis que Ngo Bao Chau a d\'emontr\'e le lemme fondamental. On va  appliquer cette conjecture. Posons
$$C_{\sharp}(n_{1},\eta[L_{1},\gamma],n_{2},\eta[L_{2},\gamma])=\gamma_{\psi_{F}}(\mathfrak{g}_{n_{1},\eta[L_{1},\gamma],iso}(F))\gamma_{\psi_{F}}(\mathfrak{g}_{n_{2},\eta[L_{2},\gamma],iso}(F))\gamma_{\psi_{F}}(\mathfrak{g}_{\sharp}(F))^{-1},$$
avec les d\'efinitions de \cite{W4} VIII.5. Pour $j=1,2$ et pour un \'el\'ement r\'egulier $Y_{j}$ de $\mathfrak{g}_{n_{j},\eta[L_{j},\gamma],iso}(F)$, notons $z_{iso}(Y_{j})$ le nombre de classes de conjugaison par $G_{n_{j},\eta[L_{j},\gamma],iso}(F)$ dans la classe de conjugaison stable de $Y_{j}$.  Rappelons que $X^{\zeta_{1}}(a^{L_{1}},w')$ est un \'el\'ement d'une certaine  classe de conjugaison stable, cf. 2.6.  Un ensemble de  classes de conjugaison par $G_{n_{1},\eta[L_{1},\gamma],iso}(F)$ dans cette classe de conjugaison stable est form\'e d'\'el\'ements $X^{\zeta_{1}}(a^{L_{1}},e_{1},u_{1})$, o\`u $e_{1}$ et $u_{1}$ d\'ecrivent des ensembles ${\cal E}_{1}$ et ${\cal U}_{1}$ analogues \`a ${\cal E}$ et ${\cal U}$, avec la restriction $\sharp(a^{L_{1}},e_{1},u_{1})=iso$. On d\'efinit pour ces \'el\'ements une fonction $\hat{i}^{\zeta_{1}}_{iso}[a^{L_{1}},e_{1},u_{1}]$ similaire \`a $\hat{i}^{\zeta}_{\sharp}[a,e,u]$. Si la condition $\sharp(a^{L_{1}},e_{1},u_{1})=iso$ n'est pas v\'erifi\'ee, on pose $\hat{i}^{\zeta_{1}}_{iso}[a^{L_{1}},e_{1},u_{1}]=0$. On pose des d\'efinitions analogues en rempla\c{c}ant l'indice $1$ par $2$ (et $w'$ par $w''$). Alors la conjecture 2 de \cite{W4} VIII.7 fournit l'\'egalit\'e
 $$(4) \qquad E_{\sharp}^{\zeta_{1},\zeta_{2}}[a,L_{1},L_{2}](Y)=C_{\sharp}(n_{1},\eta[L_{1},\gamma],n_{2},\eta[L_{2},\gamma])\sum_{ e_{1}\in {\cal E}_{1},u_{1}\in {\cal U}_{1},e_{2}\in {\cal E}_{2},u_{2}\in {\cal U}_{2}}\sum_{Y_{1},Y_{2}}$$
 $$z_{iso}(Y_{1})^{-1}z_{iso}(Y_{2})^{-1}\Delta_{n_{1},\eta[L_{1},\gamma],n_{2},\eta[L_{2},\gamma]}((Y_{1},Y_{2}),Y)\hat{i}_{iso}^{\zeta_{1}}[a^{L_{1}},e_{1},u_{1}](Y_{1})\hat{i}_{iso}^{\zeta_{2}}[a^{L_{2}},e_{2},u_{2}](Y_{2}),$$
 o\`u l'on somme sur les $(Y_{1},Y_{2})$, elliptiques r\'eguliers dans $\mathfrak{g}_{n_{1},\eta[L_{1},\gamma],iso}(F) \times \mathfrak{g}_{n_{2},\eta[L_{2},\gamma],iso}(F)$, \`a conjugaison pr\`es par $G_{n_{1},\eta[L_{1},\gamma],iso}(F)\times G_{n_{2},\eta[L_{2},\gamma],iso}(F)$. 
 
  Calculons $C_{\sharp}(n_{1},\eta[L_{1},\gamma],n_{2},\eta[L_{2},\gamma])$. D'apr\`es \cite{MW} 3.15, on a les \'egalit\'es 

$\gamma_{\psi_{F}}(\mathfrak{g}_{\sharp}(F))=\gamma_{\psi_{F}}^{n}sgn(\eta\varpi^{-val_{F}(\eta)})$, si $val_{F}(\eta)$ est pair;

$\gamma_{\psi_{F}}(\mathfrak{g}_{\sharp}(F))=\gamma_{\psi_{F}}^{n-1}$ si $val_{F}(\eta)$ est impair.

Le terme $\gamma_{\psi_{F}}$ est une constante de Weil \'el\'ementaire. Il v\'erifie l'\'egalit\'e $\gamma_{\psi_{F}}^2=sgn(-1)$. On voit tout d'abord que les termes ci-dessus ne d\'ependent pas de l'indice $\sharp$. La constante $C_{\sharp}(n_{1},\eta[L_{1},\gamma],n_{2},\eta[L_{2},\gamma])$ non plus, on la note simplement $C(n_{1},\eta[L_{1},\gamma],n_{2},\eta[L_{2},\gamma])$.
Pour $j=1,2$, on a des formules analogues pour $\gamma_{\psi_{F}}(\mathfrak{g}_{n_{j},\eta[L_{j},\gamma],iso}(F))$. On obtient les \'egalit\'es

$$ C(n_{1},\eta[L_{1},\gamma],n_{2},\eta[L_{2},\gamma])=\left\lbrace\begin{array}{cc}1,& pour\,\ val_{F}(\eta[L_{1},\gamma]) \,\,pair\\ &et\,\,val_{F}(\eta[L_{2},\gamma])\,\,pair;\\ sgn(\eta[L_{1},\gamma]\varpi^{-val_{F}(\eta[L_{1},\gamma])}),&pour\,\ val_{F}(\eta[L_{1},\gamma]) \,\,pair\\ &et\,\,val_{F}(\eta[L_{2},\gamma])\,\,impair;\\ sgn(\eta[L_{2},\gamma]\varpi^{-val_{F}(\eta[L_{2},\gamma])}),&pour\,\ val_{F}(\eta[L_{1},\gamma]) \,\,impair\\ &et\,\,val_{F}(\eta[L_{2},\gamma])\,\,pair;\\ sgn(-\eta\varpi^{-val_{F}(\eta)}),&pour\,\ val_{F}(\eta[L_{1},\gamma]) \,\,impair\\ &et\,\,val_{F}(\eta[L_{2},\gamma])\,\,impair.\\ \end{array}\right.$$

 Ins\'erons l'\'egalit\'e   (4) dans la formule (3). Pour $j=1,2$, posons
  $$\hat{i}_{iso}[a^{L_{j}},e_{j},u_{j}]=\frac{1}{2}(\hat{i}_{iso}^{+}[a^{L_{j}},e_{j},u_{j}]+\hat{i}_{iso}^{-}[a^{L_{j}},e_{j},u_{j}]).$$
On obtient
 $$  E_{\sharp}[a,L_{1},L_{2}](\bar{Y}_{1},\bar{Y}_{2})= C(n_{1},\eta[L_{1},\gamma],n_{2},\eta[L_{2},\gamma]) \sum_{ e_{1}\in {\cal E}_{1},u_{1}\in {\cal U}_{1},e_{2}\in {\cal E}_{2},u_{2}\in {\cal U}_{2}}\sum_{Y_{1},Y_{2}}$$
 $$z_{iso}(Y_{1})^{-1}z_{iso}(Y_{2})^{-1}S_{\sharp,L_{1},L_{2},\gamma}(\bar{Y}_{1},\bar{Y}_{2},Y_{1},Y_{2})\hat{i}_{iso}[a^{L_{1}},e_{1},u_{1}](Y_{1})\hat{i}_{iso}[a^{L_{2}},e_{2},u_{2}](Y_{2}),$$
 o\`u
 $$S_{\sharp,L_{1},L_{2},\gamma}(\bar{Y}_{1},\bar{Y}_{2},Y_{1},Y_{2})=\sum_{Y}\Delta_{\bar{n}_{1},\eta_{1},\bar{n}_{2},\eta_{2}}((\bar{Y}_{1},\bar{Y}_{2}),Y)\Delta_{n_{1},\eta[L_{1},\gamma],n_{2},\eta[L_{2},\gamma]}((Y_{1},Y_{2}),Y).$$
 Rappelons que $Y$ d\'ecrit ici les \'el\'ements elliptiques r\'eguliers de $\mathfrak{g}_{\sharp}(F)$, \`a conjugaison pr\`es par $G_{\sharp}(F)$. Puisque $E[a,L_{1},L_{2}](\bar{Y}_{1},\bar{Y}_{2})$ est la somme de $E_{iso}[a,L_{1},L_{2}](\bar{Y}_{1},\bar{Y}_{2})$ et de $E_{an}[a,L_{1},L_{2}](\bar{Y}_{1},\bar{Y}_{2})$, on a
  $$(5) \qquad   E[a,L_{1},L_{2}](\bar{Y}_{1},\bar{Y}_{2})= C(n_{1},\eta[L_{1},\gamma],n_{2},\eta[L_{2},\gamma]) \sum_{ e_{1}\in {\cal E}_{1},u_{1}\in {\cal U}_{1},e_{2}\in {\cal E}_{2},u_{2}\in {\cal U}_{2}}\sum_{Y_{1},Y_{2}}$$
 $$z_{iso}(Y_{1})^{-1}z_{iso}(Y_{2})^{-1}S_{L_{1},L_{2},\gamma}(\bar{Y}_{1},\bar{Y}_{2},Y_{1},Y_{2})\hat{i}_{iso}[a^{L_{1}},e_{1},u_{1}](Y_{1})\hat{i}_{iso}[a^{L_{2}},e_{2},u_{2}](Y_{2}),$$
 o\`u
 $$S_{L_{1},L_{2},\gamma}(\bar{Y}_{1},\bar{Y}_{2},Y_{1},Y_{2})=S_{iso,L_{1},L_{2},\gamma}(\bar{Y}_{1},\bar{Y}_{2},Y_{1},Y_{2})+S_{an,L_{1},L_{2},\gamma}(\bar{Y}_{1},\bar{Y}_{2},Y_{1},Y_{2}).$$
 On a suppos\'e $n_{1}\not=0$ et $n_{2}\not=0$. Supposons que ces conditions ne sont pas v\'erifi\'ees, par exemple $n_{2}=0$. On reprend le calcul. Les seules diff\'erences sont que les termes index\'es par $2$ doivent dispara\^{\i}tre. En particulier les $\zeta_{2}$.  Les $\frac{1}{4}$ figurant dans les calculs  sont remplac\'es par des $\frac{1}{2}$. On obtient encore la formule (5), dont on fait dispara\^{\i}tre les termes index\'es par $2$. 
 
 Les propri\'et\'es d'inversion des facteurs de transfert nous disent que $S_{L_{1},L_{2},\gamma}(\bar{Y}_{1},\bar{Y}_{2},Y_{1},Y_{2})=0$ sauf si $\bar{n}_{1}=n_{1}$, $\eta_{1}=\eta[L_{1},\gamma]$, $\bar{n}_{2}=n_{2}$ et $\eta_{2}=\eta[L_{2},\gamma]$.  

Supposons que $(\bar{n}_{1},\bar{n}_{2})\not=(n_{1},n_{2})$ ou que $(\bar{n}_{1},\bar{n}_{2})=(n_{1},n_{2})$ mais que $(\eta_{1},\eta_{2})$ ne v\'erifie pas la condition (1) de 2.2  (c'est-\`a-dire $val_{F}(\eta_{1})\equiv t'_{1}=t''_{1}\,\,mod \,\,2{\mathbb Z}$ et $val_{F}(\eta_{2})\equiv t'_{2}=t''_{2}\,\,mod \,\,2{\mathbb Z}$). Remarquons que les couples $(\eta[L_{1},\gamma],\eta[L_{2},\gamma])$ d\'ependent de $L_{1}$, $L_{2}$ et $\gamma$ mais v\'erifient par construction cette condition. Alors $S_{L_{1},L_{2},\gamma}(\bar{Y}_{1},\bar{Y}_{2},Y_{1},Y_{2})$ est toujours nulle. Donc aussi $E[a,L_{1},L_{2}](\bar{Y}_{1},\bar{Y}_{2})$. La formule (2) implique $J^{endo}(\bar{Y}_{1},\bar{Y}_{2},f)=0$. Ceci \'etant vrai pour tous $\bar{Y}_{1},\bar{Y}_{2}$, on a $transfert_{\bar{n}_{1},\eta_{1},\bar{n}_{2},\eta_{2}}(f)=0$. Cela d\'emontre le (ii) du lemme 2.2.

\bigskip

\subsection{D\'emonstration du (i) du lemme 2.2}
On poursuit le calcul en supposant que  $\bar{n}_{1}=n_{1}$, $\bar{n}_{2}=n_{2}$ et que $(\eta_{1},\eta_{2})$ v\'erifient la condition (1) de 2.2. 

Le terme $E[a,L_{1},L_{2}](\bar{Y}_{1},\bar{Y}_{2})$ n'est non nul que si $\eta[L_{1},\gamma]=\eta_{1}$ et $\eta[L_{2},\gamma]=\eta_{2}$. Le couple $(L_{1},L_{2})$ \'etant fix\'e, on note $\Gamma[L_{1},L_{2}]$ l'ensemble des $\gamma\in \Gamma$ pour lesquels cette condition est v\'erifi\'ee. Supposons qu'elle le soit. Alors $S_{L_{1},L_{2},\gamma}(\bar{Y}_{1},\bar{Y}_{2},Y_{1},Y_{2})$ est non nulle si et seulement si $( Y_{1},Y_{2})$ est stablement conjugu\'e \`a un couple d\'eduit de $(\bar{Y}_{1},\bar{Y}_{2})$ par un automorphisme de la donn\'ee endoscopique d\'efinie par $(n_{1},n_{2})$. Si $(Y_{1},Y_{2})$ est un tel \'el\'ement, $S_{L_{1},L_{2},\gamma}(\bar{Y}_{1},\bar{Y}_{2},Y_{1},Y_{2})$  est \'egal au nombre de termes de la sommation, c'est-\`a-dire au nombre  $z(\bar{Y}_{1},\bar{Y}_{2})$   de classes de conjugaison dans la classe  totale de conjugaison  stable de $\mathfrak{g}_{iso}(F)\cup \mathfrak{g}_{an}(F)$ correspondant \`a celle de $(\bar{Y}_{1},\bar{Y}_{2})$. 

 Il faut ici faire attention. Parce qu'on travaille simultan\'ement avec les deux groupes $G_{iso}$ et $G_{an}$, on a raffin\'e la notion de donn\'ee endoscopique. Il faut raffiner aussi celle d'automorphismes de cette donn\'ee. On voit que, si $n_{1}$ et $n_{2}$ sont tous deux non nuls, il y a deux automorphismes de notre donn\'ee: outre l'identit\'e, il y a l'action d'un \'el\'ement de $O^-(Q_{n_{1},\eta_{1},iso})\times O^-(Q_{n_{2},\eta_{2},iso})$.  
 
 {\bf Remarque.} Dans le cas o\`u $n_{1}=n_{2}$ et $\eta_{1}=\eta_{2}$, la permutation des deux groupes $G_{n_{1},\eta_{1},iso}$ et $G_{n_{2},\eta_{2},iso}$ n'est pas un automorphisme pour notre notion de donn\'ee endoscopique. 
 
 \bigskip
 
 Soit $(Y_{1},Y_{2})$ un couple stablement conjugu\'e \`a $(\bar{Y}_{1},\bar{Y}_{2})$. Notons $(\underline{Y}_{1},\underline{Y}_{2})$ le couple d\'eduit de $( Y_{1},Y_{2})$ par l'automorphisme non trivial. Il n'est pas stablement conjugu\'e \`a $(\bar{Y}_{1},\bar{Y}_{2})$. 
 Mais les fonctions $\hat{i}_{iso}[a^{L_{1}},e_{1},u_{1}]$ et $\hat{i}_{iso}[a^{L_{2}},e_{2},u_{2}]$ sont invariantes par conjugaison non seulement par les groupes sp\'eciaux orthogonaux, mais par les groupes orthogonaux tout entiers. Cela vient de la d\'efinition de ces fonctions: par exemple, $\hat{i}_{iso}[a^{L_{1}},e_{1},u_{1}]$ est la somme de deux fonctions associ\'ees \`a $X^+(a^{L_{1}},e_{1},u_{1})$  et $X^-(a^{L_{1}},e_{1},u_{1})$  et ces deux \'el\'ements sont justement conjugu\'es par un \'el\'ement du groupe orthogonal de d\'eterminant $-1$. On en d\'eduit l'\'egalit\'e
 $$\hat{i}_{iso}[a^{L_{1}},e_{1},u_{1}](\underline{Y}_{1})\hat{i}_{iso}[a^{L_{2}},e_{2},u_{2}](\underline{Y}_{2})=\hat{i}_{iso}[a^{L_{1}},e_{1},u_{1}]( Y_{1})\hat{i}_{iso}[a^{L_{2}},e_{2},u_{2}](Y_{2}).$$
   Dans la formule (5) du paragraphe pr\'ec\'edent, les couples $(Y_{1},Y_{2})$ sont simplement  ceux qui sont stablement conjugu\'es \`a $(\bar{Y}_{1},\bar{Y}_{2})$  ainsi que leurs images $(\underline{Y}_{1},\underline{Y}_{2})$ par automorphisme.  Ce qui pr\'ec\`ede montre que  ces derniers donnent la m\^eme contribution que les premiers. On peut donc sommer seulement sur les premiers, en multipliant le tout par $2$.  On calcule facilement $z(\bar{Y}_{1},\bar{Y}_{2})=4z_{iso}(\bar{Y}_{1})z_{iso}(\bar{Y}_{2})$. En posant $\beta=1$ (sous notre hypoth\`ese $n_{1}\not=0$, $n_{2}\not=0$), la formule (5) devient
$$E[a,L_{1},L_{2}](\bar{Y}_{1},\bar{Y}_{2})=2^{1+2\beta} C(n_{1},\eta _{1},n_{2},\eta _{2}) \sum_{(Y_{1},Y_{2})}$$
$$\sum_{ e_{1}\in {\cal E}_{1},u_{1}\in {\cal U}_{1},e_{2}\in {\cal E}_{2},u_{2}\in {\cal U}_{2}}\hat{i}_{iso}[a^{L_{1}},e_{1},u_{1}](Y_{1})\hat{i}_{iso}[a^{L_{2}},e_{2},u_{2}](Y_{2}),$$
o\`u on  somme sur les couples $(Y_{1},Y_{2})$ stablement conjugu\'es \`a $(\bar{Y}_{1}
,\bar{Y}_{2})$, pris \`a conjugaison pr\`es.

Si  par exemple $n_{2}=0$, notre donn\'ee endoscopique n'a pas d'autre automorphisme que l'identit\'e. Mais on a cette fois l'\'egalit\'e $z(\bar{Y}_{1})=2z_{iso}(\bar{Y}_{1})$. Le calcul conduit \`a la m\^eme \'egalit\'e que ci-dessus, o\`u maintenant $\beta=0$.

La formule (2) de 2.7 devient
$$(1) \qquad J^{endo}(\bar{Y}_{1},\bar{Y}_{2},f)=2^{1+2\beta}C(n_{1},\eta _{1},n_{2},\eta _{2})\vert {\cal L}\vert ^{-1} \sum_{(Y_{1},Y_{2})}\sum_{(L_{1},L_{2})\in {\cal L}}\sum_{\gamma\in \Gamma[L_{1},L_{2}]}C(\gamma,L_{2})\int_{A(\gamma)}$$
$$ \sum_{ e_{1}\in {\cal E}_{1},u_{1}\in {\cal U}_{1},e_{2}\in {\cal E}_{2},u_{2}\in {\cal U}_{2}}\hat{i}_{iso}[a^{L_{1}},e_{1},u_{1}](Y_{1})\hat{i}_{iso}[a^{L_{2}},e_{2},u_{2}](Y_{2})\,da,$$
la somme en $(Y_{1},Y_{2})$ \'etant comme ci-dessus.

   Rappelons qu'en 2.4, on a fix\'e des ensembles $\Gamma_{0}$ et $\Gamma_{j}$ pour $j=(R-r)/2+1,...,R$. Introduisons l'ensemble $\boldsymbol{\Gamma}$ des familles $\Gamma^*=((\Gamma_{1,j})_{j=1,...,t'_{1}},(\Gamma_{2,j})_{j=1,...,t'_{2}})$ qui v\'erifient les conditions suivantes:

pour $j=t'_{2}+1,...,t'_{1}$, $\Gamma_{1,j}=\Gamma_{j+(R-r)/2}$;

pour $j=1,...,t'_{2}$, $\Gamma_{1,j}$, resp. $\Gamma_{2,j}$, est un  sous-ensemble de $\Gamma_{0}$ \`a deux \'el\'ements, qui est un ensemble de repr\'esentants de $\mathfrak{o}^{\times}/\mathfrak{o}^{\times2}$; on  impose $\Gamma_{1,j}\cap \Gamma_{2,j}=\emptyset$.

 Pour une telle  famille, on note $\Gamma^*_{\eta_{1}}$ l'ensemble des familles $\gamma_{1}=(\gamma_{1,j})_{j=1,...,t'_{1}}\in \prod_{j=1,...,t'_{1}}\Gamma_{1,j} $ qui v\'erifient la condition

$sgn(\eta_{1}\varpi^{-val_{F}(\eta_{1})}\prod_{j=1,...,t'_{1}}\gamma_{1,j})=sgn_{CD}(w')$.

On d\'efinit l'ensemble $\Gamma^*_{\eta_{2}}$ en changeant les indices $1$ en $2 $ et $w'$ en $w''$. 

Fixons $(L_{1},L_{2})\in {\cal L}$. Pour $(\gamma_{1},\gamma_{2})\in \Gamma^*_{\eta_{1}}\times \Gamma^*_{\eta_{2}}$, on d\'efinit une famille $\gamma=(\gamma_{l})_{l=1,...,R}$ de la fa\c{c}on suivante

pour $j= t'_{2}+1,...,t'_{1}$, $\gamma_{j+(R-r)/2}=\gamma_{1,j}$;

pour $j=1,...,t'_{2}$, $\gamma_{l_{1}(j)}=\gamma_{1,j}$ et $\gamma_{l_{2}(j)}=\gamma_{2,j}$.

On v\'erifie que cette famille appartient \`a $\Gamma$. Cela d\'efinit une application
$$\pi_{L_{1};L_{2}}:\sqcup_{\Gamma^*\in \boldsymbol{\Gamma}}\Gamma^*_{\eta_{1}}\times \Gamma^*_{\eta_{2}}\to \Gamma.$$
Montrons que

(2) l'image de $\pi_{L_{1},L_{2}}$ est \'egale \`a $\Gamma[L_{1},L_{2}]$; pour $\gamma=(\gamma_{l})_{l=1,...,R}\in \Gamma[L_{1},L_{2}]$, la fibre de $\pi_{L_{1},L_{2}}$ au-dessus de $\gamma$ a pour nombre d'\'el\'ements
$$\sigma^*(\gamma)=(\frac{q-3}{4})^{t'_{2}}\prod_{l\in \hat{J}}(q-2+sgn(\gamma_{l-1}\gamma_{l})).$$
 
Soit $\Gamma^*\in \boldsymbol{\Gamma}$ et $(\gamma_{1},\gamma_{2})\in \Gamma^*_{\eta_{1}}\times \Gamma^*_{\eta_{2}}$. Posons $\gamma=\pi_{L_{1},L_{2}}(\gamma_{1},\gamma_{2})$. On a d\'efini en 2.6 des \'el\'ements $\gamma^{L_{1}}$ et $\gamma^{L_{2}}$. Ce ne sont autres que $\gamma_{1}$ et $\gamma_{2}$. Le terme $\eta[L_{2},\gamma]$ d\'efini en 2.6 est caract\'eris\'e par les relations $val_{F}(\eta[L_{2},\gamma])\equiv t'_{2}\,\,mod \,\,2{\mathbb Z}$ et $sgn(\eta[L_{2},\gamma]\varpi^{-val_{F}(\eta[L_{2},\gamma])}\prod_{j=1,...t'_{2}}\gamma_{2,j})=sgn_{CD}(w'')$. Par hypoth\`ese, $\eta_{2}$ v\'erifie la premi\`ere relation. Par d\'efinition de $\Gamma^*_{\eta_{2}}$, $\eta_{2}$ v\'erifie aussi la seconde. Donc $\eta[L_{2},\gamma]=\eta_{2}$, puis $\eta[L_{1},\gamma]=\eta_{1}$ puisque $\eta[L_{1},\gamma]\eta[L_{2},\gamma]=\eta=\eta_{1}\eta_{2}$. Par d\'efinition de $\Gamma[L_{1},L_{2}]$, $\gamma$ appartient donc \`a cet ensemble, ce qui prouve la premi\`ere assertion. Soit maintenant $\gamma\in \Gamma[L_{1},L_{2}]$. La premi\`ere partie de la preuve  montre que le nombre d'\'el\'ements de la fibre de $\pi_{L_{1},L_{2}}$ au-dessus de $\gamma$ est \'egal au nombre d'\'el\'ements $\Gamma^*\in \boldsymbol{\Gamma}$ tels que $(\gamma^{L_{1}},\gamma^{L_{2}})$ appartiennent \`a $(\prod_{j=1,...,t'_{1}}\Gamma_{1,j} )\times (\prod_{j=1,...,t'_{2}}\Gamma_{2,j})$.  Ce nombre est le produit sur $j=1,...,t'_{2}$ du nombre des paires $(\Gamma_{1,j},\Gamma_{2,j})$ possibles dont le produit contient $(\gamma_{1,j},\gamma_{2,j})$. Pour simplifier la notation, posons $x=\gamma_{1,j}$, $y=\gamma_{2,j}$. Les \'el\'ements $ x$ et $y$ appartiennent \`a $\Gamma_{0}$ et sont distincts puisqu'ils proviennent de $\gamma\in \Gamma$. Les paires  $(\Gamma_{1,j},\Gamma_{2,j})$  possibles sont de la forme $\Gamma_{1,j}=\{ x,x'\}$, $\Gamma_{2,j}=\{ y,y'\}$, o\`u $ x'$ et $y'$ sont deux \'el\'ements de $\Gamma_{0}$ qui v\'erifient les conditions suivantes:

$x'\not\in x\mathfrak{o}^{\times2}$, $y'\not\in y\mathfrak{o}^{\times2}$, $x'\not=y$, $y'\not=x$ et $x'\not=y'$. 

Notons $C'_{x}$ l'ensemble des $x'\in \Gamma_{0}$ tels que $x'\not\in x\mathfrak{o}^{\times2}$ et d\'efinissons de m\^eme $C'_{y}$. 
Supposons $y\in x\mathfrak{o}^{\times2}$, autrement dit $sgn(xy)=1$. On a $C'_{x}=C'_{y}$. Les deux premi\`eres relations  signifient que $x',y'\in C'_{x}$. Cela  entra\^{\i}ne  les deux relations suivantes $x'\not=y$ et $y'\not=x$.   La derni\`ere relation signifie que $(x',y')$ n'appartient pas \`a la diagonale de $C'_{x}\times C'_{x}$. Le nombre de $(x',y')$ possibles est donc $\vert C'_{x}\vert ^2-\vert C'_{x}\vert $. On a $\vert C'_{x}\vert =(q-1)/2$. Le nombre de couples possibles est donc $(q-1)(q-3)/4$, ou encore $(q-2+sgn(xy))(q-3)/4$.
Supposons maintenant $y\not\in x\mathfrak{o}^{\times2}$, autrement dit $sgn(xy)=-1$. On a $C'_{x}\not=C'_{y}$. Les premi\`eres et troisi\`emes conditions signifient que $x'\in C'_{x}-\{y\}$. De m\^eme, les deuxi\`emes et quatri\`emes conditions signifient que $y'\in C'_{y}-\{x\}$. La cinqui\`eme condition est automatique. le nombre de $(x',y')$ possibles est donc $(\vert C'_{x}\vert -1)(\vert C'_{y}\vert -1)=((q-3)/2)^2$, ou encore $(q-2+sgn(xy))(q-3)/4$. L'assertion (2) en r\'esulte. 

Un calcul similaire d\'emontre la relation suivante

(3) le nombre d'\'el\'ements de $\boldsymbol{\Gamma}$ est \'egal \`a $2^{-4t'_{2}}(q-1)^{2t'_{2}}(q-3)^{2t'_{2}}$. 

Pour $\Gamma^*\in \boldsymbol{\Gamma}$, $(L_{1},L_{2})\in {\cal L}$ et $(Y_{1},Y_{2})$ stablement conjugu\'e \`a $(\bar{Y}_{1},\bar{Y}_{2})$, posons
$$(4) \qquad {\cal J}_{\Gamma^*,L_{1},L_{2}}(Y_{1},Y_{2})=2^{1+2\beta}C(n_{1},\eta_{1},n_{2},\eta_{2})\sum_{\gamma_{1}\in \Gamma^*_{\eta_{1}}}\sum_{\gamma_{2}\in \Gamma^*_{\eta_{2}}}\sigma^*( \gamma)^{-1}C(\gamma,L_{2})$$
$$\int_{A(\gamma)} \sum_{ e_{1}\in {\cal E}_{1},u_{1}\in {\cal U}_{1},e_{2}\in {\cal E}_{2},u_{2}\in {\cal U}_{2}}\hat{i}_{iso}[a^{L_{1}},e_{1},u_{1}](Y_{1})\hat{i}_{iso}[a^{L_{2}},e_{2},u_{2}](Y_{2})\,da$$
o\`u, pour simplifier, on a pos\'e $\gamma=\pi_{L_{1},L_{2}}(\gamma_{1},\gamma_{2})$. 
En utilisant (2), on obtient l'\'egalit\'e
$$(5) \qquad J^{endo}(\bar{Y}_{1},\bar{Y}_{2},f)= \vert {\cal L}\vert ^{-1}\sum_{Y_{1},Y_{2}} \sum_{(L_{1},L_{2})\in {\cal L}}\sum_{\Gamma^*\in \boldsymbol{\Gamma}}{\cal J}_{\Gamma^*,L_{1},L_{2}}(Y_{1},Y_{2}).$$

Fixons $\Gamma^*\in \boldsymbol{\Gamma}$, $(L_{1},L_{2})\in {\cal L}$ et $(Y_{1},Y_{2})$ stablement conjugu\'e \`a $(\bar{Y}_{1},\bar{Y}_{2})$. Les fonctions $f_{1,\eta_{1}}$ et $f_{2,\eta_{2}}$ qui interviennent dans le lemme 2.2 sont des cas particuliers de notre fonction $f$.  Consid\'erons par exemple le cas de l'indice $1$.  On passe de $f$ \`a $f_{1,\eta_{1}}$ en rempla\c{c}ant $n$, $\eta$, $r'$, $r''$, $N'$, $N''$, $w'$  par $n_{1}$, $\eta_{1}$,  $t'_{1}$, $t''_{1}$, $N'$, $0$, $w'$. Le terme correspondant \`a $w''$ dispara\^{\i}t. L'int\'egrale orbitale $J(Y_{1},f_{1,\eta_{1},iso})$  est calcul\'ee par une formule analogue \`a celle du lemme 2.4. On peut choisir l'ensemble $\Gamma^*_{\eta_{1}}$ comme analogue de $\Gamma$. Puisque $t'_{1}=t''_{1}$, l'analogue de $\hat{J}$ est vide. Les analogues de ${\cal E}^0$ et ${\cal U}$ ne sont autres que les ensembles ${\cal E}_{1}$ et ${\cal U}_{1}$ d\'ej\`a introduits. L'analogue de $\kappa^0$ est trivial. Parce que $w''$ dispara\^{\i}t, l'analogue de $\kappa^{{\cal U}}$ est trivial. L'analogue de la constante $C(r',r'')$ vaut $2^{1-2t'_{1}}$. On note $\alpha(t'_{1},w',\eta_{1})$ l'analogue de la constante $\alpha(r',r'',w',w'')$. On a donc
$$(6) \qquad J(Y_{1},f_{1,\eta_{1},iso})=C(w')2^{1-2t'_{1}}\alpha(t'_{1},w',\eta_{1})\sum_{\gamma_{1}\in \Gamma^*_{\eta_{1}}}\sigma(\gamma_{1})$$
$$\int_{A(\gamma_{1})}\sum_{e\in {\cal E}_{1},u_{1}\in {\cal U}_{1}}\hat{i}_{iso}[a_{1},e_{1},u_{1}](Y_{1})\,da_{1}.$$
  Il faut toutefois faire attention: cette formule n'est  valable que si $n_{1}>0$.  En effet, si $n_{1}=0$,  le terme de gauche doit  dispara\^{\i}tre des calculs qui suivent, autrement dit on doit consid\'erer qu'il vaut $1$. Les sommes de droite disparaissent autrement dit valent $1$. Les deux constantes $C(w')$ et  $\alpha(t'_{1},w',\eta_{1})$ valent aussi $1$. Mais $2^{1-t'_{1}}$ vaut $2$ et l'\'egalit\'e n'est pas v\'erifi\'ee. Bien s\^ur, c'est le produit des constantes pour les indices $1$ et $2$ qui va intervenir. La puissance de $2$ qui intervient dans ce produit est donc $2^{2-2t'_{1}-2t'_{2}}$ si $n_{1}\not=0$ et $n_{2}\not=0$ et c'est $2^{1-2t'_{1}-2t'_{2}}$ sinon. Autrement dit, c'est $2^{1+\beta+2t'_{1}+2t'_{2}}$. 
  
  Soient $\gamma_{1}\in \Gamma^*_{\eta_{1}}$, $\gamma_{2}\in \Gamma^*_{\eta_{2}}$, posons  $\gamma=\pi_{L_{1},L_{2}}(\gamma_{1},\gamma_{2})$. Montrons qu'on  a l'\'egalit\'e
  
  $$(7) \qquad d(r',r'',\gamma,L_{2})\sigma(\gamma)\sigma^*(\gamma)^{-1}\sigma(\gamma_{1})^{-1}\sigma(\gamma_{2})^{-1}=C_{3} ,$$
  o\`u
  $$C_{3}=sgn(-1)^{r(R-r)/2}sgn(\eta\varpi^{-val_{F}(\eta)})^{(R-r)/2}sgn_{CD}(w')^{(R-r)/2}sgn_{CD}(w'')^{(R-r)/2+val_{F}(\eta)}$$
  $$(\frac{q-3}{4})^{-t'_{2}}sgn(-1)^{val_{F}(\eta_{2})B}sgn(\eta_{2}\varpi^{-val_{F}(\eta_{2})})^Bsgn_{CD}(w'')^B .$$
  
  Rappelons les d\'efinitions:
  $$d(r',r'',\gamma,L_{2})=sgn(\eta\varpi^{-val_{F}(\eta)})^{(R-r)/2}sgn_{CD}(w')^{(R-r)/2}sgn_{CD}(w'')^{(R-r)/2+val_{F}(\eta)}$$
$$\left(\prod_{j\in \hat{J}}sgn(\gamma_{j-1}\gamma_{j})^{j/2-1}sgn(\gamma_{j-1}-\gamma_{j}))\right)\left(\prod_{j=R-r+1,...,R}sgn(\gamma_{j})^{(R-r)/2}\right)$$
$$sgn(-1)^{val_{F}(\eta_{2})B}sgn(\eta_{2}\varpi^{-val_{F}(\eta_{2})})^Bsgn_{CD}(w'')^B;$$

$$\sigma(\gamma)=\left(\prod_{j\in \hat{J}}(q-2+sgn(\gamma_{j-1}\gamma_{j}))sgn(\gamma_{j-1}\gamma_{j}(\gamma_{j-1}-\gamma_{j}))\right)$$
$$\prod_{j=R-r+1,...R;\, j\,\, impair}sgn(-\gamma_{j});$$

 $$\sigma^*(\gamma)=(\frac{q-3}{4})^{t'_{2}}\prod_{l\in \hat{J}}(q-2+sgn(\gamma_{l-1}\gamma_{l})).$$
 
 Les formules pour $\sigma(\gamma_{1})$ et $\sigma(\gamma_{2})$ se simplifient puisque $t'_{1}=t''_{1}$ et $t'_{2}=t''_{2}$:
 
 $$\sigma(\gamma_{1})=\prod_{j=1,...,t'_{1}; \,j\,\,impair}sgn(-\gamma_{1,j}),$$
 
  $$\sigma(\gamma_{2})=\prod_{j=1,...,t'_{2};\, j\,\,impair}sgn(-\gamma_{2,j}).$$
  
  Les restrictions $j$ impair figurant dans les diff\'erents produits peuvent \^etre lev\'es en \'el\'evant le terme correspondant (qui est un signe) \`a la puissance $j$. Consid\'erons l'intervention dans le membre de gauche de (6)  d'un terme $\gamma_{j}$ pour $j=R-r+1,...,R$. Il intervient dans $d(r',r'',\gamma,L_{2})$ par un facteur $sgn(\gamma_{j})^{(R-r)/2}$ et dans $\sigma(\gamma)$ par un facteur $sgn(-\gamma_{j})^{j}$. Puisque $\gamma_{1,j-(R-r)/2}=\gamma_{1,j}$, il intervient dans $\sigma(\gamma_{1})$ par un facteur $sgn(-\gamma_{j})^{j-(R-r)/2}$. Le produit de ces contributions est $sgn(-1)^{(R-r)/2}$. Leur produit sur tous les $j=R-r+1,...,R$ est $sgn(-1)^{r(R-r)/2}$. Consid\'erons maintenant l'intervention des termes $\gamma_{j-1}$ et $\gamma_{j}$ pour un $j\in \hat{J}$. Le terme $sgn(\gamma_{j-1}-\gamma_{j})$ intervient dans $d(r',r'',\gamma,L_{2})$ et dans $\sigma(\gamma)$. Il dispara\^{\i}t. Le terme $q-2+sgn(\gamma_{j-1}\gamma_{j}$ intervient dans $\sigma(\gamma)$ et son inverse intervient dans $\sigma^*(\gamma)^{-1}$. Il dispara\^{\i}t. Il intervient $sgn(\gamma_{j-1}\gamma_{j})^{j/2-1}$ dans $d(r',r'',\gamma,L_{2})$ et $sgn(\gamma_{j-1}\gamma_{j})$ dans $\sigma(\gamma)$. Puisque $\{\gamma_{j-1},\gamma_{j}\}=\{\gamma_{1,j/2},\gamma_{2,j/2}\}$, il intervient aussi dans $\sigma(\gamma_{1})\sigma(\gamma_{2})$ le terme $sgn(\gamma_{j-1}\gamma_{j})^{j/2}$. Le produit de ces termes vaut $1$. En r\'esum\'e, la contribution des termes d\'ependant des $\gamma_{j}$ est $sgn(-1)^{r(R-r)/2}$.  Outre ce terme, il reste le produit des constantes. Le tout donne la formule   (7).

Dans la formule (4) interviennent des termes $\gamma$, $a$, $a^{L_{1}}$ et $a^{L_{2}}$. On a vu dans la preuve de (2) que $\gamma^{L_{1}}=\gamma_{1}$ et $\gamma^{L_{2}}=\gamma_{2}$. On voit que, quand $a$ d\'ecrit $A(\gamma)$, le couple $(a^{L_{1}},a^{L_{2}})$ d\'ecrit $A(\gamma_{1})\times A(\gamma_{2})$. La formule (4) se r\'ecrit
$${\cal J}_{\Gamma^*,L_{1},L_{2}}(Y_{1},Y_{2})=2^{1+2\beta}C(n_{1},\eta_{1},n_{2},\eta_{2})\sum_{\gamma_{1}\in \Gamma^*_{\eta_{1}}}\sum_{\gamma_{2}\in \Gamma^*_{\eta_{2}}}\sigma^*( \gamma)^{-1}C(\gamma,L_{2})$$
$$\int_{A(\gamma_{1})}\sum_{e_{1}\in {\cal E}_{1},u_{1}\in {\cal U}_{1}}\hat{i}_{iso}[a_{1},e_{1},u_{1}](Y_{1})\,da_{1}\int_{A(\gamma_{2})}\sum_{e_{2}\in {\cal E}_{2},u_{2}\in {\cal U}_{2}}\hat{i}_{iso}[a_{2},e_{2},u_{2}](Y_{2})\,da_{2}.$$
Soient $\gamma_{1}$, $\gamma_{2}$ et $\gamma$ intervenant dans cette formule. Rappelons la d\'efinition 
$$C(\gamma,L_{2})=d(r',r'',\gamma,L_{2})C(w')C(w'')C(r',r'')\alpha(r',r'',w',w'')\sigma(\gamma).$$
Posons 
$$C_{4}=2^{\beta+2t'_{1}+2t'_{2}}C_{3}C(n_{1},\eta_{1},n_{2},\eta_{2})C(r',r'')\alpha(r',r'',w',w'')\alpha(t'_{1},w',\eta_{1})\alpha(t'_{2},w'',\eta_{2}).$$
En vertu de (7), on a
$$2^{1+2\beta}C(n_{1},\eta_{1},n_{2},\eta_{2})\sigma^*( \gamma)^{-1}C(\gamma,L_{2})=C_{4}2^{1+\beta-2t'_{1}-2t'_{2}}C(w')C(w'')$$
$$ \alpha(t'_{1},w',\eta_{1})\alpha(t'_{2},w'',\eta_{2})\sigma(\gamma_{1})\sigma(\gamma_{2}).$$
Alors la formule pr\'ec\'edente devient
$${\cal J}_{\Gamma^*,L_{1},L_{2}}(Y_{1},Y_{2})=C_{4}2^{1+\beta-2t'_{1}-2t'_{2}}C(w')C(w'') \alpha(t'_{1},w',\eta_{1})\alpha(t'_{2},w'',\eta_{2})$$
$$\sum_{\gamma_{1}\in \Gamma^*_{\eta_{1}}} \sigma(\gamma_{1})\int_{A(\gamma_{1})}\sum_{e_{1}\in {\cal E}_{1},u_{1}\in {\cal U}_{1}}\hat{i}_{iso}[a_{1},e_{1},u_{1}](Y_{1})\,da_{1}$$
$$\sum_{\gamma_{2}\in \Gamma^*_{\eta_{2}}}\sigma(\gamma_{2})\int_{A(\gamma_{2})}\sum_{e_{2}\in {\cal E}_{2},u_{2}\in {\cal U}_{2}}\hat{i}_{iso}[a_{2},e_{2},u_{2}](Y_{2})\,da_{2}.$$
Autrement dit, en  utilisant (6) et son analogue pour l'indice $2$:
$${\cal J}_{\Gamma^*,L_{1},L_{2}}(Y_{1},Y_{2})=C_{4}J(Y_{1},f_{1,\eta_{1},iso})J(Y_{2},f_{2,\eta_{2},iso}).$$
Cette expression ne d\'epend ni de $\Gamma^*$, ni de $L_{1}$ et $L_{2}$. Reportons cette valeur dans l'\'egalit\'e (5). La somme en $(L_{1},L_{2})$, pond\'er\'ee par $\vert {\cal L}\vert ^{-1}$, dispara\^{\i}t. La somme en $\Gamma^*$ se transforme en la multiplication par $\vert \boldsymbol{\Gamma}\vert $. Enfin, la somme en $(L_{1},L_{2})$ remplace les int\'egrales orbitales par leurs versions stables. D'o\`u
$${\cal J}_{\Gamma^*,L_{1},L_{2}}(Y_{1},Y_{2})=C_{4}\vert \boldsymbol{\Gamma}\vert S(\bar{Y}_{1},f_{1,\eta_{1},iso})S(\bar{Y}_{2},f_{2,\eta_{2},iso}).$$
Consid\'erons le cas de l'indice $1$.  Le quadruplet $(t'_{1},t''_{1},N',0)$  v\'erifie l'hypoth\`ese du (i) du th\'eor\`eme 3.20 de \cite{MW}. Donc $f_{1}$ est stable au sens de cette r\'ef\'erence. On sait de plus que le nombre de classes de conjugaison dans la classe de conjugaison stable de $\bar{Y}_{1}$ est \'egal \`a celui des classes de conjugaison dans la classe de conjugaison stable dans
$\mathfrak{g}_{n_{1},\eta_{1},an}(F)$ correspondant \`a celle de $\bar{Y}_{1}$. D'o\`u
$$S(\bar{Y}_{1},f_{1,\eta_{1},iso})=\frac{1}{2}S(\bar{Y}_{1},f_{1,\eta_{1}}).$$
Une fois de plus, cette \'egalit\'e n'est vraie que si $n_{1}\not=0$. Si $n_{1}=0$, le terme $\frac{1}{2}$ dispara\^{\i}t. Le produit de ces constantes sur les indices $1$ et $2$ vaut donc $2^{-1-\beta}$. D'o\`u l'\'egalit\'e
$${\cal J}_{\Gamma^*,L_{1},L_{2}}(Y_{1},Y_{2})=2^{-1-\beta}C_{4}\vert \boldsymbol{\Gamma}\vert S(\bar{Y}_{1},f_{1,\eta_{1}})S(\bar{Y}_{2},f_{2,\eta_{2}}).$$
Ceci \'etant vrai pour tout couple $(\bar{Y}_{1},\bar{Y}_{2})$, cela d\'emontre que 
$$transfert_{n_{1},\eta_{1},n_{2},\eta_{2}}(f)=2^{-1-\beta}C_{4}\vert \boldsymbol{\Gamma}\vert f_{1,\eta_{1}}\otimes f_{2,\eta_{2}}.$$
Pour obtenir l'assertion (i) du lemme 2.2, il reste \`a d\'emontrer l'\'egalit\'e
$$2^{-1-\beta}C_{4}\vert \boldsymbol{\Gamma}\vert =C_{\eta_{1},\eta_{2}}.$$
Il suffit pour cela de reprendre les d\'efinitions de toutes nos constantes. On laisse ce calcul p\'enible mais sans myst\`ere au lecteur. $\square$

  \bigskip

\section{D\'emonstration de la proposition 1.2}

\bigskip

\subsection{Descente du facteur de transfert}
On revient \`a nos deux groupes $G_{iso}$ et $G_{an}$ du paragraphe 1. On dit qu'un \'el\'ement $u\in G_{\sharp}(F)$ est topologiquement unipotent si ses valeurs propres, qui appartiennent \`a une certaine extension finie $F'$ de $F$, sont congrues \`a $1$ modulo l'id\'eal maximal de l'anneau des entiers de $F'$. L'application $E:X\mapsto E(X)=(1+X/2)(1-X/2)$ est une bijection entre l'ensemble des \'el\'ements  $X\in \mathfrak{g}_{\sharp}(F)$ topologiquement nilpotents et l'ensemble des \'el\'ements topologiquement unipotents de $G_{\sharp}(F)$. 

Soit $\sharp=iso$ ou $an$ et soit $x\in G_{\sharp}(F)$ un \'el\'ement elliptique et fortement r\'egulier.
On peut d\'ecomposer $x$ en un unique produit $x=su$, o\`u $s$ et $u$ commutent entre eux, $s$ est un \'el\'ement dont toutes les valeurs propres sont des racines de l'unit\'e d'ordre premier \`a $p$ et  $u$ est topologiquement unipotent. Au lieu de $u$, on consid\`ere plut\^ot   l'\'el\'ement topologiquement nilpotent de $X\in \mathfrak{g}_{\sharp}(F)$ tel que $u=E(X)$. On \'ecrit donc $x=sE(X)$. Parmi les valeurs propres de $s$, on distingue le nombre  $1$ qui intervient avec une multiplicit\'e impaire not\'ee $1+2n_{+}$ et le nombre $-1$ qui intervient avec une multiplicit\'e paire not\'ee $2n_{-}$ (qui peut \^etre nulle). On note $V_{+}$ et $V_{-}$ les espaces propres associ\'es, $Q_{+}$ et $Q_{-}$ les restrictions de $Q_{\sharp}$ \`a ces espaces et on pose $\eta_{+}=\eta(Q_{+})$ et $\eta_{-}=\eta(Q_{-})$. On note $G_{+}=SO(Q_{+})$, $G_{-}=SO(Q_{-})$. 
Le groupe de Galois $Gal(\bar{F}/F)$ agit sur l'ensemble des valeurs propres de $s$ diff\'erentes de $\pm 1$, en pr\'eservant leurs multiplicit\'es. Notons $I$ l'ensemble des orbites. Soit $i\in I$, fixons un \'el\'ement $s_{i}$ de cette orbite et notons $d_{i}$ sa multiplicit\'e. Posons $E_{i}=F[s_{i}]$. C'est une extension non ramifi\'ee de $F$.  Parce que $s$ appartient \`a un groupe orthogonal, il existe une sous-extension $E^{\natural}_{i}$ de $F$ telle que $[E_{i}:E^{\natural}_{i}]=2$ et que $norme_{E_{i}/E^{\natural}_{i}}(s_{i})=1$. Notons $V_{i}$ la somme des espaces propres associ\'es \`a des valeurs propres  conjugu\'ees de $s_{i}$ et notons $s_{\vert V_{i}}$ la restriction de $s$ \`a cet espace.  L'alg\`ebre $F[s_{\vert V_{i}}]$ s'identifie \`a $E_{i}$ par $s_{\vert V_{i}}\mapsto s_{i}$. Puisque $V_{i}$ est naturellement un  $F[s_{\vert V_{i}}]$-module, il devient un $E_{i}$-espace vectoriel, dont   la dimension est  $d_{i}$. On montre qu'il existe une unique forme hermitienne non d\'eg\'en\'er\'ee $Q_{i}$ sur $V_{i}$ (relative \`a l'extension $E_{i}/E_{i}^{\natural}$)  de sorte que, pour tout $v,v'\in V_{i}$, on ait l'\'egalit\'e
$$Q_{\sharp}(v,v')=[E_{i}:F]^{-1}trace_{E_{i}/F}(Q_{i}(v,v')).$$

Notons $G_{i}$ le groupe unitaire de $Q_{i}$. 
La composante neutre $G_{s}$ du commutant de $s$ dans $G_{\sharp}$ est le groupe
$$G_{+}\times G_{-}\times\prod_{i\in I}G_{i}.$$
On a $X\in \mathfrak{g}_{s}(F)$. C'est un \'el\'ement elliptique r\'egulier.

{\bf Remarques.} (1)  La condition d'ellipticit\'e emp\^eche le couple $(n_{-},\eta_{-})$ d'\^etre \'egal \`a $(1,1)$.

(2)  Fixons un \'el\'ement $\xi\in \mathfrak{o}^{\times}-\mathfrak{o}^{\times2}$. Pour $i\in I$, notons $Q_{\sharp,\vert V_{i}}$ la restriction de $Q_{\sharp}$ \`a $V_{i}$. On calcule facilement $\eta(Q_{\sharp,\vert V_{i}})=\xi^{d_{i}}$. Parce que l'on a suppos\'e $\eta(Q_{\sharp})=1$, on en d\'eduit l'\'egalit\'e
$$\eta_{+}\eta_{-}\xi^d=1,$$
o\`u $d=\sum_{i\in I}d_{i}$.
\bigskip

 Si on \'ecrit $X=(X_{+},X_{-},(X_{i})_{i\in I})$, on  d\'efinit la classe totale de conjugaison stable de $X_{+}$, $X_{-}$ et $X_{i}$ pour $i\in I$ de la m\^eme  fa\c{c}on  qu'en 1.2.  On d\'efinit la classe totale de conjugaison stable de $X$  comme le produit de ces classes. 

Consid\'erons maintenant un \'el\'ement $x'$ dans la classe totale de conjugaison stable de $x$. Il se  d\'ecompose en $x'=s'E(X')$. Les donn\'ees $n_{+}$, $\eta_{+}$, $n_{-}$, $\eta_{-}$, $I$ et $d_{i}$ ne changent pas quand on remplace $s$ par $s'$. Les formes $Q_{+}$, $Q_{-}$ et $Q_{i}$ pour $i\in I$ peuvent changer. Par exemple, soit $Q'_{+}$ la forme rempla\c{c}ant $Q_{+}$. Ou bien $Q'_{+}$ est isomorphe \`a $Q_{+}$, ou bien le couple $(Q'_{+},Q_{+})$ est \'egal, \`a l'ordre pr\`es, \`a l'un de nos couples $(Q_{iso},Q_{an})$ du paragraphe 2.1. On voit qu'\`a la classe de conjugaison stable de $X$ dans $\mathfrak{g}_{s}(F)$ correspond une classe de conjugaison stable dans $\mathfrak{g}_{s'}(F)$. L'\'el\'ement $X'$ appartient \`a cette classe, autrement dit, $X'$ appartient \`a la classe totale de conjugaison stable de $X$. Plus pr\'ecis\'ement, l'application $x'\mapsto X'$ se quotiente en une bijection entre l'ensemble des  classes de conjugaison contenues dans la classe totale de conjugaison stable de $x$ et l'ensemble  des  classes de conjugaison contenues dans la classe totale de conjugaison stable de $X$. 

Comme on vient de le dire, $Q_{+}$ est l'une des formes $Q_{iso}$ ou $Q_{an}$ de 2.1. Posons  $sgn^*(X_{+})=1$ si c'est $Q_{iso}$ et $sgn^*(X_{+})=-1$ si c'est $Q_{an}$. On d\'efinit de m\^eme $sgn^*(X_{-})$ et $sgn^*(X_{i})$ pour $i\in I$. On pose
$$sgn^*(X)=sgn^*(X_{+})sgn^*(X_{-})\prod_{i\in I}sgn^*(X_{i}).$$
On rappelle que, pour notre indice $\sharp$ fix\'e plus haut, on a d\'efini $sgn_{\sharp}=1$ si $\sharp=iso$ et $sgn_{\sharp}=-1$ si $\sharp=an$. Montrons que

(3) on a l'\'egalit\'e 
$$sgn_{\sharp}= (-1)^{d\, val_{F}(\eta_{-})}sgn^*(X ) .$$

Preuve. Dans chacun de nos espaces $V_{+}$, $V_{-}$ et $V_{i}$ pour $i\in I$, on fixe un r\'eseau presque autodual (cf. \cite{W5} 1.1) que l'on note $L_{+}$, $L_{-}$ et $L_{i}$. On note $L$ la somme directe de ces r\'eseaux. C'est un r\'eseau presque autodual de $V_{\sharp}$. On pose $l''_{+}=L_{+}^*/L_{+}$ et on d\'efinit de m\^eme $l''_{-}$,$l''_{i}$ et $l''$. Chacun de ces espaces est muni d'une forme quadratique dont on note les d\'eterminants normalis\'es $\eta''_{+}$, $\eta''_{-}$, $\eta''_{i}$ et $\eta''$ (rappelons que par exemple $\eta''$ est le produit du d\'eterminant habituel par $(-1)^{[dim_{{\mathbb F}_{q}}(l'')/2]}$). Par d\'efinition de $V_{\sharp}$, la dimension de $l''$ sur ${\mathbb F}_{q}$ est paire et on a $sgn_{\sharp}=sgn(\eta'')$. On v\'erifie facilement que, pour $i\in I$, la dimension de $l''_{i}$ est paire. D'apr\`es la d\'efinition de \cite{W5} 1.1, on a $sgn^*(X_{i})=sgn(\eta''_{i})$. Supposons $val_{F}(\eta_{+})$ paire. Le terme $val_{F}(\eta_{-})$ est \'egalement pair d'apr\`es (2). Alors les dimensions de $l''_{+}$ et $l''_{-}$ sont paires et on a $sgn^*(X_{+})=sgn(\eta''_{+})$, $sgn^*(X_{-})=sgn(\eta''_{-})$. Le d\'eterminant non normalis\'e de la forme sur $l''$ est le produit des d\'eterminants non normalis\'es des formes sur $l''_{+}$, $l''_{-}$ et $l''_{i}$ pour $i\in I$. Puisque toutes les dimensions sont paires, il en est de m\^eme des d\'eterminants normalis\'es: on a $\eta''=\eta''_{+}\eta''_{-}\prod_{i\in I}\eta''_{i}$. Avec les formules ci-dessus, on obtient 
$$sgn_{\sharp}=  sgn^*(X) $$
ce qui co\"{\i}ncide avec (3) puisque $val_{F}(\eta_{-})$ est paire.
 Supposons maintenant $val_{F}(\eta_{+})$ impaire, donc aussi $val_{F}(\eta_{-})$ impaire. Alors les dimensions de $l''_{+}$ et $l''_{-}$ sont impaires. D'apr\`es les d\'efinitions de \cite{W5} 1.1, on a $sgn^*(X_{+})=sgn(\eta'_{+})$ et $sgn^*(X_{-})=sgn(\eta'_{-})$, o\`u $\eta'_{+}$ et $\eta'_{-}$ sont les d\'eterminants normalis\'es des formes quadratiques sur $L_{+}/\varpi L_{+}^*$ et $L_{-}/\varpi L_{-}^*$. Parce que la dimension de $V_{+}$ est impaire et  que celle de $V_{-}$ est paire, on a $\eta'_{+}\eta''_{+}=\eta_{+}\varpi^{-val_{F}(\eta_{+})}$ et $\eta'_{-}\eta''_{-}=-\eta_{-}\varpi^{-val_{F}(\eta_{-})}$. Parce que $l''_{+}$ et $l''_{-}$ sont de dimension impaire, on voit qu'il se glisse un signe $-$ dans le produit des d\'eterminants normalis\'es. C'est-\`a-dire $\eta''=-\eta''_{+}\eta''_{-}\prod_{i\in I}\eta''_{i}$. D'o\`u
$$\eta''=\eta_{+}\eta_{-}\varpi^{-val_{F}(\eta_{+})-val_{F}(\eta_{-})}\eta'_{+}\eta'_{-}\prod_{i\in I}\eta''_{i}.$$
La somme $val_{F}(\eta_{+})+val_{F}(\eta_{-})$ est nulle d'apr\`es (2).  Avec les formules ci-dessus, on obtient
$$sgn_{\sharp}=  sgn(\eta_{+}\eta_{-})sgn^*(X) .$$
Mais $sgn(\eta_{+}\eta_{-})=(-1)^d$ d'apr\`es (2). D'o\`u encore l'\'egalit\'e (3) puisque $val_{F}(\eta_{-})$ est impaire. 
$\square$

Soit $(n_{1},n_{2})\in D(n)$, auquel est associ\'ee une donn\'ee endoscopique de $G_{iso}$ et $G_{an}$, dont le groupe endoscopique est $G_{n_{1},iso}\times G_{n_{2},iso}$. Soit $(x_{1},x_{2})\in G_{n_{1},iso}(F)\times G_{n_{2},iso}(F)$ un couple $n$-r\'egulier form\'e d'\'el\'ements elliptiques. Pour $\sharp=iso$ ou $an$, soit $x\in G_{\sharp}(F)$ un \'el\'ement de la classe totale de conjugaison stable correspondant \`a $(x_{1},x_{2})$.  On lui associe les donn\'ees $\eta_{+}$, $n_{+}$, $\eta_{-}$ etc... comme ci-dessus. Pour $j=1,2$, on pose $x_{j}=s_{j}E(\underline{X}_{j})$. On adopte la notation $\underline{X}_{j}$ plut\^ot que $X_{j}$ pour une raison qui appara\^{\i}tra plus loin. On d\'efinit de m\^eme les donn\'ees $\eta_{j,+}$, $n_{j,+}$, $\eta_{j,-}$  etc...   On v\'erifie que $n_{+}=n_{1,+}+n_{2,+}$, $\eta_{+}=\eta_{1,+}\eta_{2,+}$, $n_{-}=n_{1,-}+n_{2,-}$, $\eta_{-}=\eta_{1,-}\eta_{2,-}$, que $I=I_{1}\cup I_{2}$ et que, pour $i\in I$, $d_{i}=d_{1,i}+d_{2,i}$ (avec $d_{j,i}=0$ si $i\not\in I_{j})$. Le couple $(n_{1,+},n_{2,+})$, resp. le quadruplet $(n_{1,-},\eta_{1,-},n_{2,-},\eta_{2,-})$, resp. le couple $(d_{1,i},d_{2,i})$ pour $i\in I$, d\'efinit une donn\'ee endoscopique de $G_{+}$, resp. $G_{-}$, $G_{i}$. Avec les d\'efinitions que l'on a donn\'ees, le couple $(\underline{X}_{1,+},\underline{X}_{2,+})$ (par exemple) n'a pas de raison d'appartenir au groupe endoscopique associ\'e \`a $(n_{1,+},n_{2,+})$ car ce dernier est $G_{n_{1,+},iso}\times G_{n_{2,+},iso}$ tandis que $(\underline{X}_{1,+},\underline{X}_{2,+})$ appartient \`a un certain produit $G_{n_{1,+},\sharp_{1,+}}(F)\times G_{n_{2,+},\sharp_{2,+}}(F)$. Mais  seule la classe totale de conjugaison stable de $(\underline{X}_{1,+},\underline{X}_{2,+})$ interviendra dans la suite et on peut fixer un \'el\'ement $(X_{1,+},X_{2,+})$ de cette classe qui appartient \`a $G_{n_{1,+},iso}(F)\times G_{n_{2,+},iso}(F)$.  On d\'efinit de m\^eme $(X_{1,-},X_{2,-})$ et $(X_{1,i},X_{2,i})$ pour $i\in I$.  On pose $X_{1}=X_{1,+}\oplus X_{1,-}\oplus\oplus_{i\in I}X_{1,i}$ et on d\'efinit de m\^eme $X_{2}$. 

On dispose d'un facteur de transfert $\Delta((x_{1},x_{2}),x)$ associ\'e \`a la donn\'ee endoscopique d\'efinie par $(n_{1},n_{2})$. On dispose d'un facteur de transfert, que nous noterons simplement $\Delta_{+}((X_{1,+},X_{2,+}),X_{+})$ associ\'e \`a la donn\'ee endoscopique d\'efinie par $(n_{1,+},n_{2,+})$ et, de m\^eme, de facteurs de transfert $\Delta_{-}((X_{1,-},X_{2,-}),X_{-})$ et $\Delta_{i}((X_{1,i},X_{2,i}),X_{i})$. Ces trois derniers facteurs sont normalis\'es comme dans les paragraphes 2.1, 2.2 et 2.3.  On pose
$$\Delta((X_{1},X_{2}),X)=\Delta_{+}((X_{1,+},X_{2,+}),X_{+})\Delta_{-}((X_{1,-},X_{2,-}),X_{-})\prod_{i\in I}\Delta_{i}((X_{1,i},X_{2,i}),X_{i})$$
et
$$d_{2}=\sum_{i\in I_{2}}d_{2,i}.$$

\ass{Lemme}{On a l'\'egalit\'e
$$\Delta((x_{1},x_{2}),x)=(-1)^{d_{2}val_{F}(\eta_{-})} \Delta((X_{1},X_{2}),X).$$}

Preuve.   Notons $K_{+}$  l'ensemble des orbites de l'action de $Gal(\bar{F}/F)$ dans l'ensemble des valeurs propres de $X_{+}$ diff\'erentes de $0$ (la valeur propre $0$ intervient avec multiplicit\'e $1$), $K_{-}$  l'ensemble des orbites de l'action de $Gal(\bar{F}/F)$ dans l'ensemble des valeurs propres de $X_{-}$ et, pour $i\in I$, $K_{i}$  l'ensemble des orbites de l'action de $Gal(\bar{F}/E_{i})$ dans l'ensemble des valeurs propres de $X_{i}$, vu comme un endomorphisme $E_{i}$-lin\'eaire de $V_{i}$.  Pour $k\in K_{+}$, resp. $k\in K_{-}$, $k\in K_{i}$, on fixe $\Xi_{k}\in k$. On pose $E(\Xi_{k})=(1+\Xi_{k}/2)(1-\Xi_{k}/2)^{-1}$. 
Pour $k\in K_{+}\cup K_{-}$, on pose $F_{k}=F[\Xi_{k}]$. Il existe une sous-extension $F^{\natural}_{k}$ de $F$ telle que $[F_{k}:F^{\natural}_{k}]=2$ et que $trace_{F_{k}/F^{\natural}_{k}}(\Xi_{k})=0$.  Pour $i\in I$ et $k\in K_{i}$, on pose $F_{k}=E_{i}[\Xi_{k}]$. Il existe une extension $F^{\natural}_{k}$ de $E^{\natural}_{i}$ de sorte que $F_{k}$ soit le compos\'e des extensions $E_{i}$ et  $F^{\natural}_{k}$ de $E^{\natural}_{i}$ et que $trace_{F_{k}/F^{\natural}_{k}}(\Xi_{k})=0$. On note $K$ l'union disjointe de $K_{+}$, $K_{-}$ et des $K_{i}$ pour $i\in I$. L'ensemble des valeurs propres de $x$ diff\'erentes de $1$ est l'ensemble des conjugu\'es par $Gal(\bar{F}/F)$ des \'el\'ements

$\xi_{k}=E(\Xi_{k})$ pour $k\in K_{+}$; $\xi_{k}=-E(\Xi_{k})$ pour $k\in K_{-}$; $\xi_{k}=s_{i}E(\Xi_{k})$ pour $i\in I$ et $k\in K_{i}$. 

Pour $k\in K$, notons $V_{k}$ la somme des espaces propres de $x$ associ\'es \`a des conjugu\'es de $\xi_{k}$. On peut identifier $V_{k}$ \`a $F_{k}$ de sorte que l'action de $x$ dans $V_{k}$ s'identifie \`a la multiplication par $\xi_{k}$. Il existe un \'el\'ement $c_{k}\in F^{\natural}_{k}$ de sorte que la restriction de $Q_{\sharp}$ \`a $V_{k}$ s'identifie \`a la forme quadratique
$$(v,v')\mapsto [F_{k}:F]^{-1}trace_{F_{k}/F}(c_{k}\tau_{k}(v')v)$$
sur $F_{k}$, o\`u $\tau_{k}$ est l'\'el\'ement non trivial de $Gal(F_{k}/F^{\natural}_{k})$. Seule compte la classe de $c_{k}$ modulo les normes de l'extension $F_{k}/F_{k}^{\natural}$. On note $sgn_{k}$ le caract\`ere quadratique de $F_{k}^{\natural,\times}$ associ\'e \`a  cette extension. 
On note  $P_{k}$ le polyn\^ome caract\'eristique de $\xi_{k}$ sur $F$. On note $P$ le produit de ces polyn\^omes. Pour tout polyn\^ome $R$, on note $R'$ son polyn\^ome d\'eriv\'e. On pose
$$C_{k}=(-1)^{n+1}2[F_{k}:F]c_{k}P'(\xi_{k})P(-1)\xi_{k}^{1-n}(1+\xi_{k})(\xi_{k}-1)^{-1}.$$

On d\'efinit $K_{1}$ et $K_{2}$ comme on a d\'efini $K$. L'ensemble $K$ est l'union disjointe de $K_{1}$ et $K_{2}$. L'union est disjointe car $(x_{1},x_{2})$ est $n$-r\'egulier. D'apr\`es \cite{W1} proposition 1.10, on a l'\'egalit\'e
$$(4) \qquad \Delta((x_{1},x_{2}),x)=\prod_{k\in K_{2}}sgn_{k}(C_{k}).$$

{\bf Remarque.} Dans la formule de \cite{W1} d\'efinissant $C_{k}$, il n'y a pas le facteur $[F_{k}:F]$. Mais le terme $c_{k}$ y est d\'efini diff\'eremment: c'est  notre $c_{k}$ multipli\'e par $[F_{k}:F]^{-1}$. 

\bigskip

Pour  $k\in K_{+}\cup K_{-}$, on note $\mathfrak{P}_{k}$ le polyn\^ome caract\'eristique de $\Xi_{k}$ sur $F$. On pose

$\mathfrak{C}_{k}= (-1)^{n_{+}}\eta_{+}[F_{k}:F]c_{k}\Xi_{k}\mathfrak{P}_{k}'(\Xi_{k})\prod_{k'\in K_{+},k'\not=k}\mathfrak{P}_{k'}(\Xi_{k})$, si $k\in K_{+}$;

$\mathfrak{C}_{k}=(-1)^{n_{-}} \eta_{-}[F_{k}:F]c_{k}\Xi_{k}^{-1}\mathfrak{P}_{k}'(\Xi_{k})\prod_{k'\in K_{-},k'\not=k}\mathfrak{P}_{k'}(\Xi_{k})$, si $k\in K_{-}$.

Pour $i\in I$ et $k\in K_{i}$, on note $\mathfrak{P}_{i,k}$ le polyn\^ome caract\'eristique de $\Xi_{k}$ sur $E_{i}$. On fixe un \'el\'ement $\eta_{i}\in E_{i}$ qui est une unit\'e et v\'erifie $\tau_{i}(\eta_{i})=(-1)^{d_{i}+1}\eta_{i}$, o\`u $\tau_{i}$ est l'unique \'el\'ement non trivial de $Gal(E_{i}/E_{i}^{\natural})$. On  pose

$\mathfrak{C}_{k}= \eta_{i}[F_{k}:E_{i}]c_{k} \mathfrak{P}_{i,k}'(\Xi_{k})\prod_{k'\in K_{i},k'\not=k}\mathfrak{P}_{i,k'}(\Xi_{k})$.

D'apr\`es \cite{W3} lemme X.7, on a l'\'egalit\'e
 $$(5) \qquad  \Delta((X_{1},X_{2}),X)=\prod_{k\in K_{2}}sgn_{k}(\mathfrak{C}_{k}).$$
 On va d\'emontrer les propri\'et\'es suivantes:
 
 (6) pour $k\in K_{+}\cup K_{-}$, $sgn_{k}(C_{k})=sgn_{k}(\mathfrak{C}_{k})$;
  
 (7) pour $i\in I$ et $k\in K_{i}$, $sgn_{k}(C_{k})=(-1)^{[F_{k}:E_{i}]val_{F}(\eta_{-})} sgn_{k}(\mathfrak{C}_{k})$. 
 
  En les admettant, on voit que le rapport entre les membres de droite de (4) et (5) est $(-1)^{D\,val_{F}(\eta_{-})}$, o\`u $D=\sum_{i\in I}\sum_{k\in K_{i}\cap K_{2}}[F_{k}:E_{i}]$. La somme int\'erieure en $k$ vaut $d_{2,i}$, donc $D=d_{2}$. Alors les \'egalit\'es (4) et (5) impliquent celle de l'\'enonc\'e.
  
  Pour d\'emontrer (6) et (7), on a besoin  de quelques ingr\'edients. Soit $Z\in \bar{F}$ et $\left(\begin{array}{cc}a&b\\c&d\\ \end{array}\right)\in GL(2;F)$. Supposons $cZ+d\not=0$ et posons $z=\frac{az+b}{cz+d}$. Notons $P_{z}(T)$ et $P_{Z}(T)$ les polyn\^omes caract\'eristiques de $z$ et $Z$ et notons $m$ leur degr\'e. Alors on a les \'egalit\'es
  $$(8) \qquad (cT+d)^mP_{z}(\frac{aT+b}{cT+d})=c^mP_{z}(\frac{a}{c})P_{Z}(T);$$
  $$(9) \qquad (ad-bc)(cZ+d)^{m-2}P'_{z}(z)=c^mP_{z}(\frac{a}{c})P'_{Z}(T).$$
  Cf. \cite{Li} lemme 7.4.1. D'autre part, pour $k\in K$, on calcule facilement le d\'eterminant non normalis\'e  $det_{k}\in F^{\times}/F^{\times2}$ de la forme quadratique d\'efinie plus haut sur $V_{k}$: on a
  $$(10)\qquad det_{k}=norme_{F_{k}/F}(\Xi_{k}).$$
  
 Fixons une extension galoisienne  finie $F'$ de $F$ contenant tous les corps $F_{k}$. Notons $\mathfrak{o}'_{1}$ le groupe multiplicatif des unit\'es de $F'$ congrues \`a $1$ modulo l'id\'eal maximal $\mathfrak{p}'$ de l'anneau des entiers de $F'$. Pour $k\in K_{+}\cup K_{-}$, introduisons la relation d'\'equivalence dans $F^{'\times}$: $x\equiv_{k} y$ si et seulement s'il existe $x'\in \mathfrak{o}'_{1}$ et $y'\in F_{k}^{\times}$ tels que $xy^{-1}=x'\,norme_{F_{k}/F_{k}^{\natural}}(y')$.  On peut remplacer la relation (6) par
 
 (11) pour $k\in K_{+}\cup K_{-}$, $C_{k}\equiv_{k}\mathfrak{C}_{k}$.
 
 En effet, si cette relation est v\'erifi\'ee, il existe $x\in \mathfrak{o}'_{1}$ et $y\in F_{k}^{\times}$ tels que $C_{k}=x\,norme_{F_{k}/F_{k}^{\natural}}(y)\mathfrak{C}_{k}$. On sait que les termes $C_{k}$ et $\mathfrak{C}_{k}$ appartiennent \`a $F_{k}^{\natural,\times}$. Donc $x\in \mathfrak{o}'_{1}\cap F_{k}^{\natural,\times}$. Parce que $p$ est grand, le caract\`ere $sgn_{k}$ est mod\'er\'ement ramifi\'e, donc $sgn_{k}(x)=1$. On a aussi $sgn_{k}(norme_{F_{k}/F_{k}^{\natural}}(y))=1$. L'\'egalit\'e $sgn_{k}(C_{k})=sgn_{k}(\mathfrak{C}_{k})$ s'ensuit.
 
 Soit $k\in K_{+}$. Pour $k'\in K_{-}$ ou $k'\in K_{i}$ pour $i\in I$, la contribution de $k'$ \`a $C_{k}$ est $P_{k'}(\xi_{k})P_{k'}(-1)$. On a $\xi_{k}\in \mathfrak{o}'_{1}$ tandis que les racines de $P_{k'}$ n'appartiennent pas \`a ce groupe. On en  d\'eduit $P_{k'}(\xi_{k})\equiv_{k}P_{k'}(1)\equiv_{k}P_{k'}(1)^{-1}$. Le produit $P_{k'}(1)^{-1}P_{k'}(-1)$ est \'egal \`a $norme_{F_{k'}/F}(1-\xi_{k'})^{-1}(-1-\xi_{k'})$. Parce que $norme_{F_{k'}/F_{k'}^{\natural}}(\xi_{k'})=1$, on voit que $(1-\xi_{k'})^{-1}(-1-\xi_{k'})$ est un \'el\'ement de $F_{k'}$ dont la trace dans $F_{k'}^{\natural}$ est nulle. Il existe donc $x\in F_{k'}^{\natural,\times}$ tel que $(1-\xi_{k'})^{-1}(-1-\xi_{k'})=x\Xi_{k'}$. D'apr\`es (10), on a donc $$P_{k'}(1)^{-1}P_{k'}(-1)\equiv_{k}norme_{F_{k'}/F}(x)det_{k'}\equiv_{k}norme_{F_{k'}^{\natural}/F}(x)^2det_{k'}\equiv_{k}det_{k'}.$$
 
 {\bf Remarque.} Le terme $det_{k'}$ n'est d\'efini que modulo $F^{\times2}$. Disons qu'on en choisit un repr\'esentant dans $F^{\times}$. La classe d'\'equivalence  de ce repr\'esentant pour la relation $\equiv_{k}$ est bien d\'efinie. 
 
 \bigskip
 
 Soit $k'\in K_{+}$ avec $k'\not=k$. La contribution de $k'$ \`a $C_{k}$ est $P_{k'}(\xi_{k})P_{k'}(-1)$. On utilise (8) avec $\left(\begin{array}{cc}a&b\\c&d\\ \end{array}\right)=\left(\begin{array}{cc}-\frac{1}{2}&1\\-\frac{1}{2}&1\\ \end{array}\right)$, $Z=\Xi_{k'}$ et $T=\Xi_{k}$. L'entier $m=[F_{k'}:F]$ est pair. De plus, $\Xi_{k}\in \mathfrak{p}'$, donc $c\Xi_{k}+d\in \mathfrak{o}'_{1}$. On obtient $P_{k'}(\xi_{k})\equiv_{k}P_{k'}(-1)\mathfrak{P}_{k'}(\Xi_{k})$. D'o\`u
 $$P_{k'}(\xi_{k})P_{k'}(-1)\equiv_{k}\mathfrak{P}_{k'}(\Xi_{k}).$$
 Un calcul analogue s'applique au terme $P'_{k}(\xi_{k})P_{k}(-1)$, en utilisant cette fois la relation (9). On obtient
 $$P'_{k}(\xi_{k})P_{k}(-1)\equiv_{k}\mathfrak{P}'_{k}(\Xi_{k}).$$
Evidemment $\xi_{k}^{1-n}\equiv_{k}1$ et $1+\xi_{k}\equiv_{k}2$. On calcule
$$(\xi_{k}-1)^{-1}=(1-\Xi_{k}/2)\Xi_{k}^{-1}\equiv_{k}\Xi_{k}^{-1}=-\Xi_{k}(-\Xi_{k}^2)^{-1}=-\Xi_{k}norme_{F_{k}/F_{k}^{\natural}}(\Xi_{k})^{-1}\equiv_{k}-\Xi_{k}.$$
En rassemblant ces calculs, on obtient
$$C_{k}\equiv_{k} (-1)^{n+n_{+}}\eta_{+}(\prod_{k'}det_{k'})\mathfrak{C}_{k},$$
o\`u le produit porte sur les $k'\in K_{-}$ et $k'\in K_{i}$ pour $i\in I$. Le produit des $det_{k'}$ sur ces $k'$ est le d\'eterminant de la restriction de $Q_{\sharp}$ \`a la somme de $V_{-}$ et des $V_{i}$ pour $i\in I$. Il est \'egal au d\'eterminant de $Q_{\sharp}$ divis\'e par le d\'eterminant de la restriction de $Q_{\sharp}$ \`a $V_{+}$. On a fix\'e le d\'eterminant de $Q_{\sharp}$ en \cite{W5}  1.1:  c'est $(-1)^n$. Celui de  la restriction de $Q_{\sharp}$ \`a $V_{+}$ est $(-1)^{n_{+}}\eta_{+}$. L'\'equivalence ci-dessus entra\^{\i}ne alors (11). 
 
Soit $k\in K_{-}$. Pour $k'\in K_{+}$ ou $k'\in K_{i}$ pour $i\in I$, la contribution de $k'$ \`a $ C_{k}$ est $P_{k'}(\xi_{k})P_{k'}(-1)$. On a $\xi_{k}\in -\mathfrak{o}'_{1}$ et les racines de $P_{k'}$ ne sont pas congrues \`a $-1$ modulo $\mathfrak{p}'$. Donc $P_{k'}(\xi_{k})\equiv_{k}P_{k'}(-1)$, puis
$$P_{k'}(\xi_{k})P_{k'}(-1)\equiv_{k}P_{k'}(-1)^2\equiv_{k}1.$$
Pour $k'\in K_{-}$ avec $k'\not=k$, la contribution de $k'$ \`a $ C_{k}$ est $P_{k'}(\xi_{k})P_{k'}(-1)$. On utilise (8) avec  $\left(\begin{array}{cc}a&b\\c&d\\ \end{array}\right)=\left(\begin{array}{cc}-\frac{1}{2}&-1\\-\frac{1}{2}&1\\ \end{array}\right)$, $Z=\Xi_{k'}$ et $T=\Xi_{k}$. De nouveau, $m$ est pair et $c\Xi_{k}+d\equiv_{k}1$. D'o\`u $P_{k'}(\xi_{k})\equiv_{k}P_{k'}(1)\mathfrak{P}_{k'}(\Xi_{k})$. Comme plus haut, on a
$$P_{k'}(1)P_{k'}(-1)\equiv_{k}P_{k'}(1)^{-1}P_{k'}(1)\equiv_{k}det_{k'}.$$
D'o\`u
$$P_{k'}(\xi_{k})P_{k'}(-1)\equiv_{k}det_{k'}\mathfrak{P}_{k'}(\Xi_{k}).$$
Un m\^eme calcul s'applique au terme $P'_{k}(\xi_{k})P_{k}(-1)$, en utilisant cette fois la relation (9). Ici se glisse le d\'eterminant $ad-bc$ qui vaut $-1$. D'o\`u
  $$P'_{k}(\xi_{k})P_{k}(-1)\equiv_{k}det_{k}\mathfrak{P}'_{k}(\Xi_{k}).$$
Evidemment $\xi_{k}^{1-n}\equiv_{k}(-1)^{n-1}$ et $(\xi_{k}-1)^{-1}\equiv_{k}-2$. On a
$$1+\xi_{k}=-\frac{\Xi_{k}}{1-\frac{\Xi_{k}}{2}}\equiv_{k}-\Xi_{k}\equiv_{k} -\Xi_{k}norme_{F_{k}/F_{k}^{\natural}}(\Xi_{k})^{-1}\equiv_{k}\Xi_{k}^{-1}.$$
En rassemblant ces calculs, on obtient 
$$C_{k}\equiv_{k}(-1)^{n_{-}}\eta_{-}(\prod_{k'\in K_{-}}det_{k'})\mathfrak{C}_{k}.$$
Le produit intervenant ici est le d\'eterminant de la restriction de $Q_{\sharp}$ \`a $V_{-}$, c'est-\`a-dire $(-1)^{n_{-}}\eta_{-}$. D'o\`u encore (11).

Soient $i\in I$ et $k\in K_{i}$. Cette fois, on d\'efinit dans $F^{'\times}$ l'\'equivalence $x\equiv_{k}y$ si et seulement s'il existe $x'\in F_{k}^{\times}$ tel que $ xy^{-1}norme_{F_{k}/F_{k}^{\natural}}(x')$ soit une unit\'e.  Soit $k'\in K_{+}$ ou $k'\in K_{i'}$ pour un $i'\in I$ avec $i'\not=i$. Les racines du polyn\^ome $P_{k'}$ ne sont congrues ni \`a $-1$, ni \`a $\xi_{k}$ modulo $\mathfrak{p}'$. On en d\'eduit $P_{k'}(\xi_{k})P_{k'}(-1)\equiv_{k}1$. Soit $k'\in K_{-}$. On a  
$$P_{k'}(\xi_{k})P_{k'}(-1)=P_{k'}(\xi_{k})P_{k'}(1)P_{k'}(1)^{-1}P_{k'}(-1).$$
On a encore $P_{k'}(\xi_{k})P_{k'}(1)\equiv_{k}1$ et, par un calcul d\'ej\`a fait, $P_{k'}(1)^{-1}P_{k'}(-1)\equiv_{k}det_{k'}$. 
Soit $k'\in K_{i}$ avec $k'\not=k$. Comme pr\'ec\'edemment,  $P_{k'}(-1)\equiv_{k}1$. On a
$$P_{k'}(\xi_{k})=\prod_{\sigma\in Gal(F_{k'}/F)}(\xi_{k}-\sigma(\xi_{k'})).$$
Pour $\sigma\in Gal(F_{k'}/F)$, $\xi_{k}-\sigma(\xi_{k'})$ est congru modulo $\mathfrak{p}'$ \`a $s_{i}-\sigma(s_{i})$. Si $\sigma\not\in Gal(F_{k}/E_{i})$, ce terme est une unit\'e. Si $\sigma\in Gal(F_{k}/E_{i})$, on a $\xi_{k}-\sigma(\xi_{k'})=s_{i}(E(\Xi_{k})-\sigma(E(\Xi_{k'})))$. 
Notons $P_{i,k'}$ le polyn\^ome caract\'eristique de $E(\Xi_{k'})$ sur $E_{i}$. On obtient
$$P_{k'}(\xi_{k})P_{k'}(-1)\equiv_{k}P_{i,k'}(E(\Xi_{k})).$$
On calcule ce terme gr\^ace \`a (8) o\`u l'on remplace le corps $F$ par $E_{i}$. On prend $\left(\begin{array}{cc}a&b\\c&d\\ \end{array}\right)=\left(\begin{array}{cc}\frac{1}{2}&1\\-\frac{1}{2}&1\\ \end{array}\right)$, $Z=\Xi_{k'}$, $T=\Xi_{k}$. On obtient
$$P_{i,k'}(E(\Xi_{k}))\equiv_{k}P_{i,k'}(-1)\mathfrak{P}_{i,k'}(\Xi_{k})\equiv_{k}\mathfrak{P}_{i,k'}(\Xi_{k}),$$
d'o\`u
$$P_{k'}(\xi_{k})P_{k'}(-1)\equiv_{k}\mathfrak{P}_{i,k'}(\Xi_{k}).$$
Un m\^eme calcul s'applique \`a $P'_{k}(\xi_{k})P_{k}(-1)$, en utilisant cette fois la relation (9). D'o\`u
$$P'_{k}(\xi_{k})P_{k}(-1)\equiv_{k}\mathfrak{P}'_{i,k}(\Xi_{k}).$$
Evidemment $(-1)^{n+1}2\xi_{k}^{1-n}(1+\xi_{k})(\xi_{k}-1)^{-1}\equiv_{k}1$ et $\eta_{i}\equiv_{k}1$ (tous les termes sont des unit\'es). En rassemblant ces calculs, on obtient
$$C_{k}\equiv_{k}(\prod_{k\in K_{-}}det_{k})\mathfrak{C}_{k}.$$
Le produit intervenant ici vaut $(-1)^{n_{-}}\eta_{-}$ (comme plus haut, on fixe ici un repr\'esentant de $\eta^-$ dans $F^{\times}$). Il existe donc $x\in F_{k}^{\times}$ et une unit\'e $y$ de $F^{'\times}$  tels que $C_{k}=y\,norme_{F_{k}/F_{k}^{\natural}}(x)\eta_{-}\mathfrak{C}_{k}$. N\'ecessairement, $y$ appartient \`a $F_{k}^{\natural,\times}$. L'extension $F_{k}/F_{k}^{\natural}$ est non ramifi\'ee puisque $F_{k}$ est le compos\'e de $E_{i}$ et de $F_{k}^{\natural}$ sur $E_{i}^{\natural}$. Donc $sgn_{k}(y)=1$, puis $sgn_{k}(C_{k})=sgn_{k}(\eta_{-})sgn_{k}(\mathfrak{C}_{k})$.  Notons $val_{F_{k}^{\natural}}$ la valuation usuelle de $F_{k}^{\natural}$. On a $sgn_{k}(x)=(-1)^{val_{F_{k}^{\natural}}(x)}$ pour tout $x\in F_{k}^{\natural, \times}$. En notant $e(F_{k}^{\natural}/F)$ l'indice de ramification de l'extension $F_{k}^{\natural}/F$, on a $val_{F_{k}^{\natural}}(\eta_{-})=e(F_{k}^{\natural}/F)val_{F}(\eta_{-})$.  Puisque $E_{i}^{\natural}/F$ est non ramifi\'ee, on a 
$$e(F_{k}^{\natural}/F)=e(F_{k}^{\natural}/E_{i}^{\natural})=[F_{k}^{\natural}:E_{i}^{\natural}]f(F_{k}^{\natural}/E_{i}^{\natural})^{-1}=[F_{k}:E_{i}]f(F_{k}^{\natural}/E_{i}^{\natural})^{-1},$$
o\`u $f(F_{k}^{\natural}/E_{i}^{\natural})$ est le degr\'e de l'extension r\'esiduelle. Mais $E_{i}$ est l'unique extension quadratique non ramifi\'ee de $E_{i}^{\natural}$ et elle n'est pas contenue dans $F_{k}^{\natural}$ puisque $F_{k}$ est compos\'e de $E_{i}$ et de $F_{k}^{\natural}$.  Donc $f(F_{k}^{\natural}/E_{i}^{\natural})$ est impair et $e(F_{k}^{\natural}/F)$ est de la m\^eme parit\'e que $[F_{k}:E_{i}]$. D'o\`u
$$sgn_{k}(\eta_{-})=(-1)^{e(F_{k}^{\natural}/F)val_{F}(\eta_{-})}=(-1)^{[F_{k}:E_{i}]val_{F}(\eta_{-})}.$$
D'o\`u l'\'egalit\'e (7), ce qui ach\`eve la d\'emonstration. $\square$

\bigskip

\subsection{Calcul d'int\'egrales orbitales}
Fixons $(r',r'',N',N'')\in \Gamma$, $w'\in W_{N'}$ et $w''\in W_{N''}$. Comme en 2.1, on note $\varphi_{w'}$ et $\varphi_{w''}$ les fonctions caract\'eristiques des classes de conjugaison de $w'$ et $w''$. On suppose que ces classes sont param\'etr\'ees par des couples de partitions de la forme $(\emptyset,\beta')$ et $(\emptyset,\beta'')$. On pose $f=\Psi\circ k\circ \rho\iota(\varphi_{w'}\otimes \varphi_{w''})$, cf. 1.2. On d\'efinit un couple d'entiers $(r'_{+},r'_{-})$ par les \'egalit\'es

$(r'_{+},r'_{-})=(r'+1,r')$, si $r'\equiv r''\,\,mod\,\,2{\mathbb Z}$;

$(r'_{+},r'_{-})=(r',r'+1)$, si $r'\not\equiv r''\,\,mod\,\,2{\mathbb Z}$.

Soit $\sharp=iso$ ou $an$ et soit $x\in G_{\sharp}(F)$ un \'el\'ement elliptique fortement r\'egulier. 
 On l'\'ecrit $x=sE(X)$ et on associe \`a $s$ les donn\'ees du paragraphe pr\'ec\'edent. En particulier, on a
 $$G_{s}=G_{+}\times G_{-}\times \prod_{i\in I}G_{i}.$$
 
 On consid\`ere les hypoth\`eses
 
 (1)(a) $val_{F}(\eta_{-})\equiv r''\,\,mod\,\,2{\mathbb Z}$;
 
 (1)(b) $2n_{+}+1\geq r^{'2}_{+}+r^{''2}$, $2n_{-}\geq r^{'2}_{-}+r^{''2}$.

 Remarquons que, d'apr\`es 3.1(2), la condition (1)(a) \'equivaut  \`a  $val_{F}(\eta_{+})\equiv r''\,\,mod\,\,2{\mathbb Z}$.

 Supposons v\'erifi\'ees ces conditions (1)(a) et (b). On pose
 $$N_{+}=n_{+}-(r^{'2}_{+}+r^{''2}-1)/2,\,\, N_{-}=n_{-}-(r^{'2}_{-}+r^{''2})/2.$$
 Notons $D$ l'ensemble des familles ${\bf d}=(N'_{+},N'_{-},N''_{+},N''_{-},(d'_{i},d''_{i})_{i\in I})$ d'entiers positifs ou nuls v\'erifiant les conditions
 $$N'_{+}+N''_{+}=N_{+},\,\,N'_{-}+N''_{-}=N_{-}$$;
 
 $ d'_{i}+d''_{i}=d_{i}$ pour tout $i\in I$;
 
 $N'_{+}+N'_{-}+\sum_{i\in I}d'_{i}f_{i}=N'$, $N''_{+}+N''_{-}+\sum_{i\in I}d''_{i}f_{i}=N''$,
 
 \noindent o\`u on a pos\'e $f_{i}=[E_{i}^{\natural}:F]$. Pour une telle famille, posons
 $$W'({\bf d})=W_{N'_{+}}\times W_{N'_{-}}\times \prod_{i\in I}\mathfrak{S}_{d'_{i}}.$$
 Consid\'erons l'ensemble des \'el\'ements ${\bf v}'=(v'_{+},v'_{-},(v'_{i})_{i\in I})\in W'({\bf d})$ v\'erifiant les conditions suivantes:
 
 les classes de conjugaison de $v'_{+}$ et $v'_{-}$ sont param\'etr\'ees par des couples de partitions $(\emptyset,\beta'_{+})$ et $(\emptyset,\beta'_{-})$;
 
 pour $i\in I$, la classe de conjugaison de $v'_{i}$ est param\'etr\'ee par une partition $\beta'_{i}$ dont tous les termes non nuls sont impairs; on note $f_{i}\beta'_{i}$ la partition dont les termes sont ceux de $\beta'_{i}$ multipli\'es par $f_{i}$;
 
 $\beta'=\beta'_{+}\cup \beta'_{-}\cup\bigcup_{i\in I}f_{i}\beta'_{i}$.
 
 Cet ensemble est r\'eunion de classes de conjugaison par $W'({\bf d})$ et on fixe un ensemble de repr\'esentants ${\cal V}'({\bf d})$ de ces classes. Remarquons que
 
 (2) pour ${\bf v}'=(v'_{+},v'_{-},(v_{i})_{i\in I})\in W'({\bf d})$, on a l'\'egalit\'e
 $$sgn_{CD}(w')=sgn_{CD}(v'_{+})sgn_{CD}(v'_{-})(-1)^{\sum_{i\in I}d'_{i}}.$$
 
 En effet, $sgn_{CD}(w')=(-1)^{l(\beta')}$, o\`u $l(\beta')$ est le nombre de termes non nuls de $\beta'$. On a $l(\beta')=l(\beta'_{+})+l(\beta'_{-})+\sum_{i\in I}l(\beta'_{i})$. Pour $i\in I$,  $\beta'_{i}$ est une partition de $d'_{i}$ dont tous les termes non nuls   sont impairs. Donc $l(\beta'_{i})\equiv d'_{i}\,\,mod\,\,2{\mathbb Z}$. L'assertion (2) en r\'esulte.
 
 On pose les m\^emes d\'efinitions en rempla\c{c}ant les exposants $'$ par $''$. On pose ${\cal V}({\bf d})={\cal V}'({\bf d})\times {\cal V}''({\bf d})$. 
 
 Soient  ${\bf d}=(N'_{+},N'_{-},N''_{+},N''_{-},(d'_{i},d''_{i})_{i\in I})\in D$ et ${\bf v}=(v'_{+},v'_{-},(v'_{i})_{i\in I}  ;v''_{+},v''_{-},(v''_{i})_{i\in I})\in {\cal V}({\bf d})$.
 Supposons $n_{+}\geq1$. Appliquons la construction de 2.1 \`a l'entier $n_{+}$, \`a  l'\'el\'ement $\eta_{+}$ et au quadruplet $(r'_{+},\vert r''\vert ,N'_{+},N''_{+})$. Les hypoth\`eses de ce paragraphe sont v\'erifi\'ees d'apr\`es (1)(a).  On d\'efinit deux  fonctions $f^0_{+}={\cal Q}_{r'_{+},\vert r''\vert }^{Lie}\circ\rho_{N_{+}}^*\circ\iota_{N'_{+},N''_{+}}(\varphi_{v'_{+}}\otimes \varphi_{v''_{+}})$ et $f^{1}_{+} ={\cal Q}_{r'_{+},\vert r''\vert }^{Lie}\circ\rho_{N_{+}}^*\circ\iota_{N''_{+},N'_{+}}(\varphi_{v''_{+}}\otimes \varphi_{v'_{+}})$. Elles vivent sur deux alg\`ebres de Lie dont l'une est l'ag\`ebre de Lie $\mathfrak{g}_{+}$ de la premi\`ere composante $G_{+}$ de $G_{s}$. En particulier les int\'egrales orbitales  $J(X_{+},f^0_{+})$ et $J(X_{+},f^{1}_{+})$ sont bien d\'efinies.  Supposons $n_{-}\geq1$. On applique cette fois la construction de 2.2 et on d\'efinit  deux  fonctions $f^0_{-}={\cal Q}_{r'_{-},\vert r''\vert }^{Lie}\circ\rho_{N_{-}}^*\circ\iota_{N'_{-},N''_{-}}(\varphi_{v'_{-}}\otimes \varphi_{v''_{-}})$ et $f^{1}_{-} ={\cal Q}_{r'_{-},\vert r''\vert }^{Lie}\circ\rho_{N_{-}}^*\circ\iota_{N''_{-},N'_{-}}(\varphi_{v''_{-}}\otimes \varphi_{v'_{-}})$. Les int\'egrales orbitales $J(X_{-},f^0_{-})$ et $J(X_{-},f^{1}_{-})$ sont bien d\'efinies. Enfin, pour $i\in I$ on utilise les constructions de 2.3. On d\'efinit les fonctions $f^0_{i}={\cal Q}(d'_{i},d''_{i})^{Lie}\circ \rho_{i}^*\circ \iota_{d'_{i},d''_{i}}(\varphi_{v'_{i}}\otimes \varphi_{v''_{i}})$ et 
 $f^{1}_{i}={\cal Q}(d'_{i},d''_{i})^{Lie}\circ \rho_{i}^*\circ \iota_{d''_{i},d'_{i}}(\varphi_{v''_{i}}\otimes \varphi_{v'_{i}})$. Les int\'egrales orbitales $J(X_{i},f^0_{i})$ et $J(X_{i},f^{1}_{i})$ sont bien d\'efinies. On pose $f^0[{\bf d},{\bf v}]=f^0_{+}\otimes f^0_{-}\otimes \otimes_{i\in I}f^0_{i}$ et 
 $$J(X,f^0[{\bf d},{\bf v}])=J(X_{+},f^0_{+})J(X_{-},f^0_{-})\prod_{i\in I}J(X_{i},f^0_{i}).$$
 Dans ces formules, les termes index\'es par $+$ ou $-$ disparaissent si $n_{+}=0$ ou $n_{-}=0$. On d\'efinit de m\^eme $f^{1}[{\bf d},{\bf v}]$ et $J(X,f^{1}[{\bf d},{\bf v}])$.
 On pose
 
 $b=0$ si $r''>0$ ou si $r''=0$ et $r'$ est pair;
 
 $ b=1$ si $r''<0$ ou si $r''=0$ et $r'$ est impair.
 
 Posons
 $$c(r',r'')=(-1)^{n+r''}sgn(-1)^{(r^{'2}-r')/2+(r^{''2}-\vert r''\vert )/2},$$
 $$c_{\sharp}(r',r'',w',w'')=\left\lbrace\begin{array}{cc}1,&si\,\,0<r''\leq r'\,\,ou \,\,r''=0\,\,et \,\,r'\,\,est\,\,pair,\\ sgn_{CD}(w''),&si\,\,r'<r'',\\ sgn_{\sharp}&si\,\,-r'\leq r''<0\,\,ou\,\,r''=0\,\,et\,\,r'\,\,est\,\,impair,\\ sgn_{\sharp}sgn_{CD}(w'),&si \,\,r''<-r'.\\ \end{array}\right.$$

\ass{Proposition}{(i) Si les hypoth\`eses (1)(a) et (1)(b) ne sont pas v\'erifi\'ees, $J(x,f)=0$.

(ii) Supposons ces hypoth\`eses v\'erifi\'ees. Alors on a l'\'egalit\'e
$$J(x,f)=c(r',r'')c_{\sharp}(r',r'',w',w'')\sum_{{\bf d}\in D}\sum_{{\bf v}\in {\cal V}({\bf d}) }J(X,f^{b}[{\bf d},{\bf v}]).$$}

C'est la proposition 3.19 de \cite{MW}.

\bigskip

\subsection{D\'emonstration du (ii) de la proposition 1.2}
On conserve les donn\'ees $(r',r'',N',N'')\in \Gamma$, $w'\in W_{N'}$ et $w''\in W_{N''}$. On  d\'efinit $f$ comme dans le paragraphe pr\'ec\'edent.

Soit $(n_{1},n_{2})\in D(n)$, auquel est associ\'e une donn\'ee endoscopique de $G_{iso}$ et $G_{an}$ dont le groupe endoscopique est $G_{n_{1},iso}\times G_{n_{2},iso}$. Soit $(x_{1},x_{2})\in G_{n_{1},iso}(F)\times G_{n_{2},iso}(F)$ un couple $n$-r\'egulier form\'e d'\'el\'ements elliptiques.
On lui associe les donn\'ees du paragraphe 3.1: $n_{1,+}$, $\eta_{1,+}$ etc... A un \'el\'ement quelconque de la classe totale de conjugaison stable dans $G_{iso}(F)\cup G_{an}(F)$ correspondant \`a $(x_{1},x_{2})$ sont aussi associ\'ees des donn\'ees $n_{+}$, $\eta_{+}$ etc... On va calculer 
$$J^{endo}(x_{1},x_{2},f)=\sum_{x}\Delta((x_{1},x_{2}),x)J(x,f),$$
o\`u $x$ d\'ecrit cette classe totale de conjugaison stable, \`a conjugaison pr\`es. Les int\'egrales orbitales $J(x,f)$ sont calcul\'ees par la proposition pr\'ec\'edente. On en d\'eduit imm\'ediatement

(1) si les hypoth\`eses (1)(a) et (b) de 3.2 ne sont pas v\'erifi\'ees, $J^{endo}(x_{1},x_{2},f)=0$.

Supposons ces hypoth\`eses v\'erifi\'ees.
Comme on l'a expliqu\'e en 3.1, l'application $x\mapsto X$ identifie la sommation \`a la somme sur les $X\in \mathfrak{g}_{iso}(F)\cup \mathfrak{g}_{an}(F)$ dans la classe totale de conjugaison stable correspondant \`a celle de $(X_{1},X_{2})$.  Le lemme 3.1 et la proposition 3.2 expriment les termes $\Delta((x_{1},x_{2}),x)$ et $J(x,f)$  \`a l'aide de l'\'el\'ement $X$, \`a l'exception de la constante $c_{\sharp}(r',r'',w',w'')$ car l'indice $\sharp$ est celui tel que $x$ appartienne \`a $G_{\sharp}(F)$. Mais la relation 3.1(3) calcule cet indice \`a l'aide de $X$. Remarquons que, dans cette relation, on peut remplacer $val_{F}(\eta_{-})$ par $r''$ d'apr\`es l'hypoth\`ese (1)(a) de 3.2. Rappelons que $b=0$ si $r''>0$ ou si $r''=0$ et $r'$ est pair et $b=1$ si $r''<0$ ou si $r''=0$ et $r'$ est impair. D\'efinissons
$$c(r',r'',w',w'')=\left\lbrace\begin{array}{cc}1,&si\,\,0<r''\leq r'\,\,ou \,\,r''=0\,\,et \,\,r'\,\,est\,\,pair,\\ sgn_{CD}(w''),&si\,\,r'<r'',\\  (-1)^{dr''}&si\,\,-r'\leq r''<0\,\,ou\,\,r''=0\,\,et\,\,r'\,\,est\,\,impair,\\  (-1)^{dr''}sgn_{CD}(w'),&si \,\,r''<-r'.\\ \end{array}\right.$$
 On a alors l'\'egalit\'e
$$c_{\sharp}(r',r'',w',w'')=c(r',r'',w',w'') sgn^{*}(X)^b.$$
On peut aussi remplacer $(-1)^{d_{2}val_{F}(\eta_{-})}$ par $(-1)^{dr''}$ dans l'\'enonc\'e du lemme 3.1. Alors ce lemme et la proposition 3.2 entra\^{\i}nent l'\'egalit\'e
$$(2) \qquad J^{endo}(x_{1},x_{2},f)=c(r',r'')c(r',r'',w',w'')(-1)^{d_{2}r''}\sum_{{\bf d}\in D}\sum_{{\bf v}\in {\cal V}({\bf d})}$$
$$\sum_{X}
\Delta((X_{1},X_{2}),X)sgn^{*}(X)^bJ(X,f^{b}[{\bf d},{\bf v}]).$$
Fixons ${\bf d}\in D$ et ${\bf v}\in {\cal V}({\bf d})$. On a par d\'efinition une d\'ecomposition $f^{b}[{\bf d},{\bf v}]=f^b_{+}\otimes f^b_{-}\otimes \otimes_{i\in I}f^b_{i}$. Posons
$$J^{endo,*}(X_{1,+},X_{2,+},f^b_{+})=\sum_{X_{+}}\Delta_{+}((X_{1,+},X_{2,+}),X_{+})sgn^{*}(X_{+})^bJ(X_{+},f^b_{+}),$$
o\`u $X_{+}$ parcourt la classe totale de conjugaison stable correspondant \`a celle de $(X_{1,+},X_{2,+})$. On d\'efinit de m\^eme $J^{endo,*}(X_{1,-},X_{2,-},f^b_{-})$ et $J^{endo,*}(X_{1,i},X_{2,i},f^b_{i})$ pour $i\in I$.
La somme en $X$ de l'expression (2) est \'egale \`a
$$(3) \qquad J^{endo,*}(X_{1,+},X_{2,+},f^b_{+})J^{endo,*}(X_{1,-},X_{2,-},f^b_{-})\prod_{i\in I}J^{endo,*}(X_{1,i},X_{2,i},f^b_{i}).$$

Supposons d'abord $b=0$. Alors $sgn^*(X_{+})^b$ dispara\^{\i}t de la d\'efinition de $J^{endo,*}(X_{1,+},X_{2,+},f^b_{+})$ et ce terme est l'int\'egrale endoscopique $J^{endo}(X_{1,+},X_{2,+},f^0_{+})$. De m\^eme pour les autres termes de (3). On applique les (ii) des lemmes 2.1, 2.2 et 2.3: le produit de ces integrales endoscopiques est nul sauf si les conditions suivantes sont v\'erifi\'ees:

$$(4) \left\lbrace\begin{array}{c}n_{1,+}=\frac{(r'_{+}+\vert r''\vert )^2-1}{4}+N'_{+},\\n_{2,+}=\frac{(r'_{+}-\vert r''\vert )^2-1}{4}+N''_{+},\\n_{1,-}=\frac{(r'_{-}+\vert r''\vert )^2}{4}+N'_{-},\\n_{2,-}=\frac{(r'_{-}-\vert r''\vert )^2}{4}+N''_{-},\\d_{1,i}=d'(i),\,\,d_{2,i}=d''_{i}\,\,pour\,\,tout\,\,i\in I\\ \end{array}\right.$$

(5) $val_{F}(\eta_{1,-})\equiv \frac{r'_{-}+\vert r''\vert }{2}\,\, mod\,\,2{\mathbb Z}$, $val_{F}(\eta_{2,-})\equiv \frac{ r'_{-}-\vert r''\vert  }{2}\,\, mod\,\,2{\mathbb Z}$.
 
Nootre hypoth\`ese $b=0$ implique que $\vert r''\vert =r''$. Pour $j=1,2$, on a $n_{j}=n_{j,+}+n_{j,-}+\sum_{i\in I }f_{i}d_{j,i}$.  En utilisant l'\'egalit\'e $\{r'_{+},r'_{-}\}=\{r',r'+1\}$ et le fait que ${\bf d}\in D$, la condition (4)  ci-dessus implique
$$n_{1}=\frac{(r'+r'')^2+(r'+r''+1)^2-1}{4}+N',\,\, n_{2}=\frac{(r'-r'')^2+(r'-r''+1)^2-1}{4}+N''.$$
 C'est pr\'ecis\'ement le couple $(n_{1},n_{2})$ d\'efini en 1.2. Cela prouve que si notre couple $(n_{1},n_{2})$ n'est pas celui d\'efini dans ce paragraphe, (3) est nul. Ceci \'etant vrai pour tous ${\bf d}$, ${\bf v}$, on a $J^{endo}(x_{1},x_{2},f)=0$ d'apr\`es (2). Cela \'etant vrai pour tout $(x_{1},x_{2})$, le transfert de $f$ relatif \`a $(n_{1},n_{2})$ est nul.
 
Supposons maintenant $b=1$. Le couple $(n_{1,+},n_{2,+})$ d\'efinit une donn\'ee endoscopique pour les deux groupes sp\'eciaux orthogonaux impairs  dont les alg\`ebres de Lie contiennent nos \'el\'ements $X_{+}$.  Le couple $(n_{2,+},n_{1,+})$ d\'efinit aussi une telle donn\'ee.  On a donc aussi un facteur de transfert $\Delta_{+}((X_{2,+},X_{1,+}),X_{+})$. Comme on l'a dit en  \cite{W5} 2.1, on a l'\'egalit\'e
$$\Delta_{+}((X_{2,+},X_{1,+}),X_{+})=sgn^*(X_{+})\Delta_{+}((X_{1,+},X_{2,+}),X_{+}).$$
On voit alors que $J^{endo,*}(X_{1,+},X_{2,+},f^b_{+})=J^{endo}(X_{2,+},X_{1,+},f^1_{+})$. De m\^eme pour les autres termes de (3). Le raisonnement se poursuit comme ci-dessus. On permute  les r\^oles des indices $1$ et $2$; on permute  aussi $N'$ et $N''$ puisqu'on remplace la fonction $f^{0}[{\bf d},{\bf v}]$ par $f^{1}[{\bf d},{\bf v}]$; enfin, l'hypoth\`ese $b=1$ entra\^{\i}ne que $\vert r''\vert =-r''$. Ces trois modifications conduisent au m\^eme r\'esultat: le transfert de $f$ relatif \`a $(n_{1},n_{2})$ est nul si $(n_{1},n_{2})$ n'est pas le couple d\'efini en 1.2.  Cela d\'emontre le (ii) de la proposition 1.2 pour les fonctions $\varphi'=\varphi_{w'}$ et $\varphi''=\varphi_{w''}$. En faisant varier $w'$ et $w''$, cela d\'emontre cette assertion pour toutes fonctions cuspidales  $\varphi'$ et $\varphi''$. 

\bigskip

\subsection{D\'emonstration du (i) de la proposition 1.2}
On poursuit le calcul pr\'ec\'edent en supposant que $(n_{1},n_{2})$ est le couple d\'efini en 1.2. On suppose v\'erifi\'ees les hypoth\`eses (1)(a) et (1)(b) de 3.2. On fixe ${\bf d}\in D$ et ${\bf v}\in {\cal V}({\bf d})$. On suppose d'abord que $b=0$. Comme on l'a expliqu\'e, la somme en $X$ de 3.3(2) est \'egale \`a
$$(1) \qquad J^{endo}(X_{1,+},X_{2,+},f_{+}^0)J^{endo}(X_{1,-},X_{2,-},f_{-}^0)\prod_{i\in I}J^{endo}(X_{1,i},X_{2,i},f_{i}^0).$$
Ce produit est nul sauf si les conditions (4) et (5) de 3.3 sont v\'erifi\'ees. L'hypoth\`ese (5) est ind\'ependante de ${\bf d}$ et ${\bf v}$. R\'ecrivons-la (en se rappelant que $r''\geq0$ puisque $b=0$):

(2)  $val_{F}(\eta_{1,-})\equiv \frac{r'_{-}+ r'' }{2}\,\, mod\,\,2{\mathbb Z}$, $val_{F}(\eta_{2,-})\equiv \frac{ r'_{-}- r''  }{2}\,\, mod\,\,2{\mathbb Z}$.
 
Puisque $\eta_{-}=\eta_{1,-}\eta_{2,-}$, elle implique l'hypoth\`ese (1)(a) de 3.2. Si l'hypoth\`ese (4)  de 3.3 est v\'erifi\'ee, on a les in\'egalit\'es
$$(3) \qquad n_{1,+}\geq \frac{(r'_{+}+ r'')^2-1}{4},\,\, n_{2,+}\geq \frac{(r'_{+}-r'')^2-1}{4},$$
$$n_{1,-}\geq \frac{(r'_{-}+r'')^2}{4},\,\, n_{2,-}\geq \frac{(r'_{-}-r'')^2}{4}.$$
Puisque $n_{+}=n_{1,+}+n_{2,+}$, $n_{-}=n_{1,-}+n_{2,-}$, ces in\'egalit\'es impliquent l'hypoth\`ese (1)(b) de 3.2. On peut donc oublier les hypoth\`eses (1)(a) et (b) de 3.2 et supposer v\'erifi\'ees les hypoth\`eses (2) et (3)  ci-dessus. Alors les relations (4) de 3.3 d\'eterminent un unique \'el\'ement ${\bf d}$ dont on v\'erifie qu'il appartient bien  \`a $D$. On suppose d\'esormais que ${\bf d}$ est cet unique \'el\'ement.

  L'int\'egrale endoscopique  $J^{endo}(X_{1,+},X_{2,+},f_{+}^0)$ est calcul\'ee par le lemme 2.1. Adaptons les notations. On note $t'_{1,+},t''_{1,+},t'_{2,+}$ et $t''_{2,+}$ les termes not\'es $t'_{1}$ etc... en 2.1 associ\'es \`a $n_{+}$, $r'_{+}$ et $\vert r''\vert $. On note $C_{+}({\bf v})$ la constante $C$ de 2.2. On pose $f_{1,+}={\cal Q}(t'_{1,+},t''_{1,+})^{Lie}\circ\rho_{N'_{+}}\circ\iota_{N'_{+},0}(\varphi_{v'_{+}})$ et $f_{2,+}={\cal Q}(t'_{2,+},t''_{2,+})^{Lie}\circ\rho_{N''_{+}}\circ\iota_{N''_{+},0}(\varphi_{v''_{+}})$. Alors
$$J^{endo}(X_{1,+},X_{2,+},f_{+}^0)=C_{+}({\bf v})S(X_{1,+},f_{1,+})S(X_{2,+},f_{2,+}).$$
En adaptant de fa\c{c}on similaire les notations et d\'efinitions,  le lemme 2.2 fournit l'\'egalit\'e
$$J^{endo}(X_{1,-},X_{2,-},f_{-}^0)=C_{-}({\bf v})S(X_{1,-},f_{1,-})S(X_{2,-},f_{2,-}),$$
tandis que le lemme 2.3 fournit l'\'egalit\'e
$$J^{endo}(X_{1,i},X_{2,i},f_{i}^0)=S(X_{1,i},f_{1,i})S(X_{2,i},f_{2,i}).$$
 Posons $f_{1}[{\bf v}]=f_{1,+}\otimes f_{1,-}\otimes\otimes_{i\in I}f_{1,i}$ et 
 $$S(X_{1},f_{1}[{\bf v}])=S(X_{1,+},f_{1,+})S(X_{1,-},f_{1,-})\prod_{i\in I}S(X_{1,i},f_{1,i}).$$
 D\'efinissons de m\^eme $f_{2}[{\bf v}]$ et $S(X_{2},f_{2}[{\bf v}])$. Alors l'expression (1) ci-dessus vaut 
 $$C_{+}({\bf v})C_{-}({\bf v})S(X_{1},f_{1}[{\bf v}])S(X_{2},f_{2}[{\bf v}]).$$
  D\'efinissons une nouvelle constante $C(r',r'',w',w'')$ par les \'egalit\'es

  si $r''\leq r'$, 
 $$C(r',r'',w',w'')=(-1)^{d_{2}r''}sgn(-1)^{\frac{r'_{+}-r''-1}{2}};$$
 
  si $r'<r''$, 
 $$C(r',r'',w',w'')=(-1)^{d_{2}r''}sgn(-1)^{ r'+r''+1}sgn_{CD}(w'').$$
 
 Montrons que
 
 (4)  l'expression (1) vaut
 $$C(r',r'',w',w'')S(X_{1},f_{1}[{\bf v}])S(X_{2},f_{2}[{\bf v}]).$$
 
 Supposons d'abord $r''\leq r'$. Alors $r''\leq r'_{+}$ et $r''\leq r'_{-}$. En utilisant les d\'efinitions de 2.1 et 2.2, on  obtient 
 $$C_{+}({\bf v})C_{-}({\bf v})=sgn(-1)^{\frac{r'_{+}+r''-1}{2}+val_{F}(\eta_{+})}sgn(\eta_{2,+}\varpi^{-val_{F}(\eta_{2,+})})^{val_{F}(\eta_{+})}$$
 $$sgn(\eta_{2,-}\varpi^{-val_{F}(\eta_{2,-})})^{val_{F}(\eta_{-})}.$$
 D'apr\`es 3.1(2) et l'hypoth\`ese (1)(a) de 3.2 (qui est v\'erifi\'ee), $val_{F}(\eta_{+})=-val_{F}(\eta_{-})\equiv r''\,\,mod\,\,2{\mathbb Z}$. D'apr\`es la m\^eme relation 3.1(2) appliqu\'ee aux donn\'ees index\'ees par $2$, on a $val_{F}(\eta_{2,+})+val_{F}(\eta_{2,-})=0$ et $sgn(\eta_{2,+}\eta_{2,-})=(-1)^{d_{2}}$. On en d\'eduit l'\'egalit\'e 
  $C_{+}({\bf v})C_{-}({\bf v})=C(r',r'',w',w'')$. Supposons maintenant $r'\leq r''-2$. Alors $r'_{+}< r''$ et $r'_{-}<r''$. On voit que $C_{+}({\bf v})C_{-}({\bf v})$ est \'egal \`a 
$$sgn(-1)^{\frac{r'_{+}+r''+1}{2}+val_{F}(\eta_{+})+val_{F}(\eta_{2,-})}$$
 $$sgn(\eta_{2,+}\varpi^{-val_{F}(\eta_{2,+})})^{1++val_{F}(\eta_{+})}sgn(\eta_{2,-}\varpi^{-val_{F}(\eta_{2,-})})^{1+val_{F}(\eta_{-})}  sgn_{CD}(v''_{+})sgn_{CD}(v''_{-}).$$
 Comme ci-dessus, le produit des deuxi\`eme et troisi\`eme termes vaut $(-1)^{ d_{2}(1+r'')}$. On a 
 $$val_{F}(\eta_{+})+val_{F}(\eta_{2,-})=-val_{F}(\eta_{-})+val_{F}(\eta_{2,-})=-val_{F}(\eta_{1,-})\equiv \frac{r'_{-}+r''}{2}\,\,mod\,\,2{\mathbb Z}$$
 d'apr\`es (2) ci-dessus.  Puisque $r'_{+}+r'_{-}=2r'+1$, on voit que  la puissance de $sgn(-1)$ dans l'expression ci-dessus co\"{\i}ncide avec celle figurant dans $C(r',r'',w',w'')$.  D'apr\`es la d\'efinition de ${\bf d}$ et  le fait que ${\bf v}\in {\cal V}({\bf d})$, on a l'\'egalit\'e $sgn_{CD}(v''_{+})sgn_{CD}(v''_{-})=sgn_{CD}(w'')sgn(\eta_{2,+}\eta_{2,-})$. Ce dernier facteur est \'egal \`a $(-1)^{d_{2}}$ ainsi qu'on l'a d\'ej\`a dit. En rassemblant ces calculs, on obtient de nouveau l'\'egalit\'e 
  $C_{+}({\bf v})C_{-}({\bf v})=C(r',r'',w',w'')$. 
   Supposons maintenant $r'=r''-1$. Alors $r'_{+}=r'<r''$ mais $r'_{-}=r'+1=r''$.  Comme dans le cas $r'\leq r''-2$, on obtiendrait l'\'egalit\'e voulue si $C_{-}({\bf v})$ \'etait d\'efinie par les formules du cas $r'_{-}< r''$ et non pas de notre cas $r'_{-}\geq r''$. Le rapport entre les deux formules est  
    $$(5) \qquad sgn(-1)^{val_{F}(\eta_{2,-})}sgn(\eta_{2,-}\varpi^{-val(\eta_{2,-})})^{val_{F}(\eta_{-})}sgn_{CD}(v''_{-}).$$
Si cette expression vaut $1$, on a fini. Supposons qu'elle vaille $-1$.   Puisque $r'_{-}=r''$, l'hypoth\`ese  (2) plus haut  implique que $val_{F}(\eta_{2,-})$ est paire, ce qui fait dispara\^{\i}tre le premier terme de notre expression . On a aussi $r'_{2,-}=r''_{2,-}=\frac{\vert r'_{-}-r''\vert }{2}=0$. La remarque de 2.2   entra\^{\i}ne que,  sous notre hypoth\`ese que (5) vaut $-1$, $f_{2,-}$ est nulle, donc aussi $f_{2}[{\bf v}]$.  Mais alors $S(X_{2},f_{2}[{\bf v}])=0$ et la constante n'a plus d'importance. Cela prouve (4).

  L'\'egalit\'e 3.3(2) devient
  $$(6) \qquad J^{endo}(x_{1},x_{2},f)=c(r',r'')c(r',r'',w',w'')C(r',r'',w',w'')(-1)^{d_{2}r''}$$
  $$\sum_{{\bf v}\in {\cal V}({\bf d})}S(X_{1},f_{1}[{\bf v}])S(X_{2},f_{2}[{\bf v}]).$$

Comme en 1.2, on d\'efinit  des entiers $r'_{1},r''_{1},r'_{2},r''_{2}$ et les fonctions $f_{1}=\Psi\circ k\circ\rho\iota(\varphi_{w'})$ et $f_{2}=\Psi\circ k\circ\rho\iota(\varphi_{w''})$.  Les int\'egrales orbitales stables $S(x_{1},f_{1})$  et $S(x_{2},f_{2})$ sont des cas particuliers d'int\'egrales endoscopiques $J^{endo}(x_{1},x_{2},f)$. Elles sont donc  nulles si les analogues des hypoth\`eses (2) et (3)  ne sont pas v\'erifi\'ees. Sinon, elles sont donn\'ees par des formules similaires \`a (6). Dans toutes ces formules, les termes index\'es par l'indice $2$ disparaissent ainsi que ceux faisant intervenir $N''$ et $w''$. D'autre part, l'analogue de l'entier $b\in \{0,1\}$ est toujours $0$. En effet, d'apr\`es les d\'efinitions, on a $r''_{1},r''_{2}\geq0$ et si,  par exemple, $r''_{1}=0$, on a aussi $r'_{1}=0$ donc $r'_{1}$ pair. On voit que l'analogue de  (2) pour $S(x_{1},f_{1})$ est $val_{F}(\eta_{1,-})\equiv\frac{r'_{1,-}+r''_{1}}{2}\,\,mod\,\,2{\mathbb Z}$, tandis que l'analogue pour $S(x_{2},f_{2})$ est $val_{F}(\eta_{2,-})\equiv\frac{r'_{2,-}+r''_{2}}{2}\,\,mod\,\,2{\mathbb Z}$. D'apr\`es la relation (3) du paragraphe 4 ci-apr\`es,  $r'_{1,-}+r''_{1}=r'_{-}+r''$, $r'_{2,-}+r''_{2}=\vert r'_{-}-r''\vert $. Il en r\'esulte que la conjonction des deux analogues de  (2) est \'equivalente \`a cette relation  (2) elle-m\^eme. 
L'analogue de (3) ci-dessus pour $S(x_{1},f_{1})$ est
$$n_{1,+}\geq \frac{(r'_{1,+}+r''_{1})^2-1}{4},\,\, n_{1,-}\geq\frac{(r'_{1,-}+r''_{1})^2}{4},$$
tandis que l'analogue pour $S(x_{2},f_{2})$ est
$$n_{2,+}\geq \frac{(r'_{2,+}+r''_{2})^2-1}{4},\,\, n_{1,-}\geq\frac{(r'_{2,-}+r''_{2})^2}{4}.$$
De nouveau, le calcul montre que la conjonction de ces conditions est \'equivalente \`a la relation (3) elle-m\^eme. Cela d\'emontre que, si (2)  et (3)  ne sont pas v\'erifi\'ees, $S(x_{1},f_{1})S(x_{2},f_{2})=0$. Puisqu'on a d\'ej\`a vu que $J^{endo}(x_{1},x_{2},f)=0$, on obtient $J^{endo}(x_{1},x_{2},f)=S(x_{1},f_{1})S(x_{2},f_{2})$ dans ce cas.

On suppose maintenant v\'erifi\'ees nos hypoth\`eses (2) et (3).   De m\^eme que l'on a d\'etermin\'e un unique ${\bf d}\in D$, on d\'etermine un unique ${\bf d}_{1}$ relatif \`a $S(x_{1},f_{1})$ et un unique ${\bf d}_{2}$ relatif \`a $S(x_{2},f_{2})$. Ecrivons ${\bf d}=(N'_{+},N'_{-},N''_{+},N''_{-},(d'_{i},d''_{i})_{i\in I})$. Tous ces entiers sont d\'etermin\'es par la relation 3.3(4). De m\^eme, ${\bf d}_{1}$ et ${\bf d}_{2}$ sont d\'etermin\'es par les analogues de cette relation. Par les m\^emes calculs que ci-dessus, on voit que
${\bf d}_{1}=(N'_{+},N'_{-},0,0, (d'_{i},0)_{i\in I})$ tandis que ${\bf d}_{2}=(N''_{+},N''_{-},0,0,
(d''_{i},0)_{i\in I})$. On en d\'eduit que ${\cal V}'({\bf d})\simeq {\cal V}({\bf d}_{1})$ et ${\cal V}''({\bf d})\simeq {\cal V}({\bf d}_{2})$, d'o\`u ${\cal V}({\bf d})\simeq {\cal V}({\bf d}_{1})\times {\cal V}({\bf d}_{2})$. On identifie ces deux ensembles.  A l'aide de $r',r''$, ${\bf d}$  et d'un \'el\'ement ${\bf v}\in {\cal V}({\bf d})$, on a d\'efini une fonction $f_{1}[{\bf v}]$. Pour $j=1,2$, \`a l'aide de $r'_{j},r''_{j}$, ${\bf d}_{j}$  et d'un \'el\'ement ${\bf v}_{j}\in {\cal V}({\bf d}_{j})$, on d\'efinit de m\^eme une fonction $f_{j}[{\bf v}_{j}]$.  Les analogues de l'\'egalit\'e (6) sont
$$(7)_{1}\qquad S(x_{1},f_{1})=c(r'_{1},r''_{1})c(r'_{1},r''_{1},w')C(r'_{1},r''_{1},w') \sum_{{\bf v}_{1}\in {\cal V}({\bf d}_{1})}S(X_{1},f_{1}[{\bf v}_{1}]) ;$$
$$(7)_{2}\qquad S(x_{2},f_{2})=c(r'_{2},r''_{2})c(r'_{2},r''_{2},w'')C(r'_{2},r''_{2},w'') \sum_{{\bf v}_{2}\in {\cal V}({\bf d}_{2})}S(X_{2},f_{2}[{\bf v}_{2}]),$$
o\`u on a adapt\'e de fa\c{c}on \'evidente la notation des constantes.

Soit ${\bf v}=({\bf v}_{1},{\bf v}_{2})\in {\cal V}({\bf d})$. 
Montrons que

(8)  $f_{1}[{\bf d}]=f_{1}[{\bf d}_{1}]$, $f_{2}[{\bf d}]=f_{2}[{\bf d}_{2}]$. 

On traite le cas de l'indice $1$. D'apr\`es la  d\'efinition donn\'ee plus haut, $f_{1}[{\bf v}]=f_{1,+}\otimes f_{1,-}\otimes \otimes_{i\in I}f_{1,i}$, avec par exemple $f_{1,+}={\cal Q}(t'_{1,+},t''_{1,+})^{Lie}\circ\rho_{N'_{+}}\circ\iota_{N'_{+},0}(\varphi_{v'_{+}})$. Quand on remplace $r',r''$, ${\bf d}$  et ${\bf v}$ par $r'_{1},r''_{1}$, ${\bf d}_{1}$  et ${\bf v}_{1}$, les termes $N'_{+}$ et $v'_{+}$ ne changent pas. Les entiers $t'_{1,+}$ et $t''_{1,+}$ sont d\'efinis par $\{t'_{1,+},t''_{1,+}\}=\{\frac{r'_{+}+r''+1}{2},\frac{r'_{+}+r''-1}{2}\}$ et $t'_{1,+}\equiv 1+val_{F}(\eta_{1,+})\,\,mod\,\,2{\mathbb Z}$. La relation (3) du paragraphe 4 ci-apr\`es montre que l'ensemble $ \{\frac{r'_{+}+r''+1}{2},\frac{r'_{+}+r''-1}{2}\}$ ne change par quand on remplace $r',r''$ par $r'_{1},r''_{1}$. La congruence exig\'ee non plus. Donc $t'_{1,+}$ et $t''_{1,+}$ ne changent pas et  $f_{1,+}$ non plus. 
  Un calcul analogue s'applique aux composantes $f_{1,-}$ et $f_{1,i}$ pour $i\in I$. Cela prouve (8).

De (6), $(7)_{1}$ et $(7)_{2}$ se d\'eduit l'\'egalit\'e
$$J^{endo}(x_{1},x_{2},f)=CS(x_{1},f_{1})S(x_{2},f_{2}),$$
o\`u
$$C=c(r',r'')c(r',r'',w',w'')C(r',r'',w',w'')(-1)^{d_{2}r''}c(r'_{1},r''_{1})c(r'_{1},r''_{1},w')C(r'_{1},r''_{1},w') $$
$$c(r'_{2},r''_{2})c(r'_{2},r''_{2},w'')C(r'_{2},r''_{2},w'').$$
Remarquons que l'on n'a pas besoin d'inverser les coefficients de $(7)_{1}$ et $(7)_{2}$: ce sont tous des signes. On a l'\'egalit\'e

(9) $C=1$.

En reprenant les d\'efinitions des diverses constantes, on obtient l'\'egalit\'e
$$c(r',r'')c(r',r'',w',w'')C(r',r'',w',w'')(-1)^{d_{2}r''}=(-1)^{n}U,$$
o\`u $U$ est d\'efini au paragraphe 4 ci-apr\`es. Les autres termes sont les analogues relatifs aux donn\'ees index\'ees par $1$ et $2$. Evidemment, $n=n_{1}+n_{2}$ et l'\'egalit\'e $C=1$ r\'esulte de l'\'egalit\'e $U=U_{1}U_{2}$, cf. 4(4) ci-dessous.

A l'aide de (9), on obtient $J^{endo}(x_{1},x_{2},f)=S(x_{1},f_{1})S(x_{2},f_{2})$. Cette \'egalit\'e est donc v\'erifi\'ee  avec ou sans les hypoth\`eses (2) et (3) et elle l'est pour tous $x_{1},x_{2}$. Donc $f_{1}\otimes f_{2}$ est le transfert de $f$.  Cela d\'emontre le (i) de la proposition 1.2 (sous l'hypoth\`ese $b=0$) pour les fonctions $\varphi'=\varphi_{w'}$ et $\varphi''=\varphi_{w''}$. Comme dans le paragraphe pr\'ec\'edent, cela entra\^{\i}ne  la m\^eme assertion pour toutes fonctions cuspidales $\varphi'$ et $\varphi''$. 

On a suppos\'e $b=0$. Supposons maintenant $b=1$. Comme on l'a dit en 3.3,  la somme en $X$ de 3.3(2) est  	alors \'egale \`a
$$(1) \qquad J^{endo}(X_{2,+},X_{1,+},f_{+}^1)J^{endo}(X_{2,-},X_{1,-},f_{-}^1)\prod_{i\in I}J^{endo}(X_{2,i},X_{1,i},f_{i}^1).$$
Le calcul se poursuit comme ci-dessus en permutant les r\^oles des indices $1$ et $2$, en permutant $N'$ et $N''$ et en tenant compte de l'\'egalit\'e $\vert r''\vert =-r''$. Le r\'esultat est le m\^eme, on laisse les d\'etails au lecteur. Cela ach\`eve de prouver la proposition 1.2. $\square$

\bigskip

\section{Annexe}

On rassemble ici quelques  calculs \'el\'ementaires utilis\'es dans l'article. Soient $r'\in {\mathbb N}$ et $r''\in {\mathbb Z}$.  On pose

  $r'_{1}=sup([\frac{r'+r''}{2}],-[\frac{r'+r''}{2}]-1)$, $r''_{1}=\vert [\frac{r'+r''+1}{2}]\vert $, $r'_{2}=sup([\frac{r'-r''}{2}],-[\frac{r'-r''}{2}]-1)$, $r''_{2}=\vert [\frac{r'-r''+1}{2}]\vert $.
  
  Remarquons que, si l'on change $r''$ en $-r''$, on \'echange les couples $(r'_{1},r''_{1})$ et $(r'_{2},r''_{2})$.  

  (1) On a $r''_{1}+r''_{2}\equiv r''\,\,mod\,\,2{\mathbb Z}$.
  
    Preuve. Puisque $m\equiv -m\,\,mod\,\,2{\mathbb Z}$ pour tout $m\in {\mathbb Z}$, on peut remplacer  $r''_{1}$ par  $[\frac{r'+r''+1}{2}] $ et $r''_{2}$ par $[\frac{r'-r''+1}{2}] $. On a $[\frac{r'-r''+1}{2}]=[\frac{r'+r''+1}{2}-r'']=[\frac{r'+r''+1}{2}]-r''$, d'o\`u (1). $\square$
    
    Pour faciliter les calculs qui suivent, on pose $a=0$ si $r'+r''$ est pair et $a=1$ si $r'+r''$ est impair, $A=0$ si $\vert r''\vert \leq r'$ et $A=1$ si $r'<\vert r''\vert $.
    
    (2) On a les \'egalit\'es
    $$r^{'2}_{1}+r'_{1}+r^{''2}_{1}=\frac{(r'+r'')^2+(r'+r''+1)^2-1}{4},$$
     $$r^{'2}_{2}+r'_{2}+r^{''2}_{2}=\frac{(r'-r'')^2+(r'-r''+1)^2-1}{4}.$$
     
     Preuve. Pour $m\in {\mathbb Z}$, $m^2+m$ est invariant par la transformation $m\mapsto -m-1$ et $m^2$ est invariant par $m\mapsto -m$. On peut donc remplacer $r'_{1}$ par $x'=[\frac{r'+r''}{2}]$ et $r''_{1}$ par $x''=[\frac{r'+r''+1}{2}]$. On  a $x'=\frac{r'+r''-a}{2}$, $x''=\frac{r'+r''+a}{2}$. Un calcul alg\'ebrique  montre que
     $$x^{'2}+x'+x^{''2}=\frac{(r'+r'')^2+(r'+r''+1)^2-1}{4}+ \frac{a^2-a}{2}.$$
     D'o\`u la premi\`ere \'egalit\'e puisque $a^2=a$. La seconde \'egalit\'e en r\'esulte, en changeant $r''$ en $-r''$. $\square$

On d\'efinit $r'_{+}$ et $r'_{-}$ par

 si $r'\equiv r''\,\,mod\,\,2{\mathbb Z}$, $r'_{+}=r'+1$, $r'_{-}=r'$;
 
  si $r'\not\equiv r''\,\,mod\,\,2{\mathbb Z}$, $r'_{+}=r'$, $r'_{-}=r'+1$. 
  
  Autrement dit, $r'_{+}=r'+1-a$ et $r'_{-}=r'+a$. En rempla\c{c}ant le couple $(r',r'')$ par $(r'_{1},r''_{1})$ ou $(r'_{2},r''_{2})$, on d\'efinit de m\^eme $r'_{1,+}$ et $r'_{1,-}$ ou $r'_{2,+}$ et $r'_{2,-}$.

  (3) On a les \'egalit\'es
  
  $r'_{1,+}+r''_{1}=\vert r'_{+}+r''\vert $, $r'_{1,-}+r''_{1}=\vert r'_{-}+r''\vert $, $r'_{2,+}+r''_{2}=\vert r'_{+}-r''\vert $, $r'_{2,-}+r''_{2}=\vert r'_{-}-r''\vert $.
  
  Preuve. De m\^eme que l'on a d\'efini $a$ et $A$ pour $(r',r'')$, on d\'efinit $a_{1}$ et $A_{1}$ pour $(r'_{1},r''_{1})$ et $a_{2}$ et $A_{2}$ pour $(r'_{2},r''_{2})$.  Supposons $r''\geq0$. Alors $r'_{1}=[\frac{r'+r''}{2}]=\frac{r'+r''-a}{2}$ et $r''_{1}=[\frac{r'+r''+1}{2}]=\frac{r'+r''+a}{2}$. Si $a=0$, $r'_{1}=r''_{1}$ donc $a_{1}=A_{1}=0$. Si $a=1$, $r'_{1}=r''_{1}-1$ donc $a_{1}=A_{1}=1$.  On a donc $a_{1}=A_{1}=a$ quel que soit $a$. D'o\`u $r'_{1,+}=r'_{1}+1-a_{1}=\frac{r'+r''-a}{2}+1-a$ et $r'_{1,-}=r'_{1}+a_{1}=\frac{r'+r''-a}{2}+a$.  On calcule alors
  $$r'_{1,+}+r''_{1}= \frac{r'+r''-a}{2}+1-a+\frac{r'+r''+a}{2}=r'+r''+1-a=r'_{+}+r''=\vert r'_{+}+r''\vert ,$$
 $$r'_{1,-}+r''_{1}= \frac{r'+r''-a}{2}+a+\frac{r'+r''+a}{2}=r'+r''+a=r'_{-}+r''=\vert r'_{-}+r''\vert.$$
 Ce sont les deux premi\`eres \'egalit\'es de (3). Supposons $b=0$. Alors $r'_{2}=[\frac{r'-r''}{2}]=\frac{r'-r''-a}{2}$ et $r''_{2}=[\frac{r'-r''+1}{2}]=\frac{r'-r''+a}{2}$. Le calcul est le m\^eme que pr\'ec\'edemment, on a simplement chang\'e $r''$ en $-r''$ dans les formules (en particulier, on a $a_{2}=A_{2}=a$ si $b=0$). Supposons $b=1$. Alors $r'_{2}=-1-[\frac{r'-r''}{2}]=-1+\frac{r''-r'+a}{2}$ et $r''_{2}=-[\frac{r'-r''+1}{2}]=\frac{r''-r'-a}{2}$. On constate que cette fois, $a_{2}=A_{2}=1-a$. Alors $r'_{2,+}=r'_{2}+1-a_{2}=r'_{2}+a=\frac{r''-r'+a}{2}+a-1$ et $r'_{2,-}=r'_{2}+a_{2}=r'_{2}+1-a=\frac{r''-r'+a}{2}-a$. On calcule alors
 $$r'_{2,+}+r''_{2}= \frac{r''-r'+a}{2}+a-1+\frac{r''-r'-a}{2}=r''-r'+a-1=r''-r'_{+}=\vert r'_{+}-r''\vert ,$$
 $$r'_{2,-}+r''_{2}= \frac{r''-r'+a}{2}-a+\frac{r''-r'-a}{2}=r''-r'-a=r''-r'_{-}=\vert r'_{-}-r''\vert.$$ 
 Ce sont encore les deux derni\`eres \'egalit\'es de (3). On a suppos\'e $r''\geq0$. Si $r''\leq 0$, on a remarqu\'e que nos couples $(r'_{1},r''_{1})$ et $(r'_{2},r''_{2})$ sont les m\^emes que ceux d\'eduits du couple $(r',-r'')$, \`a ceci pr\`es que l'on doit \'echanger les indices $1$ et $2$. Les \'egalit\'es (3)  pour $(r',r'')$ se d\'eduisent par ce changement  des m\^emes \'egalit\'es que l'on vient de d\'emontrer pour le couple $(r',-r'')$. $\square$
 
 D\'efinissons un  nombre $u$ par
 
 si $\vert r''\vert \leq r'$, $u= \frac{r^{'2}-r'}{2}+\frac{r^{''2}-\vert r''\vert }{2}+\frac{r'_{+}-\vert r''\vert -1}{2}$,
 
 si $r'< \vert r''\vert $, $u= \frac{r^{'2}-r'}{2}+\frac{r^{''2}-\vert r''\vert }{2}+r'+r''+1$.
 
 Posons $U=(-1)^{r''}sgn(-1)^{u}$. On d\'efinit de m\^eme $u_{1}$, $U_{1}$ et $u_{2}$ et $U_{2}$. 
 
 (4) On a l'\'egalit\'e $U=U_{1}U_{2}$. 
 
 Preuve. Comme dans la preuve pr\'ec\'edente, on peut supposer $r''\geq0$. On a $(-1)^{r''}=(-1)^{r''_{1}+r''_{2}}$ d'apr\`es (1) et il suffit de d\'emontrer que $u\equiv u_{1}+u_{2}\,\,mod\,\,2{\mathbb Z}$. Consid\'erons l'expression qui d\'efinit $u$ dans le cas $r''\leq r'$. On y remplace $r'_{+}$ par $r'+1-a$. On peut ajouter $r'+r''+a$ qui est pair et on calcule
 $u\equiv \frac{r^{'2}+r^{''2}+r'}{2}+\frac{r'+a}{2}\,\,mod \,\, 2{\mathbb Z}$. La diff\'erence entre l'expression de $u$ dans le cas $r'< r''$ et cette expression dans le cas $r''\leq r'$ est 
 $$r'+r''+1-\frac{r'_{+}- r'' -1}{2}\equiv r'-r''+1-\frac{r'-a-r''}{2}\equiv  \frac{r'+a-r''}{2}+1\,\,mod\,\,2{\mathbb Z}.$$
 On obtient  alors la congruence
 $$(5) \qquad u\equiv \frac{r^{'2}+r^{''2}+r'}{2}+\frac{r'+a}{2}+A(\frac{r'+a-r''}{2}+1)\,\,mod\,\,2{\mathbb Z}$$
 en tout cas. La demi-somme des membres de droite de (2) vaut $\frac{r^{'2}+r^{''2}+r'}{2}$. En utilisant (2), on obtient
 $$u\equiv \frac{r^{'2}_{1}+r'_{1}+r^{''2}_{1}+r^{'2}_{2}+r'_{2}+r^{''2}_{2}}{2}+\frac{r'+a}{2}+A(\frac{r'+a-r''}{2}+1)\,\,mod\,\,2{\mathbb Z},$$
 puis, en utilisant les analogues de (5) pour $u_{1}$ et $u_{2}$,
 $$u\equiv u_{1}+u_{2}+\frac{r'+a-r'_{1}-a_{1}-r'_{2}-a_{2}}{2}+A(\frac{r'+a-r''}{2}+1)$$
 $$-A_{1}(\frac{r'_{1}+a_{1}-r''_{1}}{2}+1)-A_{2}(\frac{r'_{2}+a_{2}-r''_{2}}{2}+1)\,\,mod\,\,2{\mathbb Z}.$$
 Cette expression se simplifie: puisque $r'_{1}=r''_{1}$ ou $r''_{1}-1$, on a toujours $r'_{1,-}=r''_{1}$. De m\^eme, $r'_{2,-}=r''_{2}$. On a aussi remarqu\'e que $a_{1}=A_{1}=a$. D'o\`u
 $$u\equiv u_{1}+u_{2}+\frac{r'-r'_{1}-r'_{2}-a_{2}}{2}+A(\frac{r'+a-r''}{2}+1)- a-A_{2} \,\,mod\,\,2{\mathbb Z}.$$
 Si $A=0$, on a aussi $a_{2}=A_{2}=a$, $r'_{1}=\frac{r'+r''-a}{2}$, $r''_{1}=\frac{r'-r''-a}{2}$.    Si $A=1$,  on a $a_{2}=A_{2}=1-a$, $r'_{1}=\frac{r'+r''-a}{2}$, $r''_{1}=-1-\frac{r'-r''-a}{2}$. On calcule l'expression ci-dessus dans les deux cas et on conclut $u\equiv u_{1}+u_{2}\,\,mod\,\,2{\mathbb Z}$. Cela prouve (4). $\square$ 
 
 \bigskip
 
 {\bf Index des notations}
 
 $C_{cusp}^{\infty}(G_{\sharp}(F))$ 1.2;  ${\mathbb C}[\hat{\mathfrak{S}}_{m}]_{U-cusp}$ 2.3; ${\cal E}$ 2.4; $f_{\pi}$ 1.1; $\varphi_{w}$ 2.1; $\hat{i}_{\sharp}[a,e,u]$ 2.4; $\Psi$ 1.1; ${\cal Q}(r',r'')^{Lie}$ 2.1, 2.2, 2.3; $sgn$ 2.1; $\Theta_{\pi}$ 1.1; ${\cal U}$ 2.4; $X^{\zeta}(a,e,u)$ 2.4.
 
 \bigskip

  {\bf Index des notations de \cite{W5}}
 
${\mathbb C}[X]$ 1.4; $C'_{n'}$ 1.5; $C^{''\pm}_{n'',\sharp}$ 1.5; $C''_{n''}$ 1.5; $C^{GL(m)}$ 1.5; ${\mathbb C}[\hat{W}_{N}]_{cusp}$ 1.8; $D(n)$ 1.2; $D_{iso}(n)$ 1.2;  $D_{an}(n)$ 1.2; $D$ 1.7; $D^{par}$ 1.7; $\eta(Q)$ 1.1; $\eta^+(Q)$ 1.1; $\eta^-(Q)$ 1.1; $Ell_{unip}$ 1.4; $\mathfrak{Ell}_{unip}$ 1.4; $\mathfrak{Endo}_{tunip}$ 2.1; $\mathfrak{Endo}_{unip-quad}$ 2.2; $\mathfrak{Endo}_{unip-quad}^{red}$ 2.2; $\mathfrak{Endo}_{unip,disc}$ 2.4;  ${\cal F}^L$ 1.9; ${\cal F}^{par}$ 1.9; ${\cal F}$ 2.3: $\mathfrak{F}^{par}$ 2.3; $G_{iso}$ 1.1; $G_{an}$ 1.1; $\Gamma$ 1.8; $\boldsymbol{\Gamma}$ 1.8; $\tilde{GL}(2n)$ 2.1; $Irr_{tunip}$ 1.3; $\mathfrak{Irr}_{tunip}$ 1.3; $Irr_{unip-quad}$1.3; $\mathfrak{Irr}_{unip-quad}$ 1.3; $Jord(\lambda)$ 1.3; $Jord_{bp}(\lambda)$ 1.3; $Jord_{bp}^{k}(\lambda)$ 1.4; $K_{n',n''}^{\pm}$ 1.2; $k$ 1.9; $L^*$ 1.1; $L_{n',n''}$ 1.2; $l(\lambda)$ 1.3; $mult_{\lambda}$ 1.3; $\mathfrak{o}$ 1.1; $O^+(Q)$ 1.1; $O^-(Q)$ 1.1; $\varpi$ 1.1; $\pi_{n',n''}$ 1.3; ${\cal P}(N)$ 1.3; ${\cal P}^{symp}(2N)$ 1.3; $\boldsymbol{{\cal P}^{symp}}(2N)$ 1.3; $\pi(\lambda,s,\epsilon)$ 1.3; $\pi(\lambda^+,\epsilon^+,\lambda^-,\epsilon^-)$ 1.3; $\pi_{ell}(\lambda^+,\epsilon^+,\lambda^-,\epsilon^-)$ 1.4; $proj_{cusp}$ 1.5; ${\cal P}(\leq n)$ 1.5; ${\cal P}_{k}(N)$ 1.8; $\Pi(\lambda,s,h)$ 2.1; $\Pi^{st}(\lambda^+,\lambda^-)$ 2.4; ${\cal P}^{symp,disc}(2n)$ 2.4; $Q_{iso}$ 1.1; $Q_{an}$ 1.1; $\rho_{\lambda}$ 1.3;  ${\cal R}^{par}$ 1.5; ${\cal R}^{par,glob}$ 1.5; ${\cal R}^{par}_{cusp}$ 1.5; ${\cal R}^{par,glob}_{{\bf m}}$ 1.5; ${\cal R}^{par}_{{\bf m},cusp}$ 1.5; $res'_{m}$ 1.5; $res''_{m}$ 1.5; $res_{m}$ 1.5 et 1.8; $res_{{\bf m}}$ 1.5; ${\cal R}$ 1.8; ${\cal R}(\gamma)$ 1.8; ${\cal R}(\boldsymbol{\gamma})$ 1.8; ${\cal R}^{glob}$ 1.8; ${\cal R}_{cusp}$ 1.8; $Rep$ 1.9; $\rho\iota$ 1.10; $S(\lambda)$ 1.3; $\mathfrak{S}_{N}$ 1.8; $\hat{\mathfrak{S}}_{N}$ 1.8; $sgn$ 1.8; $sgn_{CD}$ 1.8; ${\cal S}_{n}$ 1.11; $\mathfrak{St}_{tunip}$ 2.1; $\mathfrak{St}_{unip-quad}$ 2.4; $\mathfrak{St}_{unip,disc}$ 2.4;  $sgn_{iso}$ 2.6; $sgn_{an}$ 2.6; $val_{F}$ 1.1; $V_{iso}$ 1.1; $V_{an}$ 1.1; $W_{N}$ 1.8; $\hat{W}_{N}$ 1.8; $w_{\alpha}$ 1.8; $w_{\alpha,\beta}$ 1.8; $w_{\alpha,\beta',\beta''}$ 1.8; $Z(\lambda)$ 1.3; $Z(\lambda,s)$ 1.3; ${\bf Z}(\lambda,s)$ 1.3; ${\bf Z}(\lambda,s)^{\vee}$ 1.3; $\vert .\vert _{F}$ 1.1.

CNRS IMJ-PRG

4 place Jussieu

75005 Paris

 jean-loup.waldspurger@imj-prg.fr

\end{document}